\documentclass[11pt]{article}

\usepackage{titlesec}
\usepackage{lipsum} 

\usepackage[T1]{fontenc}
\usepackage[latin9]{inputenc}
\usepackage{float}
\usepackage{placeins}
\usepackage{textcomp}
\usepackage{amsmath}
\usepackage{amsthm}
\usepackage{thm-restate}  
\usepackage{amssymb}
\usepackage{graphicx}
\usepackage{wasysym}
\usepackage{babel}
\usepackage[margin=1in]{geometry}
\usepackage{setspace}
\usepackage{multirow}
\makeatletter
\usepackage{xcolor}

\usepackage[numbers]{natbib}
\usepackage{comment}
\usepackage{caption}
\usepackage{ragged2e}

\usepackage{hyperref}
\usepackage[normalem]{ulem}

\titleformat{\paragraph}[runin]{\bfseries}{}{0pt}{}
\titlespacing{\paragraph}{0pt}{1ex plus 0.5ex minus 0.2ex}{1em}

\hypersetup{colorlinks=true, citecolor=blue, linkcolor=red, urlcolor=green}

\newcommand{\updot}[1]{\raisebox{0.9pt}{$\stackrel{\bullet}{#1}$}} 

\providecommand{\tabularnewline}{\\}

\declaretheorem[name=Theorem]{theorem}  
\declaretheorem[sibling=theorem, name=Lemma]{lemma}

\declaretheorem[sibling=theorem, name=Proposition]{proposition}
\declaretheorem[sibling=theorem, style=definition, name=Definition]{definition}

\global\long\def\P{\mathbb{P}}%

\global\long\def\E{\mathbb{E}}%

\global\long\def\M{\mathbf{M}}%

\global\long\def\A{\mathcal{A}}%

\global\long\def\F{\mathcal{F}}%

\global\long\def\var{\text{Var}}%

\makeatother
\begin{document}
\title{Community Bail Fund Systems: Fluid Limits and Approximations}

\author{%
  \begin{tabular}[t]{c}
    Yidan Zhang\\
     Cornell University\\
    \texttt{yz2689@cornell.edu}
  \end{tabular}
  \qquad                        
  \begin{tabular}[t]{c}
    Jamol Pender\\
   Cornell University\\
    \texttt{jjp274@cornell.edu}
  \end{tabular}
}

\maketitle
\begin{abstract}
Community bail funds (CBFs) assist individuals who
have been arrested and cannot afford bail, preventing unnecessary
pretrial incarceration along with its harmful or sometimes fatal consequences.
By posting bail, CBFs allow defendants to stay at home and maintain
their livelihoods until trial.  This paper introduces new stochastic models
that combine queueing theory with classic insurance risk models
to capture the dynamics of the remaining funds in a CBF. We first analyze a model where \emph{all bail requests are accepted}. Although the remaining fund balance can go negative, this model provides insight for CBFs that are not financially constrained.  We then apply the Skorokhod map to make sure the CBF balance does not go negative and show that the Skorokhod map produces a model where requests are \emph{partially fulfilled}.   Finally, we analyze a model where bail requests can be blocked if there is not enough money to satisfy the request upon arrival.  Although the blocking model prevents the CBF from being negative, the blocking feature gives rise to new analytical challenges for a direct stochastic analysis.  Thus, we prove a functional law of large numbers or a fluid limit for the blocking model and show that the fluid limit is a distributed delay equation. We assess the quality of our fluid limit via simulation and show that the fluid limit accurately describes the large-scale stochastic dynamics of the CBF. Finally, we prove stochastic ordering results for the CBF processes we analyze.
\end{abstract}

\section{Introduction}

\subsection{Background: Community Bail Funds}

A community bail fund (CBF) is a pool of community resources to help post bail for individuals accused of crimes, thereby allowing them to be
released from jail while awaiting their trial. Historically, bail
funds in the United States date back to the 1920s, when the ACLU created
a fund to release individuals arrested for sedition during the First
Red Scare, see for example \citet{steinberg2018freedom, goldman2021freedom, simonson2016democratizing}. Since 2012, the popularity of community bail funds has surged dramatically and has become an important part of the judicial system.

Bail funds serve an important personal and social purpose for pretrial,
incarcerated defendants. Being incarcerated before a trial can lead
to job loss, family separation, and pressure to accept plea deals,
even if the person is innocent. Kalief Browder's unfortunate
suicide exemplifies the profound impact of excessive incarceration
on individuals and their families, see for example \citet{jones2015bronx, johnson2018justice}. By
providing financial assistance, bail funds strive to ensure that individuals
are not incarcerated simply due to their financial circumstances. 

Furthermore, community bail funds often advocate for broader criminal
justice reforms, including the reduction or elimination of cash bail
systems. They argue that these systems disproportionately affect low-
income individuals and perpetuate financial and ethnic/racial disparities
in the justice system. For example, Joseph E. Krakora, New Jersey's
public defender, criticizes the cash bail system for forcing defendants
to either endure months in jail or plead guilty for release, highlighting
its inherent injustice by linking freedom to wealth rather than risk, see for example
\citet{Bail2019}. By providing direct support and raising awareness
about the issue, bail funds play a crucial role in challenging the
reliance on cash bail and promoting a more equitable judicial system.
For more on these and other examples, see the introduction of \citet{gunluk2023simulating} and the references therein. 

In this paper, we develop stochastic models for understanding the dynamics of CBFs. These models could serve as the foundation for future work that guides decision-making in the judicial system to improve fairness, quality of life, and reduce recidivism. Moreover, it integrates with a growing literature applying operations research towards understanding the criminal justice system.

\subsection{Connection to Criminal Justice and Operations Research}

Traditionally, operations research has focused on decision-making to maximize the efficiency or profitability of some businesses and organizations. However, recent work in operations research instead applies stochastic modeling to understand the dynamics of the criminal justice system and its societal impact. Instead of profitability or efficiency, the focus shifts to understanding and optimizing the impact that the justice system has on the rehabilitation of accused individuals. For example, \citet{bakshi2024service} uses a two-stage queueing model and simulation analysis to provide a deeper understanding of the judicial delays caused by resource limitations in the Supreme Court of India and how to use the simulation analysis to improve the system. 
 \citet{li2024combining} that uses machine learning and queueing simulation to model incarceration diversion programs and to predict program census and staffing requirements under varying admission criteria. \citet{wang2018analyzing} use probabilistic modeling and statistical analysis to evaluate the U.S. backlog of untested sexual assault kits. They propose a framework to optimize resources and priorities and find that while prioritization offers little advantage, testing all kits is cost-effective with significant societal benefits.
\citet{blumstein1969models} develops integrated steady-state linear and feedback models of the criminal justice system to analyze offender flows, costs, and recidivism. Their approach assesses resource allocation, criminal career trajectories, and system sensitivity to policy changes. \citet{kleinberg2018human} shows machine learning can improve judicial bail decisions. Using gradient-boosted decision trees and quasi-random case assignment, they find judges mis-rank defendants. The study also highlights inefficiencies in judges' risk assessments and demonstrates how algorithms optimize decision-making while mitigating equity concerns.
These and other papers are part of a growing literature \citep{attari2021simulation, taxman2013simulation, usta2015assessing, 
bray2016multitasking, azaria2023justice, azaria2024alleviating,
kleinberg2018discrimination} that applies operations research to the improvement of the justice system.

In this spirit, the goal of this paper is to introduce new stochastic models that describe the dynamics of community bail funds. As will be described in the subsequent sections, we will combine stochastic queueing models with classic insurance risk models to build, simulate, and derive theoretical properties of CBFs.

\section{Constructing Stochastic Models for CBFs}

\begin{figure}[H]
\centering
\includegraphics[scale=0.4]{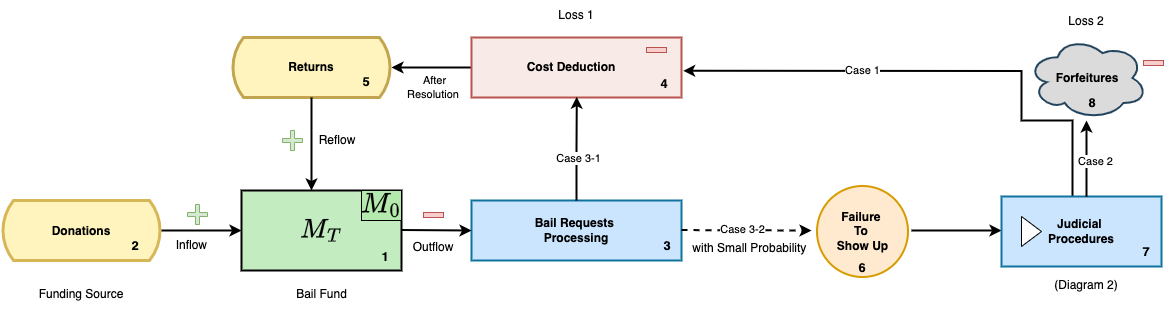}

\caption{\label{fig:main}Bail Fund Money Flow Operation System}
\end{figure}

In this section, we are going to introduce the three CBF stochastic models that we will analyze in the rest of this paper. As we will see later on, a CBF is affected by multiple
random factors, hence it is a stochastic system. In Figure \ref{fig:main}, we provide a visualization of the community
bail fund's money flow.
Each component in Figure \ref{fig:main} is labeled by a number index and each component is described using its assigned number, with sequences of these numbers explicitly showing how each specific component participates and interacts within the bail fund operations system.
Figure \ref{fig:main} also includes two components that are not in the model of this paper (6 and 7), but we will consider the analysis of those components in future work.

\begin{itemize}
\item {Initial Capital Endowment $(1)$}

{Based on an initial capital endowment, a bail fund
commences operations. }
\item {Donations $(2)$}

{During its operation, donations may arrive randomly,
flowing into the fund and augmenting its cash flow, i.e., $\left(2\rightarrow1\right)$.}
\item {Bail Requests $(3)$}

{Upon establishment, the fund opens to receive a
sequence of bail requests. The fund\textquoteright s acceptance policy
identifies and accepts defendants from incoming requests who meet
specific criteria, resulting in the corresponding funds to flow out
of the bail fund, i.e., $\left(1\rightarrow3\right)$.}
\item {Cost Deductions and Forfeiture $(4,8)$}

{Accepted bail amounts undergo cost deductions during
the judicial process ($3\rightarrow4$).
If a defendant who has been granted bail fails to appear in court
and falls into state $(6)$, a more complex judicial procedure is
triggered. Various factors in judicial systems can result in partial
cost deductions from the bail amount for administrative fees, or in
the forfeiture of the entire bail amount, i.e., $\left(6\rightarrow7\rightarrow4\right)$
or $\left(6\rightarrow7\rightarrow8\right)$. Although not modeled
directly, this cost can be absorbed into a \textquotedblleft poundage\textquotedblright{}
fee term in the models considered in this paper. However, for future
explorations, we note that it may be also interesting to model these
aspects separately.}
\item {Bail Residual Returns $(5)$}

{Once defendants fulfill their court obligations,
any residual bail, after deducting fees incurred during the judicial
process, flows back into the bail fund, i.e., $\left(5\rightarrow1\right)$.
Consequently, the acceptance policy also determines the sequence of
bail returns: whose residuals return, how much is returned after fees,
and when they re-enter the fund.}
\end{itemize}
The fund\textquoteright s balance fluctuates randomly, decreasing
due to outflow of accepted bail requests and increasing from inflow
of donations and the subsequent return of bail residuals. Moreover,
variations in the judicial system\textquoteright s operations add
another factor disrupting the bail return process. According to \citet{graef2023systemic}, systemic failures to show up, including the nonappearance of police officers, victim absenteeism, defendant absences, as well as discrepancies among
different judicial areas, significantly prolong case proceedings. This creates unnecessary delays
and uncertainty, which can indirectly complicate related processes
such as bail returns. 

A bail fund inherently operates as a stochastic
process with balance-dependent jumps.  In practice, a bail
fund\textquoteright s policy accounts for limited resources by using
the current balance as a constraint to decide whether to accept new bail requests. For example, the
acceptance of the first bail request reduces the fund's balance,
thereby affecting the decision on the second request and making it less likely to be accepted. As each defendant fulfills their obligations and residual funds are returned, the balance is adjusted again, affecting
subsequent bail decisions. Consequently, whether each bail request
is fulfilled depends on the entire sample path of the CBF, beginning with the initial bail decision. The fact that the CBF is path-dependent makes exact mathematical analysis of the CBF challenging.

In fact, the model depicted in Figure \ref{fig:main} is intractable not only because of the multitude
of interacting random factors, but also the path-dependence. Thus, our goal in this work is to construct modified CBF models that are  more tractable and amenable to stochastic analysis. For our modified
CBF models, we prove fluid limits or functional law of large
numbers for the CBF. We show that the limit of our
scaled stochastic models converge to functional differential equations. In the sequel, we describe our three modified CBF models and how they are special cases of the Figure \ref{fig:main}.

\subsection{Infinite Acceptance Model}

We initially simplify our model by eliminating
the complexities introduced by its state dependent
jumps. We construct an \emph{infinite acceptance model}
with no blocking, assuming that the bail fund could
borrow money without limited credit. This model
infinitely accepts all bail requests without constraints, removing
the dependence on prior events affecting the fund's
balance. This simplification enabled straightforward analysis, leading
to the almost sure convergence of its scaled process to its fluid
limit, which reveals its large-scale dynamics. Of course,
a natural consequence of this model is that the balance becomes negative
when bail requests exceed donations. 

\begin{figure}[H]
\centering
\includegraphics[scale=0.35]{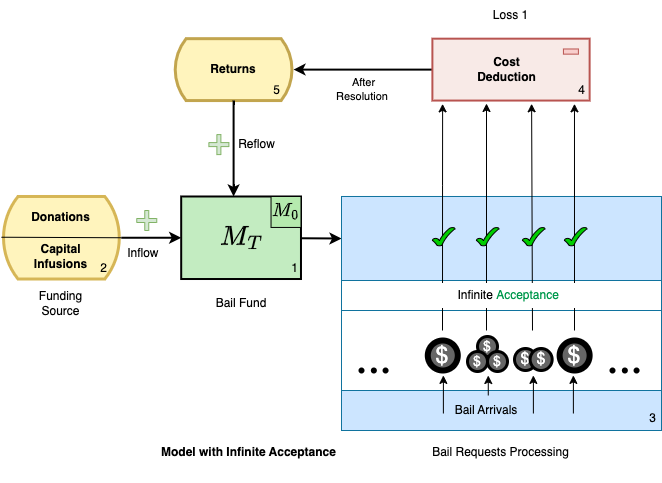}\caption{Infinite Acceptance Model}
\end{figure}

\subsection{The Partial Fulfillment Model}
Since the infinite acceptance model could become negative, this motivates us to ensure the balance
remains non-negative. Thus, we apply the Skorokhod map to the infinite acceptance model, transforming
it into an always non-negative process. Our results (Theorem \ref{thm:partial} and Corollary \ref{col:partial}) prove that the
Skorokhod map\textquoteright s output, applied to a generalized class
of input processes with independent jumps, leads to a model in which each defendant's bail request 
is partially fulfilled if there is not enough money available upon request. These results sharpen our understanding
of the Skorokhod map's interpretable output process, unveiling its
mathematical structure and implications. 

\begin{figure}[H]
\centering
\includegraphics[scale=0.35]{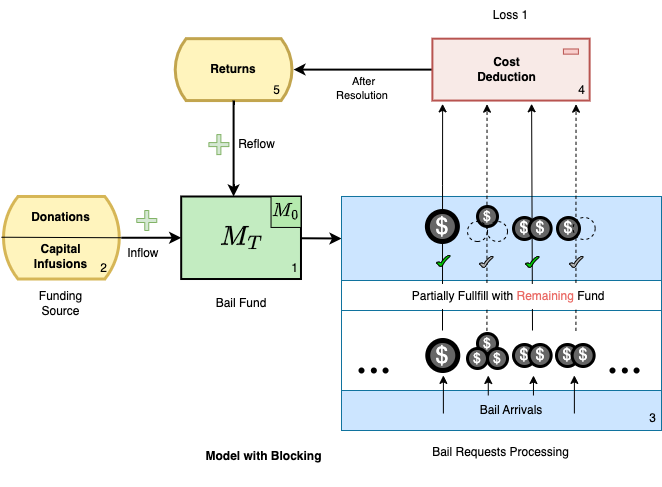}\caption{Partial Fulfillment Model}
\end{figure}

\subsection{The Blocking Model}

In the case where defendants can not be partially funded, we design a model where defendants are blocked from the CBF if the balance is lower than their request. This model captures the state dependent nature of
the bail fund by introducing a blocking decision term into the model.
This approach assumes limited fund recourse, where the current balance
solely determines whether bail requests are accepted. If a bail request
exceeds the total accumulated balance from all prior events, then
such a request is blocked and will not fulfilled. Only requests within
the remaining balance are accepted, leading to outflows matching the
corresponding amount from the bail fund. Despite the challenges of
state dependent jumps, we prove that, with appropriate scaling, our
bail process models converge to fluid limits described by distributed delay
equations.

\begin{figure}[H]
\centering
\includegraphics[viewport=0bp 0bp 662bp 487bp,scale=0.35]{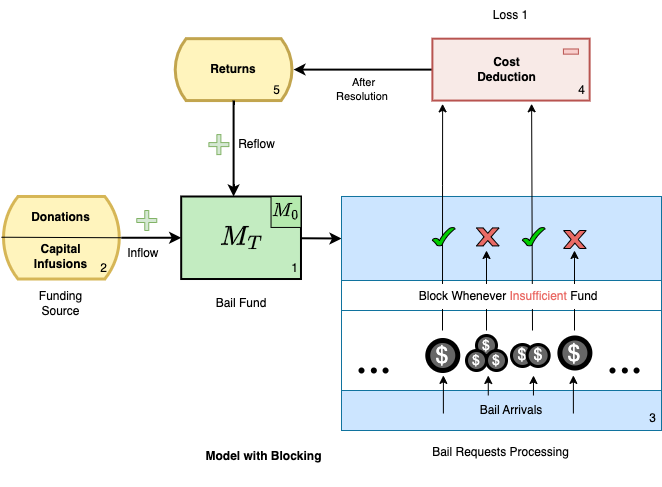}\caption{Blocking Model}
\end{figure}

\subsection{Stochastic Ordering Analysis of CBF Processes}

To gain more insight, we prove a stochastic ordering relationship between the different processes we study in this paper. This analysis compares their balance values across
models and over time. The bail returns process was found to disrupt
the inherent ordering: models without the return component maintained
consistent orderings, whereas those with it experienced disrupted
and inconsistent rankings. Consequently, we partition our
stochastic ordering results into two cases: one with returns and one without returns.

\subsection{Notation of Paper}

\label{def:SideNotations}We establish the following notations to
be used throughout this paper: 

\begin{center}
\begin{tabular}{|c|c|}
\hline 
$M_{0}$ & initial capital money value\tabularnewline
\hline 
$\left\{ d_{i}\right\} $ & $i^{th}$ donation jump size following some exponential distribution\tabularnewline
\hline 
$\left\{ b_{j}\right\} $ & $j^{th}$ bail request jump size following some exponential distribution\tabularnewline
\hline 
$d^{*}$ & $\E\left[d_{i}\right]$\tabularnewline
\hline 
$b^{*}$ & $\E\left[b_{j}\right]$\tabularnewline
\hline 
$p^{*}$ & $\E\left[p_{j}\right]$\tabularnewline
\hline 
$N^{d}\left(t\right)$ & Poisson process describing the arrival of donations\tabularnewline
\hline 
$N^{b}\left(t\right)$ & Poisson process describing the arrival of bail requests\tabularnewline
\hline 
$\lambda_{d}$ & arrival rate of the Poisson process $N^{d}\left(t\right)$\tabularnewline
\hline 
$\lambda_{b}$ & arrival rate of the Poisson process $N^{b}\left(t\right)$\tabularnewline
\hline 
$\left\{ a_{j}\right\} $ & arrival time for the $j^{th}$ bail request\tabularnewline
\hline 
$\left\{ s_{j}\right\} $ & trial time for the $j^{th}$ bail request\tabularnewline
\hline 
$\left\{ p_{j}\right\} $ & poundage rate for the $j^{th}$ bail request\tabularnewline
\hline 
$f_{b}\left(\cdot\right)$ & probability density function of $\left\{ b_{j}\right\} $\tabularnewline
\hline 
$F_{b}\left(\cdot\right)$ & cumulative distribution function of $\left\{ b_{j}\right\} $\tabularnewline
\hline 
$f_{s}\left(\cdot\right)$ & probability density function of $\left\{ s_{j}\right\} $\tabularnewline
\hline 
$F_{s}\left(\cdot\right)$ & cumulative distribution function of $\left\{ s_{j}\right\} $\tabularnewline
\hline 
\end{tabular}
\end{center}

\section{\label{sec:Inf}Infinite Acceptance Model}

In this section, we analyze the \textbf{Infinite Acceptance Model}. The infinite acceptance model assumes that every bail request
is accepted immediately upon arrival and the CBF has no financial constraints.  This model is convenient since it would provide insight about how much money is needed to serve all clients that need access to the bail fund. As such, the infinite acceptance model has the following stochastic process representation

\begin{align*}\label{Model:Inf}
\underbrace{M^{\infty,R^{\infty}}(t)}_{\text{Total Funds Available}} & =\underbrace{M_{0}}_{\text{Initial Capital}}+\underbrace{\sum_{i=1}^{N^{d}(t)}d_{i}}_{\text{Total Donations}}-\underbrace{\sum_{j=1}^{N^{b}(t)}b_{j}}_{\text{Total Bail Paid Out}} +\underbrace{\sum_{j=1}^{N^{b}(t)}\left(1-p_{j}\right)b_{j}\left\{ t>a_{j}+s_{j}\right\} }_{\text{Total Funds Returned}}.
\end{align*}
The first two terms on the right-hand side correspond to the initial
funds and donations. The third term indicates the total amount of bail
paid for defendants awaiting trial up to time $t$, assuming bail requests
are infinitely accepted upon the defendant's request. The final term accounts
for the total bail returned by the courts by time $t$, excluding
any forfeitures.

This process models donations and bail requests as independent, unlimited jumps. The simplicity of the process allows for direct derivations of its properties such as expectation, fluid limit, and variance. The
derivation of these results is presented in \citet{gunluk2023simulating}.  In later sections, we will examine more realistic processes at the cost of mathematical tractability. For this simplified model, we will derive the expectation and fluid limit of this model which quantifies its deterministic average and limiting behaviors, respectively. 
For the fluid limit in particular, we first define a properly scaled version of our model
in the next section.

\subsection{\label{subsec:InfScaling} Scaling of the Infinite Acceptance Process}

We scale the infinite acceptance model with the following two step procedure. First, the arrival rate
parameters for donations and bail requests, along with the time-to-trial
service rate, are increased by a factor of $\eta$. Meanwhile, the
magnitude of each individual jump for both donations and bail requests
is proportionally reduced by $\frac{1}{\eta}$. We now define the
$\eta$-scaled version of $M^{\infty,R^{\infty}}(t)$ as follows:

\begin{equation}
M_{\eta}^{\infty,R^{\infty}}(t)=\frac{1}{\eta}\left(\eta M_{0}\right)+\frac{1}{\eta}\sum_{i=1}^{N_{\eta}^{d}( t)}d_{i}-\frac{1}{\eta}\sum_{j=1}^{N_\eta^{b}( t)}b_{j}+\frac{1}{\eta}\sum_{j=1}^{N_{\eta}^{b}( t)}\left(1-p_{j}\right)b_{j}\left\{t>a_{j,\eta}+s_{j}\right\} .\label{eq:MIR}
\end{equation}

We chose this scaling framework to study large-scale
bail funds. By increasing the arrival rate and time-to-trial service
rate, it captures a high-frequency scenario in which the bail fund
manages larger volumes of donations and the bail requests as $\eta$
grows. This approach is especially relevant for large-scale nonprofit
organizations, such as The Bail Project and the National Bail Fund
Network, which operate under similar high-volume conditions. While we scale the volume of demand up, we scale down the magnitude of each individual's contribution to the bail fund. In next section, we show how the scaling we have defined converges to the fluid limit.
\subsection{\label{subsec:FLConvergeInf}Convergence to Fluid Limit}

Our goal in this Section \ref{subsec:FLConvergeInf} is to prove the
convergence of our $\eta$-scaled model defined in Section \ref{subsec:InfScaling}
to its fluid limit. We will establish this by proving:
\begin{enumerate}
\item Each model component---donations, bail requests, and returns, converges
almost surely (Section \ref{subsubsec:CPConverge}).
\item The model converges to its fluid limit (Section \ref{subsubsec:InfConverge}).
\end{enumerate}

\subsubsection{\label{subsubsec:CPConverge}Each Component Converges Almost Surely}

To help obtain the almost sure convergence of each component, we will
first prove each component is a martingale, and then apply Doob's
maximal inequality for continuous-time submartingales.
\begin{restatable}[label={lem:MgMG}]{lemma}{MgMGLemma}
Consider the process with scale parameter $\eta$
\[
\M_{\bowtie}^{\eta}(t):=\frac{1}{\eta}\sum_{i=1}^{N_\eta\left( t\right)}v_{i,\eta}-\lambda v_{\eta}^{*}t.
\]
Further, assume: 
\begin{itemize}
\item the summands $v_{i,\eta}$ conditional on $N_\eta\left(t\right)$
are i.i.d distributed random variables, 
\item and $\E\left[v_{i,\eta}\Bigg|N_\eta\left( t\right)\right]=v_{\eta}^{*}$
is a constant depending only on $\eta$ . 
\end{itemize}
Then, $\M_{\bowtie}^{\eta}(t)$ is a martingale with respect to the
filtration $\left\{ \F_{s}\right\} ,$ which contains all arrival
and jump size information until time $s$: 
\[
\F_{s}:=\sigma\left(N_\eta\left( s\right),\left\{ v_{i,\eta}\right\} _{i=1}^{N_\eta\left( s\right)}\right),
\]
\end{restatable}
\begin{proof}
This proof is in Appendix \ref{app:MgMG}. 
\end{proof}

\begin{restatable}[label={lem:Gcpconverge}]{lemma}{GcpconvergeLemma}
Consider a general process with scale parameter
$\eta$ defined by 
\[
\sum_{i=1}^{N_{\eta}\left( t\right)}v_{i,\eta}
\]
where $N_{\eta}\left( t\right)$ is the Poisson process with rate $\eta\lambda$
at time $t$, and we assume the following: 
\begin{itemize}
\item the summands $v_{i,\eta}$ conditional on $N_{\eta}\left( t\right)$
are i.i.d distributed random variables, 
\item and $\E\left[v_{i,\eta}\Bigg|N_{\eta}\left( t\right)\right]=v_{\eta}^{*}$
is a constant depending only on $\eta$ . 
\end{itemize}
Then, for any fixed time T, as $\eta\rightarrow\infty$,

\[
\sup_{0\leq t\leq T}\left|\frac{1}{\eta}\sum_{i=1}^{N_{\eta}\left( t\right)}v_{i,\eta}-\text{\ensuremath{\lambda v_{\eta}^{*}t}}\right|\underset{a.s.}{\to}0.
\]
\end{restatable}

\begin{proof}
This proof is in Appendix \ref{app:Gcpconverge}.
\end{proof}

\begin{restatable}[label={lem:ReturnAS}]{lemma}{ReturnASLemma}
 (Centered Return Component Converges Almost Surely) For
any fixed time $T$, as $\eta\rightarrow\infty$,

\[
\sup_{0\le t\leq T}\left|\frac{1}{\eta}\sum_{j=1}^{N_{\eta}^{b}(t)}\left(1-p_{j}\right)b_{j}\left\{ t>a_{j,\eta}+s_{j,\eta}\right\} -\lambda_{b}\left(1-p^{*}\right)b^{*}\int_{0}^{t}F_{s}\left(t-v\right)dv\right|\underset{a.s.}{\rightarrow}0.
\]
\end{restatable}

\begin{proof}
The proof of this result is contained in Appendix \ref{app:ReturnAS}. 
\end{proof}

\subsubsection{\label{subsubsec:InfConverge} Convergence of $\eta$-scaled No Blocking
Model to the Fluid Limit}

\begin{theorem}
\label{thm:FLMIR}(Fluid Limit of the Infinite Acceptance Process)
Consider the $\eta$-sclaed process in \eqref{eq:MIR}:
\begin{align*}
M_{\eta}^{\infty,R^{\infty}}(t) & =\frac{1}{\eta}\left(\eta M_{0}\right)+\frac{1}{\eta}\sum_{i=1}^{N_{\eta}^{d}(t)}d_{i}-\frac{1}{\eta}\sum_{j=1}^{N_{\eta}^{b}(t)}b_{j}+\frac{1}{\eta}\sum_{j=1}^{N_{\eta}^{b}(t)}\left(1-p_{j}\right)b_{j}\left\{ t>a_{j,\eta}+s_{j}\right\} ,
\end{align*}

 as well as the deterministic path $m^{\infty,R^{\infty}}\left(t\right):$
\begin{align*}
m^{\infty,R^{\infty}}\left(t\right) & :=M_{0}+d^{*}\lambda_{d}t-b^{*}\lambda_{b}t+\left(1-p^{*}\right)b^{*}\lambda_{b}\int^t_0F_s(v)dv.
\end{align*}
Then,\
\begin{align*}
\sup_{0\leq u\leq T}\left|M_{\eta}^{\infty,R^{\infty}}(t)-m^{\infty,R^{\infty}}\left(t\right)\right| & \underset{a.s.}{\rightarrow}0
\end{align*}
as $\eta\rightarrow\infty$. Therefore, $m^{\infty,R^{\infty}}\left(t\right)$
is the fluid limit of $M_{\eta}^{\infty,R^{\infty}}$as $\eta\rightarrow0$. 
\end{theorem}
\begin{proof}
By definition, we want to show that for all $\epsilon>0$, 
\begin{align*}
P\left(\lim\sup_{\eta\rightarrow\infty}\left\{ \omega\in\Omega:\sup_{t\in\left[0,T\right]}\left|M_{\eta}^{\infty,R^{\infty}}\left(t\right)-m^{\infty}\left(t\right)\right|<\epsilon\right\} \right) & =1.
\end{align*}
Notice, by Lemma \ref{lem:MgMG} ,

\begin{align*}
\sup_{t\in\left[0,T\right]}\left|M_{\eta}^{\infty,R^{\infty}}(t)-m^{\infty,R^{\infty}}\left(t\right)\right| & =\sup_{t\in\left[0,T\right]}\left|\M_{d,\eta}^{\infty}(t)+\M_{b,\eta}^{\infty}(t)+\M_{r,\eta}^{\infty}(t)\right|.
\end{align*}
Define the following notation for brevity:
\[
S_{T}(M,\epsilon):=\left\{ \omega\in\Omega:\sup_{t\in[0,T]}\left|M(t)\right|\leq\epsilon\right\},
\]
and $S^c_{T}(\cdot, \cdot)$ is its complement.

Observe,\begin{gather*}
\P\left(\lim\sup_{\eta\rightarrow\infty}\left\{ \omega\in\Omega:\sup_{t\in\left[0,T\right]}\left|M_{\eta}^{\infty,R^{\infty}}\left(t\right)-m^{\infty}\left(t\right)\right|<\epsilon\right\} \right)\\
\geq\P\left(\lim\sup_{\eta\rightarrow\infty}S_{T}\left(\M_{d,\eta}^{\infty},\frac{\epsilon}{3}\right)\text{\ensuremath{\bigcap}}\lim\sup_{\eta\rightarrow\infty}S_{T}\left(\M_{b,\eta}^{\infty},\frac{\epsilon}{3}\right)\bigcap\lim\sup_{\eta\rightarrow\infty}S_{T}\left(\M_{r,\eta}^{\infty},\frac{\epsilon}{3}\right)\right)\\
\geq1-\P\left(\lim\sup_{\eta\rightarrow\infty}S_{T}^{c}\left(\M_{d,\eta}^{\infty},\frac{\epsilon}{3}\right)\right)-\P\left(\lim\sup_{\eta\rightarrow\infty}S_{T}^{c}\left(\M_{b,\eta}^{\infty},\frac{\epsilon}{3}\right)\right)\\
-\P\left(\lim\sup_{\eta\rightarrow\infty}S_{T}^{c}\left(\M_{r,\eta}^{\infty},\frac{\epsilon}{3}\right)\right),
\end{gather*}
where the last line follows from the Bonferroni Inequality.
Then, by Lemma \ref{lem:Gcpconverge},
\begin{align*}
\sup_{0\leq t\leq T}\left|\M_{d,\eta}^{\infty}(t)\right| & \underset{a.s.}{\rightarrow}0,\\
\sup_{0\leq t\leq T}\left|\M_{b,\eta}^{\infty}(t)\right| & \underset{a.s.}{\rightarrow}0,
\end{align*}
and by Lemma \ref{lem:ReturnAS},
\begin{align*}
\sup_{0\leq t\leq T}\left|\M_{r,\eta}^{\infty}(t)\right| & \underset{a.s.}{\rightarrow}0.
\end{align*}
Therefore, 
\begin{align*}
\P\left(\lim\sup_{\eta\rightarrow\infty}\left\{ \omega\in\Omega:\sup_{t\in\left[0,T\right]}\left|M_{\eta}^{\infty,R^{\infty}}\left(t\right)-m^{\infty}\left(t\right)\right|<\epsilon\right\} \right) & \geq1-0-0-0=1.
\end{align*}
This completes the proof. 
\end{proof}

\paragraph{Remark.}
We notice this fluid limit of the infinite acceptance model matches its expectation:
\begin{eqnarray*}
  \mathbb{E} [M^\infty(t)] &=&    M_0  +  \lambda^{(d)} t d^* - \lambda^{(b)} t b^* + \lambda^{(b)} \cdot b^* \cdot ( 1- p^*)  \cdot  \int^{t}_{0} F_s(t-u) du  .
\end{eqnarray*}
See \citet{gunluk2023simulating} for the proof of this result.

\providecommand{\tabularnewline}{\\}

\makeatother
\subsection{\label{sec:MIRSim} Numerical Analysis}

\begin{figure}[h]
\centering
\includegraphics[scale=0.5]{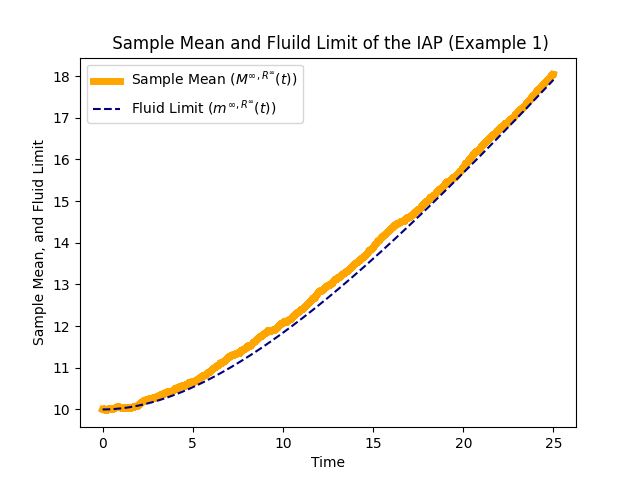}
 \includegraphics[scale=0.5]{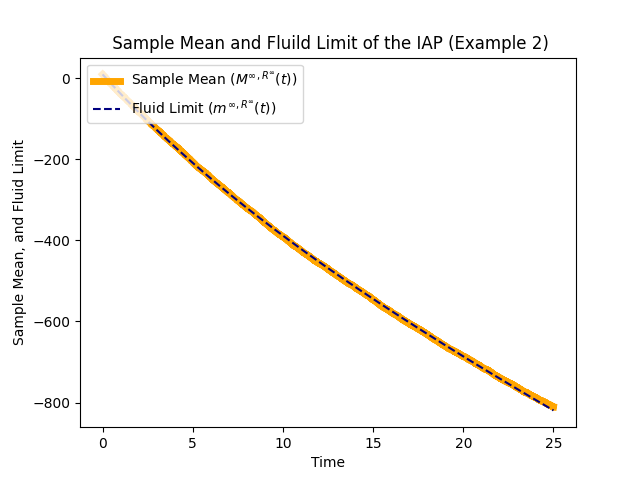}
\\
\begin{tabular}[H]
{|c|c|c|c|c|c|c|c|}
\hline 
Example/Parameters & $M_{0}$ & $\lambda_{b}$ & $\lambda_{b}$ & $d_{i}$ & $b_{j}$ & $p_{j}$ & $s_{j}$\tabularnewline
\hline 
\hline 
Example 1 & 10 & 1 & 1 & $\text{Expo}(1)$ & $\text{Expo}(1)$ & $\text{Unif}\left[0,1\right]$ & $\text{Expo}(10)$\tabularnewline
\hline 
Example 2 & 10 & 1 & 1 & $\text{Expo}(1)$  &  $\text{Expo}(50)$ & $\text{Unif}\left[0,1\right]$ & $\text{Expo}(10)$\tabularnewline
\hline 
\end{tabular}

\caption{\textbf{Mean and Fluid limit of the Infinite Acceptance Process (IAP).}
}
\label{fig:inf_mean_variance}
\end{figure}
We will first discuss the numerical results in Figure \ref{fig:inf_mean_variance}. In \textbf{Example 1}, the sample mean matches the theoretical expectation. Initially, both remain roughly constant.
Since donations $\left\{ d_{i}\right\} _{i\geq1}$ and bail requests
$\left\{ b_{j}\right\} _{j\geq1}$ share identical arrival rates ($\lambda_{d}=\lambda_{b}=1$)
and average jump sizes ($Expo\left(1\right)$) , they produce equal
mean increment per unit time. Around $t=10,$ both begin to rise.
This occurs because bail amounts, proportionally reduced by the poundage
fee $p_{j}\sim\text{Unif\ensuremath{\left[0,1\right]}}$, return after
approximately 10 time units after their initial departure, and this
is consistent with our service time $s_{j}\sim\text{Expo\ensuremath{\left(10\right)}}$. In \textbf{Example 2}, both the sample mean and theoretical expectation
decrease over time. This occurs because the average bail amount per
unit $\left(\lambda_{b}\E\left[b_{j}\right]=50\right)$ outsizes the
average donation amount per unit time. 
\FloatBarrier

Second, the simulation result in Figure \ref{fig:MBRFluid} (using the same set of parameters as the two examples above) empirically verifies that the $\eta$-scaled
infinite acceptance process $M_{\eta}^{\infty,R^{\infty}}\left(t\right)$ converges to its
fluid limit $m^{\infty,R^{\infty}}\left(t\right)$. This is mathematically proven in Theorem \ref{thm:FLMIR} of Section \ref{subsubsec:InfConverge} .

\begin{figure}[H]
\centering
\includegraphics[scale=0.5]{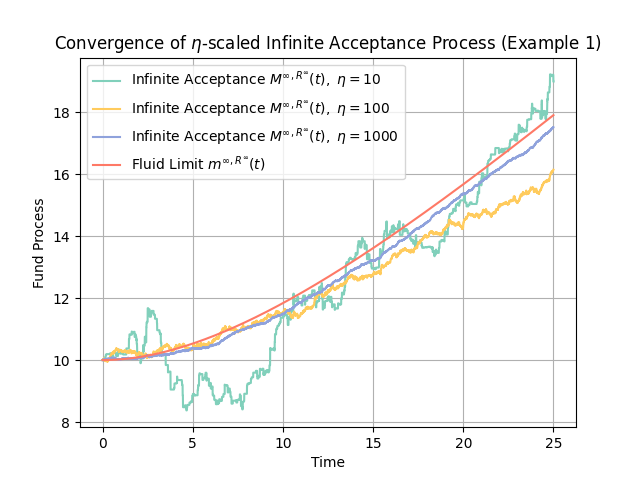}
\includegraphics[scale=0.5]{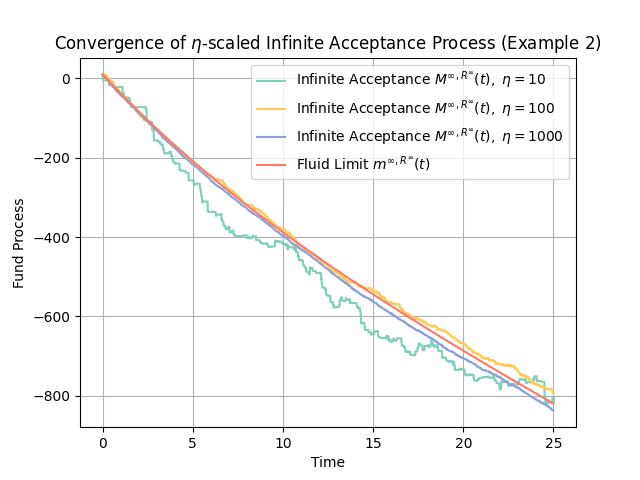}
\caption{Convergence of the Infinite Acceptance Process}

\label{fig:MInfFluid}
\end{figure}

\subsection{Steady State Behavior of  Infinite Acceptance Fluid Model}

\begin{proposition}\label{prop:steadyMIR}
Recall the fluid limit equation in Theorem \ref{thm:FLMIR}.
\[
m^{\infty,R^{\infty}}\left(t\right):=M_{0}+d^{*}\lambda_{d}t-b^{*}\lambda_{b}t+\left(1-p^{*}\right)b^{*}\lambda_{b}\int_{0}^{t}F_{s}\left(v\right)dv.
\]

The long-term steady state behavior of the fluild limit $\lim_{t\rightarrow\infty}m^{\infty,R^{\infty}}\left(t\right)$
is determined by the sign of the quantity $d^{*}\lambda_{d}-p^{*}b^{*}\lambda_{b}$
as follows:

\begin{align*}
\lim_{t\rightarrow\infty}m^{\infty,R^{\infty}}\left(t\right) & =\begin{cases}
M_{0}+d^{*}\lambda_{d}t-b^{*}\lambda_{b}t+\left(1-p^{*}\right)b^{*}\lambda_{b} & \text{if }d^{*}\lambda_{d}-p^{*}b^{*}\lambda_{b}=0\\
+\infty & \text{if }d^{*}\lambda_{d}-p^{*}b^{*}\lambda_{b}>0\\
-\infty & \text{if }d^{*}\lambda_{d}-p^{*}b^{*}\lambda_{b}<0
\end{cases}
\end{align*}
\end{proposition}
\begin{proof}
Notice the time derivative of this fluild limit is given by: 
\begin{align*}
\frac{d}{dt}\left(m^{\infty,R^{\infty}}\left(t\right)\right) & =d^{*}\lambda_{d}-b^{*}\lambda_{b}+\left(1-p^{*}\right)b^{*}\lambda_{b}F_{s}\left(t\right).
\end{align*}

Since $F_{s}\left(t\right)\rightarrow1$ as $t\rightarrow\infty,$
the fluild limit approaches a steady-state behavior described by
\begin{align*}
\lim_{t\rightarrow\infty}\frac{d}{dt}m^{\infty,R^{\infty}}\left(t\right) & =d^{*}\lambda_{d}-b^{*}\lambda_{b}+\left(1-p^{*}\right)b^{*}\lambda_{b}\\
 & =d^{*}\lambda_{d}-p^{*}b^{*}\lambda_{b}.
\end{align*}

This steady state behavior then diverges into three distinct cases,
which we now discuss individually.

\paragraph{Case 1: $d^{*}\lambda_{d}-p^{*}b^{*}\lambda_{b}=0.$} In this case, the expected
donation amount per unit time $d^{*}\lambda_{d}$ gained equals the
expected poundage fee lost from the bail amount $p^{*}b^{*}\lambda_{b}$.
Then, the fluid limit of the infinite acceptance model stabilizes
at the steady state value $M_{0}+d^{*}\lambda_{d}t-b^{*}\lambda_{b}t+\left(1-p^{*}\right)b^{*}\lambda_{b}$.

\paragraph{Case 2: $d^{*}\lambda_{d}-p^{*}b^{*}\lambda_{b}>0.$}When the expected donation
amount per unit time $d^{*}\lambda_{d}$ gained exceeds the expected
poundage fee of bail amount $p^{*}b^{*}\lambda_{b}$ lost, then the
fluid limit of the infinite acceptance model increases without bound,
diverging to positive infinity.

\paragraph{Case 3: $d^{*}\lambda_{d}-p^{*}b^{*}\lambda_{b}<0.$} When the average donation
amount per unit time $d^{*}\lambda_{d}$ gained is smaller the average
poundage fee lost from bail amount $p^{*}b^{*}\lambda_{b}$ , then
the fluid limit of the infinite acceptance model drives into a deficit
over time, diverging to negative infinity.

This completes the proof.
\end{proof}
\section{\label{sec:SkrkandParitial}Partial Fulfillment via the Skorokhod Map}

In this section, we adjust the infinite acceptance model to capture the fact that
real bail funds typically avoid deficits. The infinite acceptance model examined in the previous section could become negative, corresponding to being in deficit. In practice, it is more likely that a CBF would reject or partially fund a bail request when it cannot afford the entire amount requested. To
update our model such that is always non-negative, one immediate idea
inspired by the queueing literature is to transform the process using
the Skorkohod map defined as follows.
\begin{definition}
\label{def:SkrkMap}(Skorokhod Map \citep{skorokhod1961stochastic}) 
Let $\mathcal{D}[0,\infty)$ be the space of real-valued c\`adl\`ag functions defined on $[0,\infty)$, i.e., functions that are right-continuous and have left limits.
For an arbitrary process realization $P\left(t\right)\in \mathcal{D}[0,\infty)$, the Skorokhod map $\phi:\mathcal{D}[0,\infty)\rightarrow\mathcal{D}[0,\infty) $ is defined by:
\begin{gather*}
\phi\left[P\right]\left(t\right)=P\left(t\right)-\inf_{s\leq t}\left\{ 0,P\left(s\right)\right\} .
\end{gather*}
\end{definition}
The question now is: What is the resulting model
after the transformation? We identify later, in our main result Theorem \ref{thm:partial},
that for a class of generalized processes with necessarily independent
outflow jump sizes, the Skorokhod map produces an interpretable ``Partial
Fulfillment'' process. We call that class of generalized input processes
as \textquotedblleft simple bail processes\textquotedblright{} of
the following form.
\begin{definition}
\label{def:simple}(Simple Bail Process) We call processes of the
following form\textit{ simple bail} processes. 
\begin{align*}
M\left(t\right) & =G\left(t\right)-\sum_{j=1}^{N^{b}\left(t\right)}b_{j}\\
 & =G\left(t\right)-\sum_{j=1}^{\infty}b_{j}\left\{ a_{j}\leq t\right\} 
\end{align*}
where $G\left(t\right)$ is an arbitrary \textit{non-decreasing }process. 
\end{definition}
$G\left(t\right)$ here represents the \textquotedblleft Growth
Term\textquotedblright{} that absorbs both donations and returns.
For example, G(t) could contain all the donations received by time
t, and all the bail returns by time t. 

We wish to show that the Skorokhod map of the \textquotedblleft \textit{simple
bail} process\textquotedblright{} $M\left(t\right)$
(left hand side below) in Definition \ref{def:simple} is equal to
a process that \textquotedblleft partial fulfills\textquotedblright{}
the bail requests (right hand side below):

\[
\phi\left[M\right]\left(t\right)=G\left(t\right)-\underbrace{\sum_{i=1}^{N^{b}\left(t\right)}b_{j}\wedge\phi\left[M\right]\left(a_{j}-\right)}_{\text{Partially Fulfills Bail Requests}},
\] where $M(t)$ is a  simple bail process in Definition \ref{def:simple}.

This expression suggests that we should focus on
the connection between the $\inf_{s\leq t}\left\{ 0,M\left(s\right)\right\} $
and the $\sum_{j=1}^{N^{b}\left(t\right)}b_{j}\wedge\phi\left[M\right]\left(a_{j}-\right)$
terms. This inspires the following lemma before the main theorem.

\begin{restatable}[label={lem:shift}]{lemma}{shiftLemma}
Consider the simple bail process $M\left(t\right)$
in Definition \ref{def:simple}, and let $M^{*}\left(t\right):=\phi\left[M\right]\left(t\right)$
be the output process of the Skorokhod map on $M\left(t\right)$.
Assume that $G\left(a_{j}-\right)=G\left(a_{j}+\right)$ for all the
bail request arrival times $a_{j}$. Then, the Skorokhod shift, $\inf_{s\leq t}\left\{ 0,M\left(s\right)\right\} $,
changes only at times $t=a_{j}$ satisfying 
\begin{gather*}
M^{*}\left(a_{j}-\right)<b_{j}
\end{gather*}
and, at these times, 
\begin{gather*}
\inf_{s\leq a_{j}+}\left\{ 0,M\left(s\right)\right\} -\inf_{s\leq a_{j}-}\left\{ 0,M\left(s\right)\right\} =M^{*}\left(a_{j}-\right)-b_{j}.
\end{gather*}
\end{restatable}
\begin{proof}
From Definition \ref{def:simple}, $M\left(t\right)=G\left(t\right)-\sum_{j=1}^{N^{b}\left(t\right)}b_{j}$.
Since $G\left(t\right)$ is non-decreasing, the only times where $M\left(t\right)$
can get smaller is at the bail request arrival times $\left\{ a_{j}\right\} _{j=1}^{\infty}$.
The arrival times $t=a_{j}$ where $\inf_{s\leq t}\left\{ 0,M\left(s\right)\right\} $
changes are the following:

\begin{align*}
\left\{ a_{j}:\ \inf_{s\leq a_{j}+}\left\{ 0,M\left(s\right)\right\} <\inf_{s\leq a_{j}-}\left\{ 0,M\left(s\right)\right\} \right\}  & =\left\{ a_{j}:\ M\left(a_{j}+\right)<\inf_{s\leq a_{j}-}\left\{ 0,M\left(s\right)\right\} \right\} \\
 & =\left\{ a_{j}:\ M\left(a_{j}-\right)-b_{j}<\inf_{s\leq a_{j}-}\left\{ 0,M\left(s\right)\right\} \right\} \\
 & =\left\{ a_{j}:\ M\left(a_{j}-\right)-\inf_{s\leq a_{j}-}\left\{ 0,M\left(s\right)\right\} <b_{j}\right\} \\
 & =\left\{ a_{j}:\ M^{*}\left(a_{j}-\right)<b_{j}\right\} .
\end{align*}
The first equality follows from the fact that, if $\inf_{s\leq t}\left\{ 0,M\left(s\right)\right\} $
achieves a new lowest-value at $t=a_{j}$, this lowest value must
be $M\left(a_{j}+\right)$. The second equality follows from the fact,
since $G\left(a_{j}-\right)=G\left(a_{j}+\right)$, the the change
in $M\left(t\right)$ at $a_{j}$ is exactly $b_{j}$. This proves
the first part of this lemma.

To prove the second part, observe that at a $t=a_{j}$ where $\inf_{s\leq t}\left\{ 0,M\left(s\right)\right\} $
changes 
\begin{align*}
\inf_{s\leq a_{j}+}\left\{ 0,M\left(s\right)\right\} -\inf_{s\leq a_{j}-}\left\{ 0,M\left(s\right)\right\}  & =M\left(a_{j}+\right)-\inf_{s\leq a_{j}-}\left\{ 0,M\left(s\right)\right\} \\
 & =M\left(a_{j}-\right)-b_{j}-\inf_{s\leq a_{j}-}\left\{ 0,M\left(s\right)\right\} \\
 & =M^{*}\left(a_{j}-\right)-b_{j}.
\end{align*}
This proves the second part of this lemma. 
\end{proof}
First, the assumption that $G\left(a_{j}+\right)=G\left(a_{j}-\right)$
just says that no donations or returns happen at exactly the same
time as the arrivals of the bail requests, which is true in practice
and also true with probability-1 since the arrival times of donations
and bail requests are independent and continuous.

Additionally, this lemma 
 \ref{lem:shift} just says
that when the shift, $\inf_{s\leq t}\left\{ 0,M\left(s\right)\right\} $,
changes in the Skorokhod map, it changes by the amount $M^{*}\left(a_{j}-\right)-b_{j}$,
which is exactly the amount that is partially blocked by the partial
fulfillment process.

Using this idea, we show our main result in Theorem
\ref{thm:partial}.

\begin{restatable}[label={thm:partial}]{theorem}{partialLemma}
 (Interpretability of the Skorokhod Map) Consider the simple bail process $M\left(t\right)$
in Definition \ref{lem:shift}. Assume that $G\left(a_{j}-\right)=G\left(a_{j}+\right)$
for all the bail request arrival times $a_{j}$. Then, 
\begin{gather*}
M^{*}\left(t\right)=G\left(t\right)-\sum_{j=1}^{N^{b}\left(t\right)}b_{j}\wedge M^{*}\left(a_{j}-\right),
\end{gather*}
where $M^{*}\left(t\right):=\phi\left[M\right]\left(t\right)$ is
the output process of the Skorokhod map on $M\left(t\right)$. 
\end{restatable}

\begin{proof}
From Lemma \ref{lem:shift}, we know that the Skorokhod shift, $\inf_{s\leq t}\left\{ 0,M\left(s\right)\right\} $,
only changes at arrival time $t=a_{i}$ that satisfy $M^{*}\left(a_{i}-\right)<b_{i}$.
We denote subsequence of $\left\{ a_{j}\right\} _{j=1}^{\infty}$
satisfying this condition as $\left\{ a_{j_{k}}\right\} _{k=1}^{\infty}$:

\begin{gather*}
\left\{ a_{j_{k}}\right\} _{k=1}^{\infty}:=\left\{ a_{j}:\ M^{*}\left(a_{j}-\right)<b_{j}\right\} .
\end{gather*}
We adopt the convention $a_{0}=a_{j_{0}}=0$. We also use the notation
\begin{gather*}
M_{Q}\left(t\right):=G\left(t\right)-\sum_{j=1}^{N^{b}\left(t\right)}b_{j}\wedge M^{*}\left(a_{j}-\right).
\end{gather*}
We use induction here and want to prove that $M^{*}\left(t\right)=M_{Q}\left(t\right)$
for all $t\in\left[a_{j_{k}}+,a_{j_{k+1}}-\right]$ for all $k\geq0$.

We first consider the base case where $k=0$ and $t\in\left[a_{j_{0}}+,a_{j_{1}}-\right]=\left[0,a_{j_{1}}-\right]$.
By Lemma \ref{lem:shift}, $\inf_{s\leq a_{j_{1}}-}\left\{ 0,M\left(s\right)\right\} =0$
during this interval, so 
\begin{gather*}
M^{*}\left(t\right)=M\left(t\right)\ \forall\ t\in\left[0,a_{i_{1}}-\right].
\end{gather*}
Also by Lemma \ref{lem:shift}, $b_{j}\leq M^{*}\left(a_{j}-\right)$
for all $j<j_{1}$. So, for $t\in\left[0,a_{j_{1}}-\right]$, 
\begin{align*}
M_{Q}\left(t\right) & =G\left(t\right)-\sum_{j=1}^{N^{b}\left(t\right)}b_{j}\wedge M^{*}\left(a_{j}-\right)\\
 & =G\left(t\right)-\sum_{j=1}^{N^{b}\left(t\right)}b_{j}\\
 & =M\left(t\right).
\end{align*}
So, $M^{*}\left(t\right)=M\left(t\right)=M_{Q}$ for all $t\in\left[a_{j_{0}}+,a_{j_{1}}-\right]$,
and we have proved our base case.

Now, we use induction. Assume $M^{*}\left(t\right)=M_{P}\left(t\right)$
for all $t\in\left[a_{j_{k}}+,a_{j_{k+1}}-\right]$ for $k=l$. We
now prove it is true for $k=l+1$. In this case, we know that by induction
hypothesis
\begin{gather*}
M^{*}\left(a_{j_{l+1}}-\right)=M_{Q}\left(a_{j_{l+1}}-\right).
\end{gather*}
For $t\in\left[a_{j_{l+1}}+,a_{j_{l+2}}-\right]$, we rewrite 
\begin{align*}
M^{*}\left(t\right) & =\left[M^{*}\left(t\right)-M^{*}\left(a_{j_{l+1}}+\right)\right]+\left[M^{*}\left(a_{j_{l+1}}+\right)-M^{*}\left(a_{j_{l+1}}-\right)\right]+M^{*}\left(a_{j_{l+1}}-\right)\\
 & =\underbrace{\left[M^{*}\left(t\right)-M^{*}\left(a_{j_{l+1}}+\right)\right]}_{(A)}+\underbrace{\left[M^{*}\left(a_{j_{l+1}}+\right)-M^{*}\left(a_{j_{l+1}}-\right)\right]}_{(B)}+M_{P}\left(a_{j_{l+1}}-\right).
\end{align*}
We consider the $\left(A\right)$ and $\left(B\right)$ terms. For
$\left(B\right)$, 
\begin{align*}
(B) & =\left[M\left(a_{j_{l+1}}+\right)-\inf_{s\leq a_{j_{l+1}}+}\left\{ 0,M\left(s\right)\right\} \right]-\left[M\left(a_{j_{l+1}}-\right)-\inf_{s\leq a_{j_{l+1}}-}\left\{ 0,M\left(s\right)\right\} \right]\\
 & =-b_{j_{l+1}}+\inf_{s\leq a_{j_{l+1}}-}\left\{ 0,M\left(s\right)\right\} -\inf_{s\leq a_{j_{l+1}}+}\left\{ 0,M\left(s\right)\right\} \\
 & =-b_{j_{l+1}}+\left(b_{j_{l+1}}-M^{*}\left(a_{j_{l+1}}-\right)\right)\\
 & =-M^{*}\left(a_{j_{l+1}}-\right)\\
 & =-\left[b_{j_{l+1}}\wedge M^{*}\left(a_{j_{l+1}}-\right)\right]\\
 & =M_{Q}\left(a_{j_{l+1}}+\right)-M_{Q}\left(a_{j_{l+1}}-\right),
\end{align*}
where the second equality follows from the fact that $G\left(a_{j_{l+1}}+\right)=G\left(a_{j_{l+1}}-\right)$,
the third equality follows from Lemma \ref{lem:shift} , the fifth
equality comes from the definition of the subsequence $\left\{ a_{j_{k}}\right\} $.

For $\left(A\right)$, 
\begin{align*}
\left(A\right) & =\left[M\left(t\right)-\inf_{s\leq t}\left\{ 0,M\left(s\right)\right\} \right]-\left[M\left(a_{j_{l+1}}+\right)-\inf_{s\leq a_{j_{l+1}}+}\left\{ 0,M\left(s\right)\right\} \right]\\
 & =M\left(t\right)-M\left(a_{j_{l+1}}+\right)\\
 & =\left[G\left(t\right)-\sum_{j=1}^{N^{b}\left(t\right)}b_{j}\right]-\left[G\left(a_{j_{l+1}}+\right)-\sum_{j=1}^{N^{b}\left(a_{j_{l+1}}+\right)}b_{j}\right]\\
 & =G\left(t\right)-G\left(a_{j_{l+1}}+\right)-\sum_{j=j_{l+1}+1}^{N^{b}\left(t\right)}b_{j}.\\
 & =G\left(t\right)-G\left(a_{j_{l+1}}+\right)-\sum_{j=j_{l+1}+1}^{N^{b}\left(t\right)}b_{j}\wedge M^{*}\left(a_{j}-\right)\\
 & =\left[G\left(t\right)-\sum_{j=1}^{N^{b}\left(t\right)}b_{j}\wedge M^{*}\left(a_{j}-\right)\right]-\left[G\left(a_{j_{l+1}}+\right)-\sum_{j=1}^{N^{b}\left(a_{j_{l+1}}+\right)}b_{j}\wedge M^{*}\left(a_{j}-\right)\right]\\
 & =M_{Q}\left(t\right)-M_{Q}\left(a_{j_{l+1}}+\right).
\end{align*}
where the second equality follows from the fact that Lemma \ref{lem:shift}
, and the fifth equality follows from the fact that $t\leq a_{j_{l+2}}-$
and the definition of the subsequence $\left\{ a_{j_{k}}\right\} _{k=1}^{\infty}$.

Therefore, plugging in $\left(A\right)$ and $\left(B\right)$, we
have 
\begin{align*}
M^{*}\left(t\right) & =\left[M_{Q}\left(t\right)-M_{Q}\left(a_{j_{l+1}}+\right)\right]+\left[M_{Q}\left(a_{j_{l+1}}+\right)-M_{Q}\left(a_{j_{l+1}}-\right)\right]+M_{Q}\left(a_{j_{l+1}}-\right)\\
 & =M_{Q}\left(t\right).
\end{align*}
This completes the induction. 
\end{proof}

\begin{restatable}[label={cor:partial}]{corollary}{partialCorollary}
\label{col:partial}
(Skorokhod Map Output as Partial Fulfillment) Consider the Infinite Acceptance Model $M^{\infty,R^{\infty}}\left(t\right)$
in Equation \eqref{eq:MIR}. Then, the following result holds:
\begin{align}
M^{P^{*},R^{\infty}}\left(t\right) & =M_{0}+\sum_{i=1}^{N^{d}\left(t\right)}d_{i}-\underbrace{\sum_{j=1}^{N^{b}\left(t\right)}b_{j}\wedge M^{P^{*},R^{\infty}}\left(a_{j}-\right)}_{\text{Bail Returns Partially Fulfilled}}+\underbrace{\sum_{j=1}^{N^{b}\left(t\right)}p_{j}b_{j}\left\{ a_{j}+s_{j}\leq t\right\} }_{\text{Bail Returns from Fully Requested Amount}}\label{eq:SkrP}
\end{align}
where $M^{P^{*},R^{\infty}}\left(t\right):=\phi\left[M^{\infty,R^{\infty}}\right]\left(t\right)$
is the output process of the Skorokhod Map on $M^{\infty,R^{\infty}}\left(t\right).$
\end{restatable}

\begin{proof}
The proof is in Appendix \ref{app:partial}

\end{proof}

Observe that this Skorokhod output does not produce accurate returns because the returns do not exactly match the (potentially/partially) borrowed amount.This discrepancy occurs because the return process from the Skorokhod mapping directly inherits  input's infinite acceptance return process. Consequently, the returns exceed the actual borrowed amount. We will analyze this discrepancy and create a more realistic version of this model in the next section.
\subsection{\label{subsec:SkrReturns} Skorokhod Ignores Returns}

This section clarifies that the partial fulfillment
described in Corollary \ref{col:partial} 
diverges from that observed in
practical operations. To formalize this distinction, we first define
a notion of \textbf{Partial Fulfillment} consistent with real-world operations
in Equation \ref{eq:Partial}.

\begin{equation}\label{eq:Partial}
M^{P,R^{P}}\left(t\right)=M_{0}+\sum_{i=1}^{N^{d}(t)}d_{i}-\underbrace{\sum_{j=1}^{N^{b}(t)}b_{j}\wedge M^{P,R^{P}}\left(t\right)}_{\text{Bail Requests Partially Fulfilled}}+\underbrace{\sum_{j=1}^{N^{b}(t)}\left(1-p_{j}\right)\left(b_{j}\wedge M^{P,R^{P}}\left(t\right)\right)\left\{ t>a_{j}+s_{j}\right\} }_{\text{Bail Returns from Partially Fulfilled Amounts}}.
\end{equation}

This partial fulfillment is a realistic one as the
exact partially fulfilled amount is returned precisely at its designated
time.

Upon comparing the Equation \eqref{eq:SkrP} with
the Realistic Partial Fulfillment , we observe that the only distinction
appears in the bail returns component. That is, the returns from the
Skorokhod Output (left hand side) are over estimated compared to the
returns from the Realistic Partial Fulfillment:

\[
\underbrace{\sum_{j=1}^{N^{b}\left(t\right)}p_{j}b_{j}\left\{ a_{j}+s_{j}\leq t\right\} }_{\text{Bail Returns from Fully Requested Amount}}\leq\underbrace{\sum_{j=1}^{N^{b}(t)}\left(1-p_{j}\right)\left(b_{j}\wedge M^{P,R^{P}}\left(t\right)\right)\left\{ t>a_{j}+s_{j}\right\} }_{\text{Bail Returns from Partially Fulfilled Amounts}}
\]

\begin{figure}[H]
\centering
\includegraphics[scale=0.6]{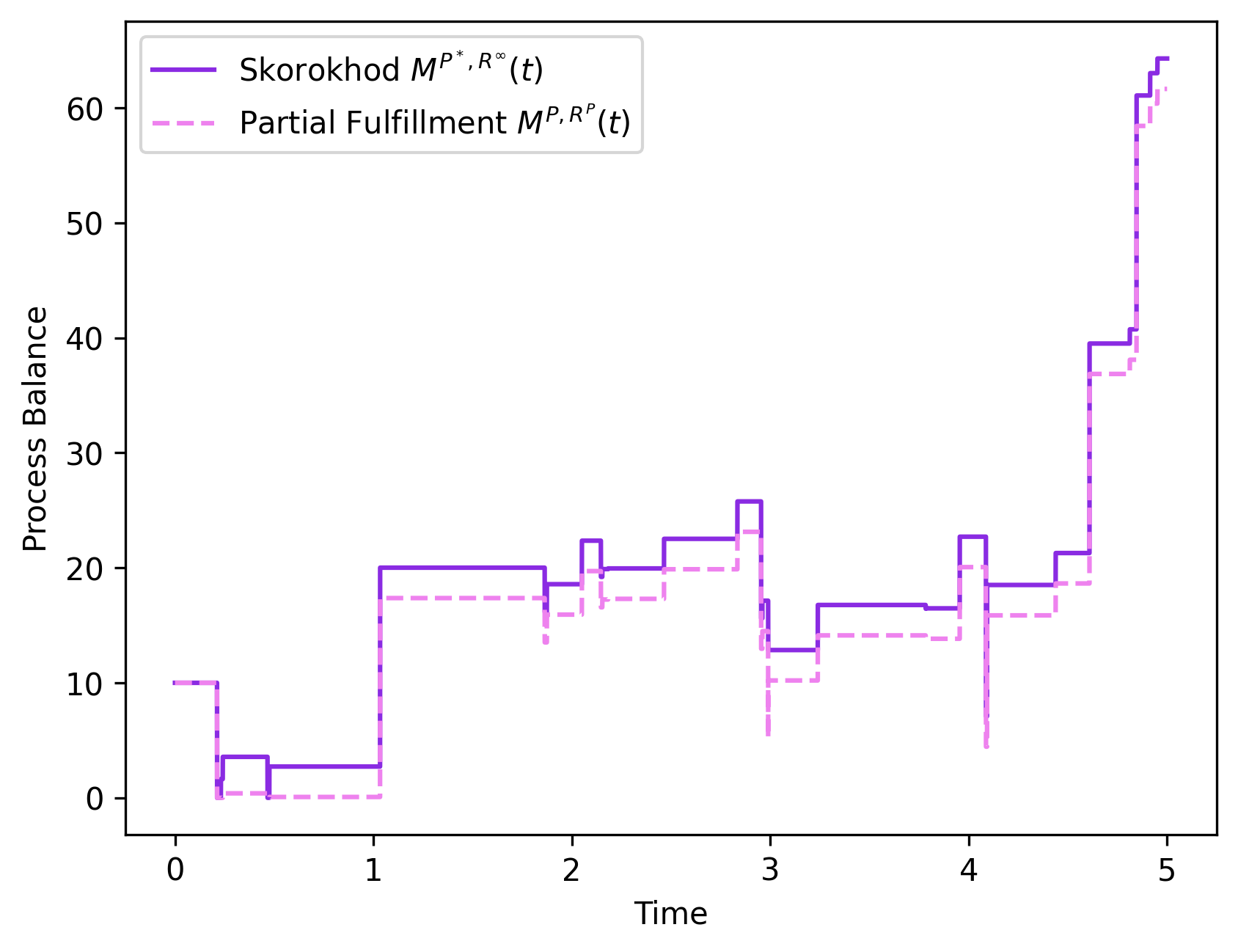}
\caption{Skorokhod Approximation is not a Partial Fulfillment}\label{fig:Skrk_VA_Partl}
\end{figure}

We include a simulation  in Figure \ref{fig:Skrk_VA_Partl} to compare the overall
differences accumulated between two processes as a result of naive/over
estimations in the return components.

This naive/over estimation of bail returns in the
Skorokhod output is due to the property of Skorokhod map itself. Theorem
\ref{thm:partial} shows while the Skorokhod map changes the bail
fulfillment component
of the simple processes to reflect partial fulfillment, it does \textit{not
} change anything in the non-decreasing component
$G\left(t\right)$ including the returns.

\begin{figure}[H]
\centering
\includegraphics[scale=0.6]
{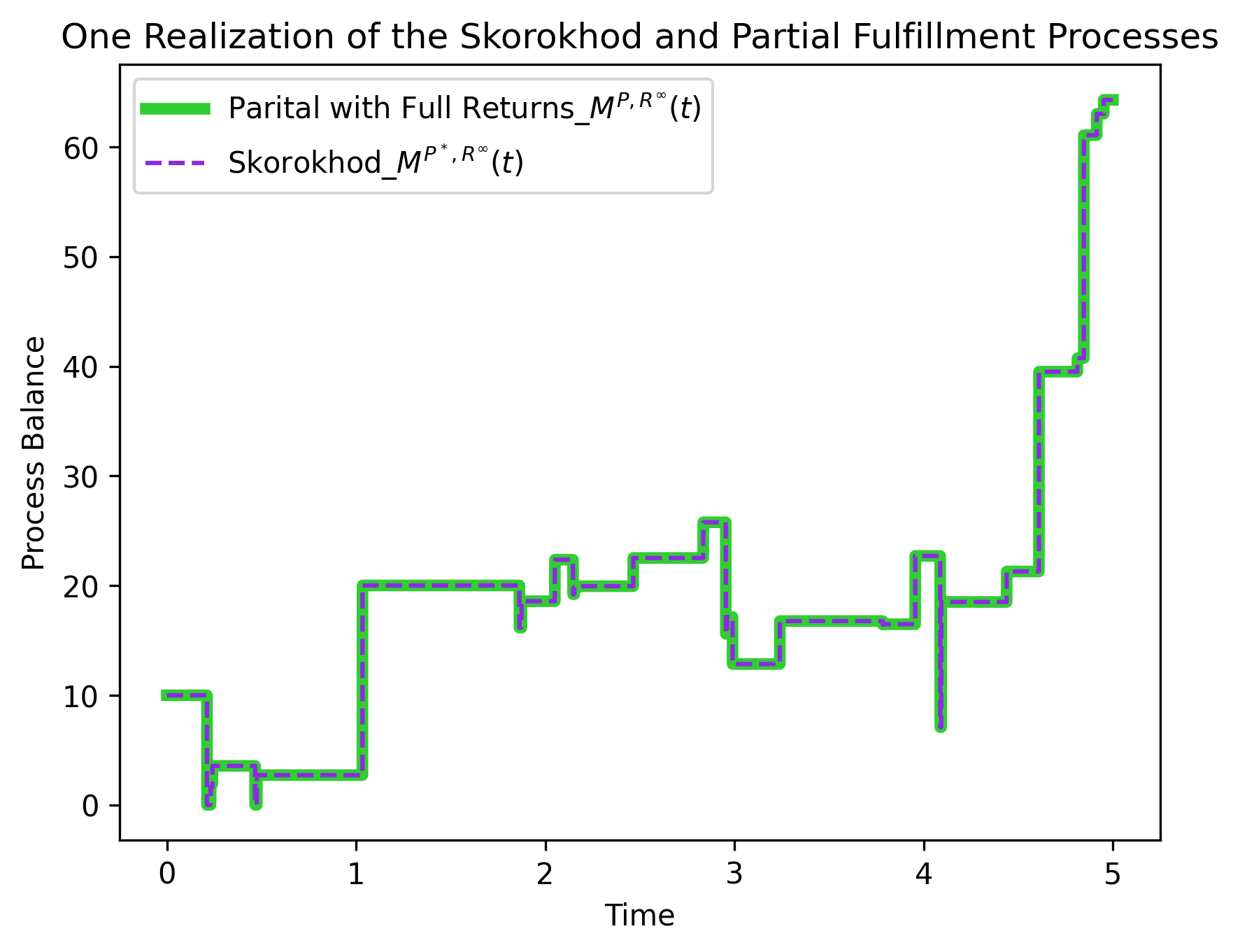}
\caption{Partial with Full Returns v.s. Skorokhod Approximation}
\label{fig:SkrkVSPlinf}
\end{figure}

Our subsequent simulation experiments in Figure  \ref{fig:SkrkVSPlinf} further verify
this. We construct a hybrid model by modifying the components as follows:
We take the Realistic Partial Fulfillment model, remove its return
component, and instead incorporate the return component from the Infinite
Acceptance Model in its place. This constructed hybrid model (solid line)
matches exactly with the approximation generated by applying the Skorokhod
Map on the Infinite Acceptance Model (dashed line).

Interestingly, if a fund process lacks a return
component entirely, the Realistic Partial Fulfillment Model without
return component will produce paths identical to those generated by
the Skorokhod Approximation applied to the Infinite Acceptance Model
without its return component. In this case, the two models are mathematically
equivalent. See the simulation result in Figure \ref{fig:SkrkNRvsPNR} below:

\begin{figure}[H]
\centering
\includegraphics[scale=0.6]{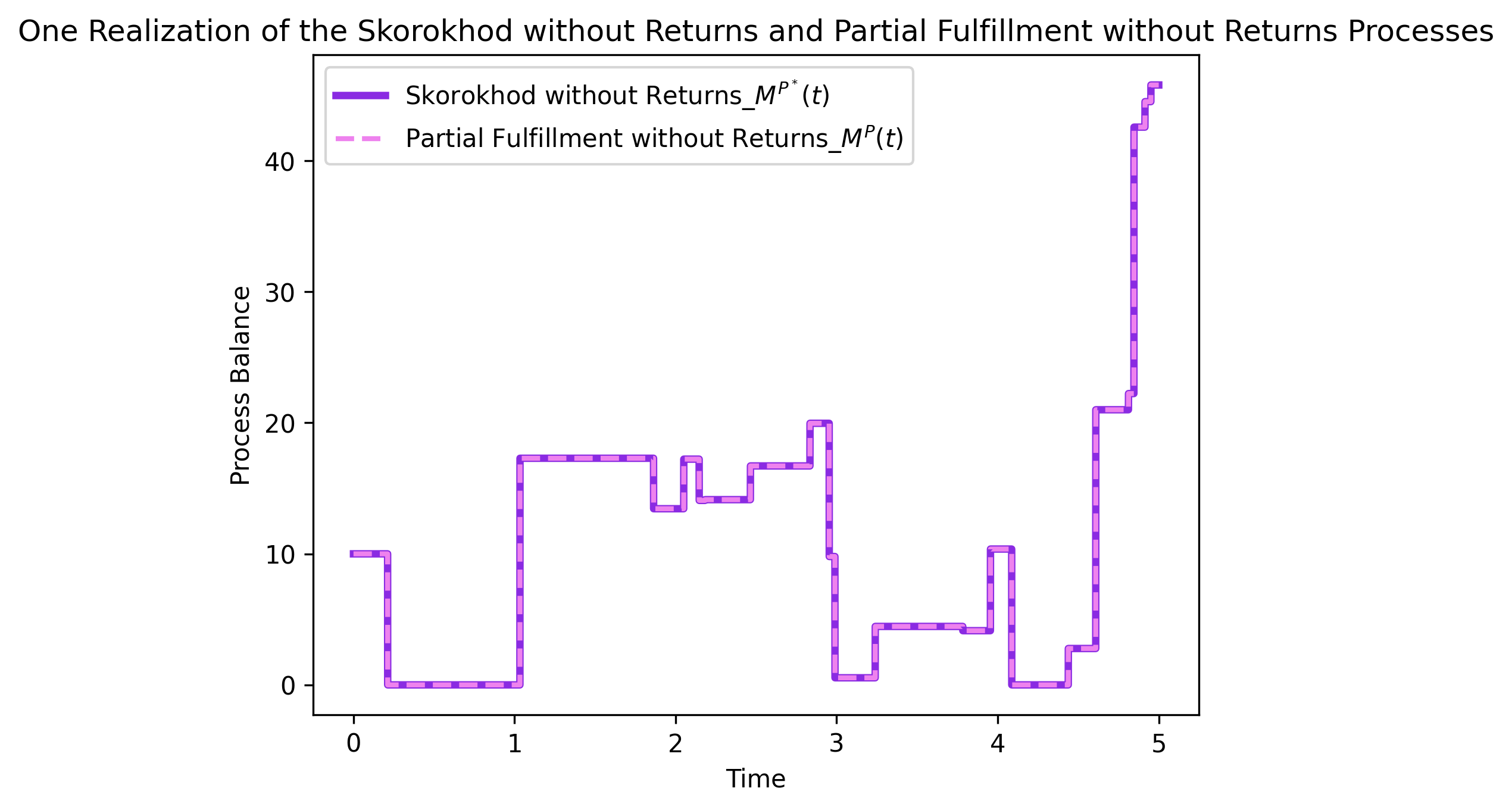}
\caption{Partial Fulfillment without Returns v.s. Skorokhod Approximation without returns}\label{fig:SkrkNRvsPNR}
\end{figure}

\subsection{\label{subsec:SkrkLimit}Skorokhod Limit}

\begin{restatable}[label={thm:SkrkLimit}]{theorem}{SkrkLimitTheorem}
\label{thm:SkrkLimit}
(Skorokhod Limit) The Skorokhod map of the Infinite Acceptance Process
approaches the following limit as $\eta\rightarrow\infty$:

{\footnotesize
\begin{gather*}
\phi\left[M_{\eta}^{\infty,R^{\infty}}\right]\left(t\right)\underset{a.s.}{\rightarrow}\phi\left[m^{\infty,R^{\infty}}\right]\text{\ensuremath{\left(t\right)}}\\
=\begin{cases}
m^{\infty,R^{\infty}}\left(t\right) & \text{if }d^{*}\lambda_{d}-b^{*}\lambda_{b}\ge0,\\
m^{\infty,R^{\infty}}\left(t\right)-\min\left\{ 0,m^{\infty,R^{\infty}}\left(t\right)\right\}  & \text{if }d^{*}\lambda_{d}+\left(1-p^{*}\right)b^{*}\lambda_{b}F_{s}\left(z\right)\leq b^{*}\lambda_{b}\\
m^{\infty,R^{\infty}}\left(t\right)-\inf_{z\leq t}\left\{ 0,m^{\infty,R^{\infty}}\left(z\right)\right\}  & \text{if }\ensuremath{b^{*}\lambda_{b}-\left(1-p^{*}\right)b^{*}\lambda_{b}\leq d^{*}\lambda_{d}\leq b^{*}\lambda_{b}.}
\end{cases}\text{for }\forall z\leq t.
\end{gather*}
}{\footnotesize\par}
\end{restatable}

\begin{proof}
Since the Skorokhod map $\phi\left[\text{\ensuremath{\cdot}}\right]$
is a continuous mapping on $\left(\mathbb{D}\left[0,\infty\right],\mathbb{R}\right)$
,

\begin{align*}
\lim_{\eta\rightarrow\infty}\phi\left[M_{\eta}^{\infty,R^{\infty}}\right]\left(t\right) & =\phi\left[\lim_{\eta\rightarrow\infty}M_{\eta}^{\infty,R^{\infty}}\right]\left(t\right)\\
 & =\phi\left[m^{\infty,R^{\infty}}\right]\left(t\right)\\
 & =m^{\infty,R^{\infty}}\left(t\right)-\inf_{z\leq t}\left\{ 0,m^{\infty,R^{\infty}}\left(z\right)\right\} \\
\end{align*}
Recall that the time derivative of the fluid limit of the infinite
acceptance process for $\forall z\leq t$ is given by:

\begin{align*}
\frac{d}{dz}m^{\infty,R^{\infty}}\left(z\right) & =d^{*}\lambda_{d}-b^{*}\lambda_{b}+\left(1-p^{*}\right)b^{*}\lambda_{b}F_{s}\left(z\right)\text{ for \ensuremath{\forall z\le t}}
\end{align*}
Observe $\left(1-p^{*}\right)b^{*}\lambda_{b}F_{s}\left(z\right)$
is monotonically increasing and can be bounded above
and below by constants, as shown below:
\[
0\leq\left(1-p^{*}\right)b^{*}\lambda_{b}F_{s}\left(z\right)\leq\left(1-p^{*}\right)b^{*}\lambda_{b}.
\]

\paragraph{Case 1.} If $d^{*}\lambda_{d}\geq b^{*}\lambda_{b},$ then the time derivative
at time $z$:
\[
\underbrace{\left(d^{*}\lambda_{d}-b^{*}\lambda_{b}\right)}_{\geq0}+\underbrace{\left(1-p^{*}\right)b^{*}\lambda_{b}F_{s}\left(z\right)}_{\geq0}\geq0\quad\forall z\leq t
\]
So this fluid limit is always increasing and never becomes negative.
Thus, Skorokhod shift $\inf_{z\leq t}\left\{ 0,m^{\infty,R^{\infty}}\left(z\right)\right\} $
will always be 0, and we have
\begin{align*}
\lim_{\eta\rightarrow\infty}\phi\left[M_{\eta}^{\infty,R^{\infty}}\right]\left(t\right) & =m^{\infty,R^{\infty}}\left(t\right)-\inf_{z\leq t}\left\{ 0,m^{\infty,R^{\infty}}\left(z\right)\right\} \\
 & =m^{\infty,R^{\infty}}\left(t\right).
\end{align*}
\paragraph{Case 2.}If $d^{*}\lambda_{d}+\left(1-p^{*}\right)b^{*}\lambda_{b}F_{s}\left(z\right)\leq b^{*}\lambda_{b}\text{ for }\forall z\leq t,$
then the time derivative 

\[
d^{*}\lambda_{d}-b^{*}\lambda_{b}+\left(1-p^{*}\right)b^{*}\lambda_{b}F_{s}\left(z\right)\leq0\quad\forall z.
\]
This inequality indicates that the combined average
rate at which donations and returns flow into the fund per unit time
is at most equal to the average rate at which bail amounts flow out.
Thus, the fluid limit is monotonically decreasing and becomes negative
for some large enough time. Then, the Skorokhod shift is described
by

\[
\inf_{z\leq t}\left\{ 0,m^{\infty,R^{\infty}}\left(z\right)\right\} =\min\left\{ 0,m^{\infty,R^{\infty}}\left(t\right)\right\} .
\]

We have,
\begin{align*}
\lim_{\eta\rightarrow\infty}\phi\left[M_{\eta}^{\infty,R^{\infty}}\right]\left(t\right) & =m^{\infty,R^{\infty}}\left(t\right)-\inf_{z\leq t}\left\{ 0,m^{\infty,R^{\infty}}\left(z\right)\right\} \\
 & =m^{\infty,R^{\infty}}\left(t\right)-\min\left\{ 0,m^{\infty,R^{\infty}}\left(t\right)\right\} 
\end{align*}

\paragraph{Case 3.} If $b^{*}\lambda_{b}-\left(1-p^{*}\right)b^{*}\lambda_{b}\leq d^{*}\lambda_{d}\leq b^{*}\lambda_{b}$,
then the sign of the time derivitive depends on the current value
of $F_{s}\left(z\right)$. Observe that this inequality follows directly
from the set of inequalities:
\begin{gather*}
d^{*}\lambda_{d}\leq b^{*}\lambda_{b},\\
0\leq\left(1-p^{*}\right)b^{*}\lambda_{b}F_{s}\left(z\right)\leq\left(1-p^{*}\right)b^{*}\lambda_{b},\\
b^{*}\lambda_{b}-d^{*}\lambda_{d}\leq\left(1-p^{*}\right)b^{*}\lambda_{b}.
\end{gather*}The first line states that average rate of donations
flowing into the fund is at most equal to the average rate of bail
amounts flowing out. The second line states that the average rate
of returns flowing into the fund is a non-negative quantity dependent
on the current time value $z,$ and it can be upper bounded by a constant
$\left(1-p^{*}\right)b^{*}\lambda_{b}.$ The third line states that
the deficit between bails outflows and donations inflows can be fully
compensated by this upper bound on returns inflow. Consequently, the
current value of $F_{s}\left(z\right)$ determines the sign of the
time derivative. We then can only use the original definition of
the Skorokhod map.

Therefore,

{\small
\begin{align*}
\lim_{\eta\rightarrow\infty}\phi\left[M_{\eta}^{\infty,R^{\infty}}\right]\left(t\right) & =m^{\infty,R^{\infty}}\left(t\right)-\inf_{z\leq t}\left\{ 0,m^{\infty,R^{\infty}}\left(z\right)\right\} \\
 & =\begin{cases}
m^{\infty,R^{\infty}}\left(t\right) & \text{if }d^{*}\lambda_{d}-b^{*}\lambda_{b}\ge0,\\
m^{\infty,R^{\infty}}\left(t\right)-\min\left\{ 0,m^{\infty,R^{\infty}}\left(t\right)\right\}  & \text{if }d^{*}\lambda_{d}+\left(1-p^{*}\right)b^{*}\lambda_{b}F_{s}\left(z\right)\leq b^{*}\lambda_{b}\\
m^{\infty,R^{\infty}}\left(t\right)-\inf_{z\leq t}\left\{ 0,m^{\infty,R^{\infty}}\left(z\right)\right\}  & \text{if }\ensuremath{b^{*}\lambda_{b}-\left(1-p^{*}\right)b^{*}\lambda_{b}\leq d^{*}\lambda_{d}\leq b^{*}\lambda_{b}.}
\end{cases}\text{for }\forall z\leq t
\end{align*}
}{\small\par}

This completes the proof.
\end{proof}

\section{Bail Fund Models with Blocking\label{sec:Models}}

\subsection{\textcolor{black}{\label{subsec:BRModels}Modeling of the Bail Fund
with Blocking}}

In this section, we define a blocking process to directly model
the limited funding constraints faced by the bail fund. In this model,
the decision rules are as follows: if, at the moment just before a
bail request arrives, the fund has a sufficient balance to cover the
bail amount, then the request is accepted and fulfilled. If the balance
is insufficient at that moment, then the bail request is blocked,
and no money flows out from the bail fund. Formally, the \textbf{blocking
process} is defined as:

\begin{align}
\underbrace{M^{B,R}\left(t\right)}_{\text{Total Funds Available}} & =\underbrace{M_{0}}_{\text{Initial Capital}}+\underbrace{\sum_{i=1}^{N^{d}\left(t\right)}d_{i}}_{\text{Total Donations}}-\underbrace{\sum_{j=1}^{N^{b}\left(t\right)}b_{j}\cdot\left\{ M^{B,R}\left(a_{j}-\right)\geq b_{j}\right\} }_{\text{Total Bail Paid Out}}\nonumber \\
 & \quad+\underbrace{\sum_{j=1}^{N^{b}(t)}\left(1-p_{j}\right)\cdot b_{j}\cdot\left\{ M^{B,R}\left(a_{j}-\right)\geq b_{j}\right\} \cdot\left\{ t>a_{j}+s_{j}\right\} }_{\text{Total Fund Returned}}.\label{eq:MBR}
\end{align}
The first two terms on the right-hand side correspond to the initial
funds and donations as in the infinite acceptance process. The third
term indicates the total bail amount paid for defendants awaiting
trial up to time $t$, which is conditional on whether enough funds
were available immediately before the $j^{\text{th}}$ bail request
arrived, using the indicator function $\left\{ M^{B,R}\left(a_{j}-\right)\geq b_{j}\right\} $.
The final term accounts for the total bail returned by the courts
by time $t$, excluding forfeitures, and is similarly conditioned
on whether the bail requests were initially approved.

Notice that in each blocking term $\left\{ M^{B,R}\left(a_{j}-\right)\geq b_{j}\right\} $
of this model, the balance of the fund itself evaluated just before
the arrival time of $j^{th}$ bail request, $M^{B,R}\left(a_{j}-\right)$,
is also considered. This indicates a balance dependence in the evolving
process, as the value of $M^{B,R}\left(a_{j}-\right)$ depends on
all preceding donations $\left\{ d_{i}\right\} _{i\geq1}$, all preceding
bail requests $\left\{ b_{j}\right\} _{j\geq1}$, all time-to-trial
values $\left\{ s_{j}\right\} _{j\geq1}$, and all poundage fees $\left\{ p_{j}\right\} _{j=1}^{N^{b}\left(a_{j}-\right)}$
before time $t=a_{j}$.

Deriving the expectation or variance of the blocking
process is challenging due to its dependent jumps tied to this blocking
term. This balance dependence introduces substantial analytical complexity,
in contrast to the Infinite Acceptance Process, where such derivations
are straightforward. However, leveraging the analytical tools available,
we establish a fluid limit theorem for this process to analyze its
average behavior. As a preliminary step, we construct a scaling model
in the next section.

\subsection{\textcolor{black}{\label{subsec:BRScaling}Scaling of the Bail Fund
Model}}

Now we scale our blocking process by the same scaling
spirit as we do for Infinite Acceptance Model in Section \ref{subsec:InfScaling}.
Define the \textbf{$\eta$-scaled version of $M^{B,R}\left(t\right)$}
from Equation \ref{eq:MBR}:

\begin{align*}
M^{B,R}_\eta\left(t\right) & =\frac{1}{\eta}\left[\eta M_{0}+\sum_{i=1}^{N_{\eta}^{(d)}(t)}d_{i}-\sum_{j=1}^{N_{\eta}^{b}(t)}b_{j}\cdot\left\{ M^{B,R}_{\eta}(a_{j,\eta}-)\geq b_{j}\right\} \right]\\
 & \quad+\frac{1}{\eta}\left[ \sum_{j=1}^{N_{\eta}^{(b)}( t)}(1-p_{j})\cdot b_{j}\cdot\left\{ M^{B,R}_{\eta}(a_{j,\eta}-)\geq b_{j}\right\} \cdot\left\{ t>a_{j,\eta}+s_{j}\right\}\right].
\end{align*}
 
In the next section, we prove this scaled version
of the blocking model converges to its fluid limit.

\subsection{\textcolor{black}{\label{subsec:FLConvergeBlk}Convergence to Fluid
Limit}}

\textcolor{black}{Our goal in this Section \ref{subsec:FLConvergeBlk}
is to prove the convergence of our models defined in Section \ref{subsec:BRScaling}
to their fluid limits. To do this, we will first prove the following
steps: }
\begin{enumerate}
\item \textcolor{black}{The donation component of models converges almost
surely as in the Infinite Acceptance Model (Section 
 \ref{subsubsec:CPConverge}).}
\item The blocking component and the return components
of model converge in probability (Section \ref{subsubsec:BRConverge}).
\begin{enumerate}
\item Prove that the bail request and return components
are martingales.
\item Derive an upper bound for the second moment of each
component to deal with blocking terms.
\item Apply Doob's supremum inequality for continuous
submartingales and use the bounded second moments to prove convergence
in probability in Lemma \ref{lem:ProbConvergence}.
\end{enumerate}

\end{enumerate}
\textcolor{black}{Combined with these above 2 steps, we prove the
models converge to their fluid limits in Section \ref{subsubsec:FLConverge}. }

\subsubsection{\label{subsubsec:BRConverge}General Bounds and
Convergence in Probability for the Bail Request and Return Components}
We now establish the convergence in probability of both the bail and return
components.

\begin{restatable}[label={lem:ProbConvergence}]{lemma}{ProbConvergenceLemma}
\label{lem:ProbConvergence}
(Generalized Martingale Convergence) Consider
any arbitrary $\eta$- dependent process $M_{\eta}\left(t\right)$,
and assume
\begin{itemize}
\item $M_{\eta}\left(t\right)$ is a martingale, and 
\item $\E\left[\left(M_{\eta}\left(t\right)\right)^{2}\right]\leq\frac{C}{\eta}$
for any $0\leq t\leq T$ and some constant $C\in\mathbb{R}$,
\end{itemize}
where $T$ is any fixed end-time. Then, as $\eta\rightarrow\infty$,
\[
\sup_{0\leq u\leq T}\left|M_{\eta}\left(u\right)\right|\overset{P}{\rightarrow}0.
\]
\end{restatable}
\begin{proof}
This proof is in Appendix \ref{app:ProbConvergence}
\end{proof}

The function $H:\mathbb{\mathbb{R^{\text{+}}}\rightarrow\mathbb{R^{\text{+}}}}$
represents the truncated expected value of the bail request jump
sizes $b_{j}$ :

\begin{equation}
H\left(m\right):=\int_{0}^{m}xf_{b}\left(x\right)dx,\label{eq:Hfcn}
\end{equation}
where $f_{b}\left(\cdot\right)$ is the probability density function
of $b_{j}$.

This term will appear in many of the subsequent
theorems.

\begin{restatable}[label={lem:BCP}]{lemma}{BCPLemma}
\label{lem:BCP}(Centered Bail Component Convergence in Probability) The centered process 

\begin{align*}
\M_{\eta,b}^{B,R}\left(t\right) & :=\frac{1}{\eta}\sum_{j=1}^{N^{b}_{\eta}\left( t\right)}b_{j}\cdot\left\{ M_{\eta}^{B,R}\left(a_{j}-\right)\geq b_{j}\right\} -\lambda\int_{0}^{t}H\left(M_{\eta}^{B,R}\left(u\right)\right)du
\end{align*}
 is a martingale with respect to the filtration $\left\{ \F_{s}\right\} $
which contains all bail arrivals and bail request sizes until time
$s$:
\[
\F_{s}:=\sigma\left(N^{b}\left(\eta s\right),\left\{ b_{j}\right\} _{j=1}^{N^{b}\left(\eta s\right)},\left\{ a_{j}\right\} _{j=1}^{N^{b}\left(\eta s\right)}\right).
\]

Additionally, as $\eta\rightarrow\infty,$

\[
\sup_{t\in[0,T]}\left|\M_{\eta,b}^{B,R}\left(t\right)\right|\overset{P}{\rightarrow}0.
\]
\end{restatable}
\begin{proof}
This proof is in Appendix \ref{app:BCP}.
\end{proof}

\begin{restatable}[label={lem:RCP}]{lemma}{RCPLemma}
\label{lem:RCP}
(Centered Return Component Convergence in Probability) The centered
process 
\begin{align*}
\M_{\eta,r}^{B,R}\left(t\right) & :=\frac{1}{\eta}\sum_{j=1}^{N_{\eta}^{b}\left( t\right)}\left(1-p_{j}\right)\cdot b_{j}\cdot\left\{ M_{\eta}^{B,R}\left(a_{j}-\right)\geq b_{j}\right\} \cdot\left\{ t>a_{j,\eta}+s_{j}\right\} \\
 & \quad-\left(1-p^{*}\right)\cdot\lambda_{b}\cdot\int_{0}^{ t}H\left(M_{\eta}^{B,R}\left(u\right)\right)F_{s}\left( t-u\right)du
\end{align*}
 is a martingale with respect to the filtration $\left\{ \F_{s}\right\} $
which contains all bail arrivals and bail request sizes until time
$s$:
\[
\F_{s}:=\sigma\left(N_{\eta}^{b}\left( s\right),\left\{ b_{j}\right\} _{j=1}^{N_{\eta}^{b}\left( s\right)},\left\{ a_{j}\right\} _{j=1}^{N_{\eta}^{b}\left( s\right)},\left\{ s_{j}\right\} _{j=1}^{N_{\eta}^{b}\left( s\right)}\right).
\]

Additionally, as $\eta\rightarrow\infty,$

\[
\sup_{t\in[0,T]}\left|\M_{\eta,r}^{B,R}\left(t\right)\right|\overset{P}{\rightarrow}0.
\]
\end{restatable}
\begin{proof}
This proof is in Appendix \ref{app:RCP}.
\end{proof}

\subsubsection{\label{subsubsec:FLConverge}Convergence of $\eta$-scaled
Model to the Fluid Limit}

We first prove function $H_{b}\left(\cdot\right)$
is Lipschitz which will help us prove Theorem \ref{thm:BRFL}.

\begin{lemma}
Consider the function 

\[
H_{b}\left(m\right):=\int_{0}^{m}xf_{b}\left(x\right)dx.
\]

$H_{b}\left(\cdot\right)$ is a Lipschitz function for distribution
functions $f_{b}\left(x\right)$ that satisfy

\[
\sup_{m\in R^{\text{+}}}mf_{b}\left(m\right)<\infty.
\]
\end{lemma}
\begin{proof}

We need only verify the difference 
\begin{gather*}
H_{b}\left(m_{1}\right)-H_{b}\left(m_{2}\right),
\end{gather*}
 for $m_{1},m_{2}\in\mathbb{R}^{+}$with $m_{1}<m_{2}$, is bounded
on $\mathbb{R}^{+}$ for the Lipschitz continuity :

\begin{align*}
H_{b}\left(m_{1}\right)-H_{b}\left(m_{2}\right) & =\int_{0}^{m_{1}}xf_{b}\left(x\right)dx-\int_{0}^{m_{2}}xf_{b}\left(x\right)dx\\
 & =\int_{m_{2}}^{m_{1}}xf_{b}\left(x\right)dx\\
 & \leq\left(m_{1}-m_{2}\right)\sup_{m\in R^{\text{+}}}mf_{b}\left(m\right)\\
 & \leq C\left(m_{1}-m_{2}\right).
\end{align*}
Therefore, $H_{b}\left(\cdot\right)$ is lipschitz whenever the distribution
function of $b_{j}$ satisfies

\[
\sup_{m\in R^{\text{+}}}mf_{b}\left(m\right)<\infty.
\]
\end{proof}
\paragraph{Remark.}
It can be shown that the condition $\sup_{m\in R^{\text{+}}}mf_{b}\left(m\right)<\infty$
can be satisfied by most standard distributions. For
example,when $b_{j}\sim\text{Expo(\ensuremath{\lambda_{d}})}$, then
$mf_{b}\left(m\right)$ achieves a finite maximum at $m=\lambda_{b}$,
which is the solution to $\frac{d}{dm}m\lambda_{b}e^{-\lambda_{b}m}=0.$

We finally establish the fluid limit of our blocking
process in the Theorem \ref{thm:BRFL}. Owing to the balance dependence
of the blocking term, an implicit form of fluid limit is demonstrated
in Theorem \ref{thm:BRFL}, along with in-probability convergence.

\begin{restatable}[label={thm:BRFL}]{theorem}{BRFLTheorem}
\label{thm:BRFL} (Fluid Limit of the Blocking Model) Consider the following
$m^{B,R}\left(t\right)$ process:

\begin{align*}
m^{B,R}\left(t\right) & =M_{0}+\lambda_{b}d^{*}t-\lambda_{b}\int_{0}^{t}H\left(m^{B,R}\left(u\right)\right)du+\left(1-p^{*}\right)\cdot\lambda_{b}\cdot\int_{0}^{t}H\left(m^{B,R}(u)\right)F_{s}\left(t-u\right)du,
\end{align*}
as well as the following $M^{B,R}_{\eta}\left(t\right)$ process:

\begin{align*}
M^{B,R}_{\eta}\left(t\right) & =\frac{1}{\eta}\left[\eta M_{0}+\sum_{i=1}^{N_{\eta}^{(d)}( t)}d_{i}-\sum_{j=1}^{N_{\eta}^{(b)}(t)}b_{j}\cdot\left\{ M^{B,R}_{\eta}(a_{j,\eta}-)\geq b_{j}\right\} \right]\\
 & \quad+\frac{1}{\eta}\left[\sum_{j=1}^{N_{\eta}^{(b)}( t)}(1-p_{j})\cdot b_{j}\cdot\left\{ M_{\eta}^{B,R}(a_{j,\eta}-)\geq b_{j}\right\} \cdot\left\{ t>a_{j,\eta}+s_{j}\right\} \right],
\end{align*}
Then,
\[
\sup_{t\in\left[0,T\right]}\left|M^{B,R}_{\eta}\left(t\right)-m^{B,R}\left(t\right)\right|\overset{P}{\to}0,
\]
as $\eta\rightarrow\infty$. Therefore, $m^{B,R}\left(t\right)$ is
the fluid limit of $M^{B,R}_{\eta}\left(t\right)$ as $\eta\rightarrow\infty$. Consequently, as $\eta \to \infty$, the dynamics of the limiting fluid model is described by:
\[
\updot{m}^{B,R}\left(t\right):=\frac{d}{dt}m^{B,R}\left(t\right)=\lambda_{d}d^{*}-\lambda_{b}H\left(m^{B,R}(t)\right)+\left(1-p^{*}\right)\cdot\lambda_{b}\cdot\left[\int_{0}^{t}H\left(m^{B,R}(u)\right)f_{s}(t-u)du\right].
\]

\end{restatable}
\begin{proof}
This proof is in Appendix \textbf{\ref{app:BRFL}}. 
\end{proof}

\subsection{Numerical Results}

\begin{figure}[htbp]
\includegraphics[scale=0.5]{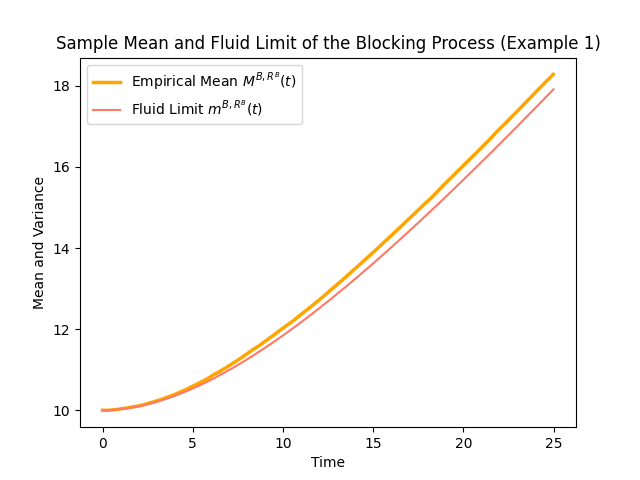}
\includegraphics[scale=0.5]{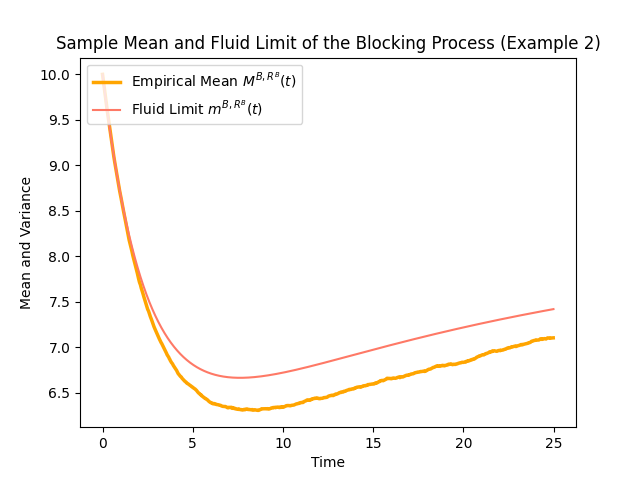}
\caption{Mean and Fluid Limit of the Blocking Process}
\begin{center}
\begin{tabular}{|c|c|c|c|c|c|c|c|}
\hline 
Example/Parameters & $M_{0}$ & $\lambda_{b}$ & $\lambda_{b}$ & $d_{i}$ & $b_{j}$ & $p_{j}$ & $s_{j}$\tabularnewline
\hline 
\hline 
Example 1 & 10 & 1 & 1 & $\text{Expo}(1)$ & $\text{Expo}(1)$ & $\text{Unif}\left[0,1\right]$ & $\text{Expo}(10)$\tabularnewline
\hline 
Example 2 & 10 & 1 & 1 & $\text{Expo}(1)$  &  $\text{Expo}(10)$ & $\text{Unif}\left[0,1\right]$ & $\text{Expo}(10)$\tabularnewline
\hline 
\end{tabular}
\end{center}
\label{fig:blk_mean_variance}
\end{figure}

First, we discuss the numerical results in Figure \ref{fig:blk_mean_variance}. In \textbf{Example 1}, fluild limit matches the sample mean. Both are initially
roughly constant and then gradually increases. Early on, donations
and bail requests yield equivalent average increments per unit time.
Around $t=10,$ matching the mean service time $s_{j}\sim\text{Expo\ensuremath{(10)}}$,
approximately half of the bail amounts start returning, proportionally
reduced by poundage rates $p_{j}\sim\text{Unif\ensuremath{\left[0,1\right]}}$.
\textbf{In Example 2}, the sample mean behaves similar to the fluid limit. In this case, we first observe a gradual drop, after which
both level off. This flattening occurs due to frequent blocking arises from the large
average bail-request jump sizes $\text{\ensuremath{\left(Expo\ensuremath{(10)}\right)}}$.

\begin{figure}[H]
\centering
\includegraphics[scale=0.5]{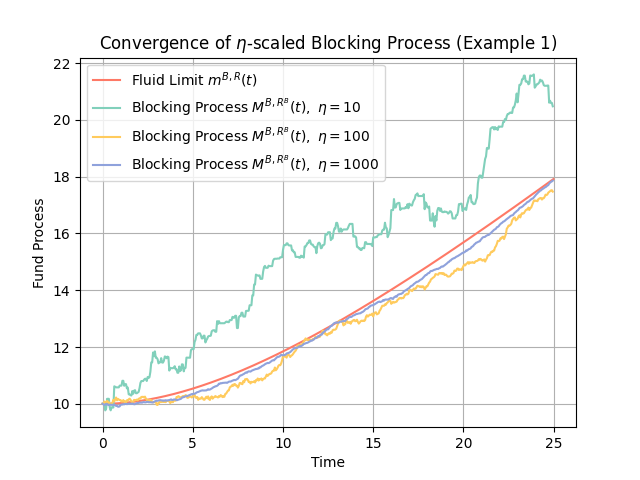}
\includegraphics[scale=0.5]{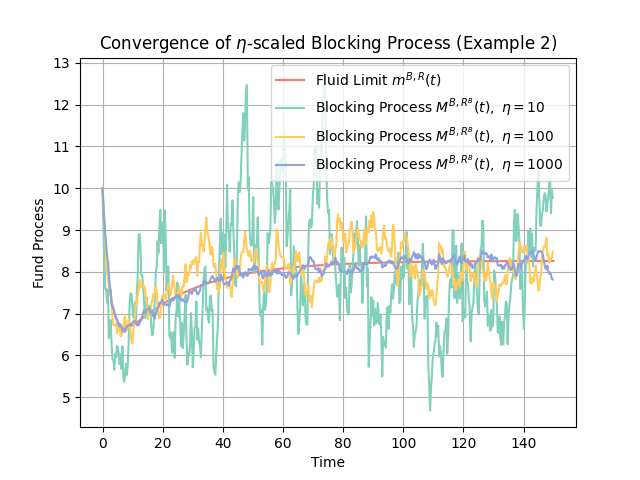}
\caption[Convergence of the Blocking Process]{
\centering Convergence of the Blocking Process\par 
\justifying
In example 2, we use a longer timescale here to show that the fluid limit $m^{B,R}(t)$ converges asymptotically to a constant amount near 8.2 in this setting.}
\label{fig:MBRFluid}
\end{figure}

Second, the simulation result in Figure \ref{fig:MBRFluid} (using the same set of parameters as the two examples above) empirically verifies that the $\eta$-scaled
blocking process $M_{\eta}^{B,R}\left(t\right)$ converges to its
fluid limit $m^{B,R}\left(t\right)$. This is mathematically proven in Theorem \ref{thm:BRFL}.

\subsection{\textcolor{black}{Insights from the Numerical Analysis}}

Observe that when the expected donation rate per
unit time, $d^{*}\lambda_{d}$ , exceeds the expected bail request
rate per unit time, $b^{*}\lambda_{b}$, the blocking process mirrors
the behavior of the Infinite Acceptance Process. In this case, the
Infinite Acceptance Process never incurs a deficit, and the Blocking
Model does not reject any bail requests. This is due to the aggregate
inflows consistently surpassing the outflows, ensuring a positive
financial equilibrium.

When the expected donation rate per unit time, $d^{*}\lambda_{d}$, is smaller than the expected bail request rate per unit time, $b^{*}\lambda_{b}$,
the blocking process diminishes over time until stabilizing at a constant
level above zero. Typically, this constant level corresponds to a
value below the average bail request amount, as it marks the point
where the blocking process must begin blocking bail requests due to
insufficient balance. Refer to the simulation result in the previous
section. We provide an example where the $\eta$-scaled blocking
process stabilizes at a value of approximately 40 as $\eta$ increases,
occurring after around $t=30.$

\section{Stochastic Ordering}

In this section, we conclude our work by establishing a stochastic
ordering as the foundational structure across four processes we have
constructed so far: the blocking $M^{B,R^{B}}\left(t\right),$ the
infinite acceptance $M^{\infty,R^{\infty}}\left(t\right),$ the Skorokhod
approximation $M^{P^*,R^{\infty}}\left(t\right),$ and the
partial fulfillment $M^{P,R^{P}}\left(t\right).$ This ordering serves
as a reference point to identify an optimal bail fund policy in the
future and build further analysis. 

Notably, the return component emerges as a key disruptor
of stochastic ordering. A simple example is listed below to explain
how returns disrupt order analysis.

\begin{restatable}{example}{ReturnsExample} \label{ex:returnsmess}
Consider an example
where returns complicate the Skorokhod approximation $M^{P^*,R^P}(t)$ and the blocking
process $M^{B,R^B}(t)$, and the ordering is not preserved, as shown in the table below. 
The infinite acceptance model $M^{\infty,R^{\infty}}(t)$ (shown in bold) is included to help explain the Skorokhod approximation. Notice how the ordering of the processes changes when the returns are added back.

\begin{center}
\begin{tabular}{|c|c|c|c|c|c|c|c|}
\hline 
\multicolumn{8}{|c|}{$s_{j}\equiv2,j=1,2$}\tabularnewline
\multicolumn{8}{|c|}{$p_{j}\equiv0,j=1,2$}\tabularnewline
\hline 
 time & $t_{0}=0$ & $t_{1}=1$ & $t_{2}=2$ & $t_{3}=3$ & $t_{4}=4$ & $t_{5}=5$ & $t_{6}=6$
\tabularnewline
\hline
events & $d_1=5$ & $b_{1}=6$ &  & $b_{2}=4$ &  $r_1$ &  & $r_2$\tabularnewline
\hline
{$\mathbf{M}^{\infty,R^{\infty}}(t)$}& \textbf{5} & \textbf{-1} & \textbf{-1} & \textbf{-5} & \textbf{1} & \textbf{1} & \textbf{5}\tabularnewline
\hline
$M^{P^*,R^P}(t)$ & 5 & 0 & 0 & 0 & 6 & 6 & 10\tabularnewline
\hline
$M^{B,R^{B}}(t)$ & 5 & 5 & 5 & 1 & 1 & 1 & 5\tabularnewline
\hline
\end{tabular}
\end{center}

\end{restatable}
The detailed description of Example \ref{ex:returnsmess} is in shown in Appendix \ref{app:exreturnsmess}.

Thus,  this example shows how the return component can change the relative ordering between the processes of interest. This leads to a more challenging analysis. To address this, our stochastic ordering
analysis is first conducted on the four processes after removing the
return component. Then, a second analysis is performed with the returns
reintroduced.

\subsection{Ordering of Processes without Returns}\label{subsecion:OrdNR}

We first define the following simplified versions
to remove return components for now. The definitions, in order, are: no-returns Infinite Acceptance Process, no-turns Blocking Process,  Skorokhod Approximation using no-returns
Infinite Acceptance Process, and the no-returns Realistic Partial
Fulfillment process.

\begin{align}
M^{\infty}\left(t\right) & =M_{0}+\sum_{i=1}^{N^{d}\left(t\right)}d_{i}-\sum_{j=1}^{N^{b}\left(t\right)}b_{j},\label{eq:M-infty}\\
M^{B}\left(t\right) & =M_{0}+\sum_{i=1}^{N^{d}\left(t\right)}d_{i}-\sum_{j=1}^{N^{b}\left(t\right)}b_{j}\left\{ M^{B}\left(a_{j}-\right)\geq b_{j}\right\},\label{eq:M-B}\\
M^{\infty*}\left(t\right) & =M_{0}+\sum_{i=1}^{N^{d}\left(t\right)}d_{i}-\sum_{j=1}^{N^{b}\left(t\right)}b_{j}\wedge M^{\infty*}\left(a_{j}-\right),\label{eq:M-infty-star}\\
M^{P}\left(t\right) & =M_{0}+\sum_{i=1}^{N^{d}\left(t\right)}d_{i}-\sum_{j=1}^{N^{b}\left(t\right)}b_{j}\wedge M^{P}\left(a_{j}-\right)\label{eq:M-P}.
\end{align}

We will show the order is 
\begin{gather}
M^{B}\left(t\right)\geq M^{\infty*}\left(t\right)=M^{P}(t)\geq M^{\infty}\left(t\right).\label{eq:AllOrderNR}
\end{gather}

\begin{figure}[H]
\centering
\includegraphics[scale=0.6]{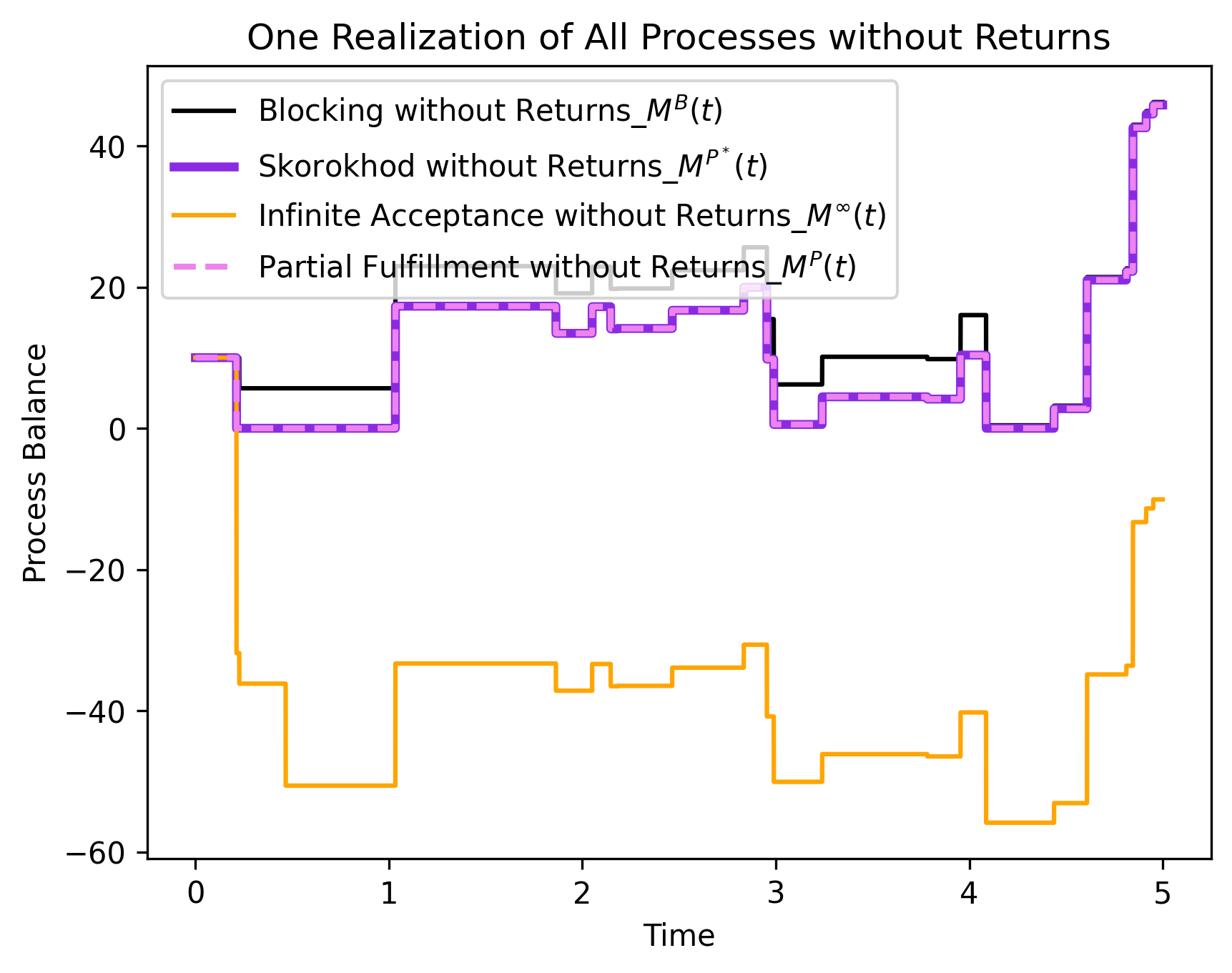}
\caption{Stochastic Ordering of the No-returns Processes}
\label{fig:OrderNR}
\end{figure}

This order is first shown by the simulation results in Figure \ref{fig:OrderNR}.

We now prove the above stochastic order in Theorem \ref{thm:orderNR}.

\begin{restatable}[label={thm:orderNR}]{theorem}{orderNRTheorem}
\label{thm:orderNR}
    Consider the processes in Equations \eqref{eq:M-infty}, \eqref{eq:M-B}, \eqref{eq:M-infty-star}, and \eqref{eq:M-P}. Then, the following is true

    \begin{gather}
M^{B}\left(t\right)\geq M^{\infty*}\left(t\right)=M^{P}(t)\geq M^{\infty}\left(t\right).\label{eq:order}
\end{gather}
\end{restatable}

\begin{proof}
 The proof of this result is contained in Appendix \ref{app:OrdNR}.
\end{proof}

\subsection{Ordering of Processes with Returns}
\label{subsection:OrdR}
We now introduce the following entities that will be the focus for this section:
\begin{align}
M^{P^{*},R^{\infty}} & :=\phi\left(M^{\infty,R^{\infty}}\left(t\right)\right)\nonumber \\
 & =M_{0}+\sum_{i=1}^{N^{d}\left(t\right)}d_{i}-\sum_{j=1}^{N^{b}\left(t\right)}b_{j}\wedge M^{P^{*},R^{\infty}}\left(a_{j}-\right)+\sum_{j=1}^{N^{b}\left(t\right)}p_{j}b_{j}\left\{ a_{j}+s_{j}\leq t\right\} .\label{eq:SkrR}
\end{align}
 Equation \eqref{eq:SkrR} is the Skorokhod Map applied to the Infinite Acceptance Model \eqref{eq:MIR}.

Compared to Equation \ref{eq:Partial}, the model in \eqref{eq:SkrR} will be mathematically simpler to analyze but will deviate from reality when a bail request is only partially fulfilled. In that case, Equation \ref{eq:Partial} will correctly add back the partially fulfilled amount as returns in the third term, but \eqref{eq:SkrR} will add back the entire requested amount even if it is not entirely fulfilled. Thus, \eqref{eq:SkrR}'s utility is that it is a mathematically simpler approximation of (10).
As we observed in
Example \ref{ex:returnsmess}, not all of the orderings in Equation \ref{eq:order}
can be preserved. After introducing returns to the model, we can only prove, for now, the following
ordering:

\[
M^{\infty,R^{\infty}}\left(t\right)\leq M^{P,R^{P}}\left(t\right)\leq M^{P^{*},R^{\infty}}\left(t\right).
\label{eq:AllOrdR}\]

\begin{figure}

\centering
\includegraphics[scale=0.6]{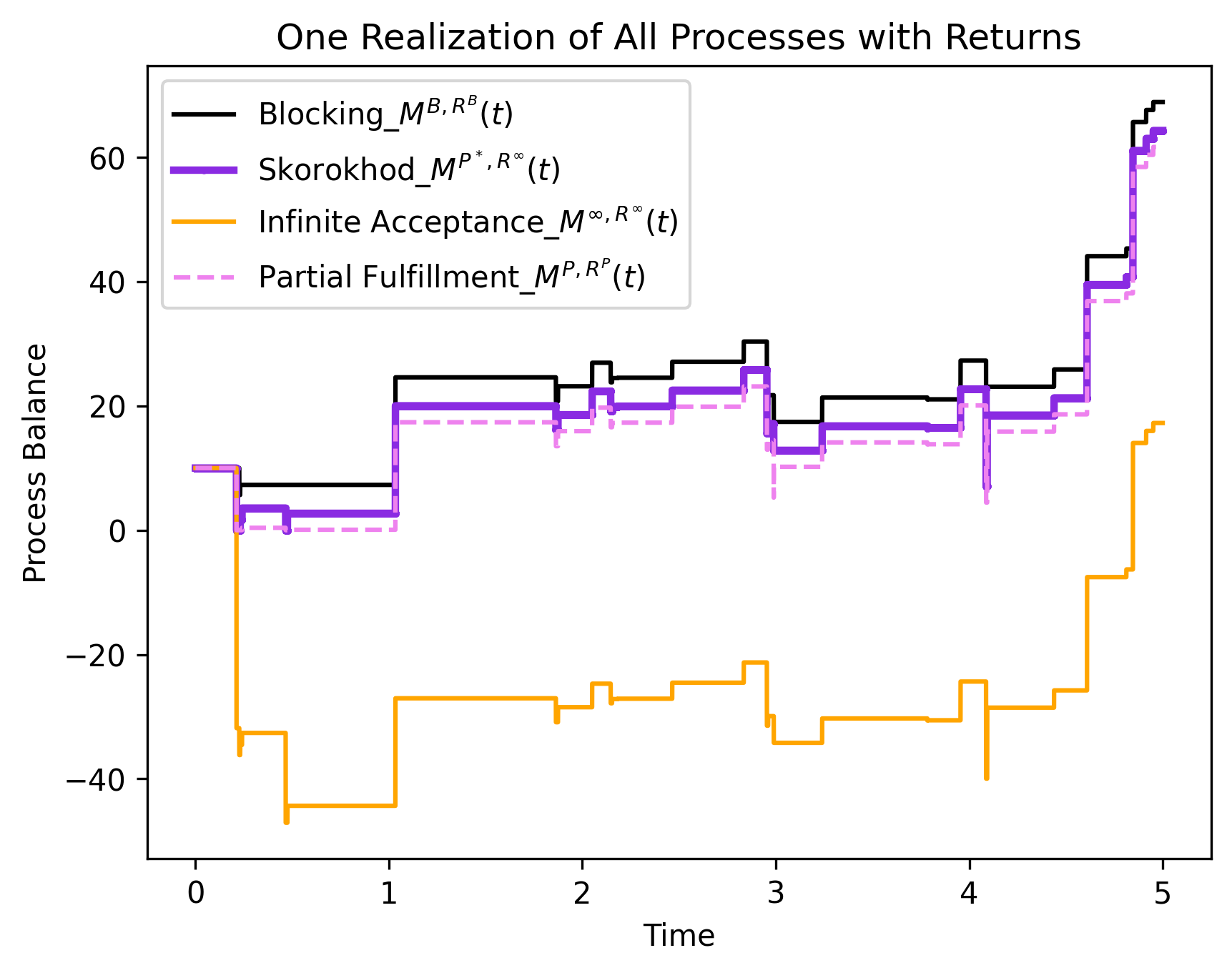}
\caption{Stochastic Ordering of the Standard Processes with Returns}
\label{fig:OrdR}
\end{figure}

Simulations of the processes with returns is shown in Figure \ref{fig:OrdR} before the proofs.

We now prove the above stochastic ordering in Theorem \ref{thm:OrdR}.

\begin{restatable}[label={thm:OrdR}]{theorem}{OrdRTheorem}
\label{thm:OrdR} Consider the processes in Equations \eqref{eq:MIR}, \eqref{eq:SkrP} and \eqref{eq:Partial}. The following holds
\[
M^{\infty,R^{\infty}}\left(t\right)\leq M^{P,R^{P}}\left(t\right)\leq M^{P^{*},R^{\infty}}\left(t\right).
\]
\end{restatable}
\begin{proof}
 The proof of this result is contained in Appendix \ref{app:OrdR}.   
\end{proof}

\section{Conclusion}

In this paper, we introduced and analyzed stochastic models to approximate the complex dynamics of community bail funds (CBFs). We developed three distinct models: an infinite acceptance model, a partial fulfillment model, and a blocking model. For these models, we derived fluid limits that describe their large-scale dynamics. Notably, to aid our analysis, we prove how the partial fulfillment model, which corresponds to an interpretable policy, coincides exactly to applying the Skorokhod mapping \citep{skorokhod1961stochastic} to the sample path of the infinite acceptance process. To support our theoretical results, we further present numerical simulations that support the accuracy of our approximations and limits. Finally, we conduct a stochastic ordering analysis of our proposed models, revealing the relationships between the balances of these CBF approximations. Insights from our models and analysis could lead to deeper understandings of CBF operations and provide a quantitative foundation for future research into optimizing CBF policies, enhancing their impact on the justice system.

\section*{Acknowledgements}

Jamol Pender would like to acknowledge the support of Tara Watford and Katy Moose of The Bail Project for their time and willingness to explain how a bail fund actually works to us. 

The authors used Overleaf's  Writefull and ChatGPT as an editor for non-mathematical content.

\bibliographystyle{plainnat}
\bibliography{references}

\appendix
\section{Appendix}

\subsection{Proof of the Results in Section \ref{subsubsec:CPConverge}}

\subsubsection{\label{app:MgMG} Proof of the Lemma \ref{lem:MgMG}.}

\MgMGLemma*

\begin{proof}
It suffices to check the Martingale property for any $s<t$: in other
words, $\E\left[\M_{\bowtie}^{\eta}(t)\Big|\F_{s}\right]=\M_{\bowtie}^{\eta}(s)$.
\begin{align*}
\E\left[\M_{\bowtie}^{\eta}(t)\Big|\F_{s}\right] & =\E\left[\M_{\bowtie}^{\eta}(t)-\M_{\bowtie}^{\eta}(s)+\M_{\bowtie}^{\eta}(s)\Big|\F_{s}\right]\\
 & =\E\left[\M_{\bowtie}^{\eta}(t)-\M_{\bowtie}^{\eta}(s)\Big|\F_{s}\right]+\E\left[\M_{\bowtie}^{\eta}(s)\Big|\F_{s}\right]\\
 & =\E\left[\M_{\bowtie}^{\eta}(t)-\M_{\bowtie}^{\eta}(s)\Big|\F_{s}\right]+\M_{\bowtie}^{\eta}(s).
\end{align*}
So, it suffices to prove $\E\left[\M_{\bowtie}^{\eta}(t)-\M_{\bowtie}^{\eta}(s)\Big|\F_{s}\right]=0$.
\begin{align*}
\E\left[\M_{\bowtie}^{\eta}(t)-\M_{\bowtie}^{\eta}(s)\Big|\F_{s}\right] & =\E\left[\frac{1}{\eta}\sum_{i=1}^{N_{\eta}\left(t\right)}v_{i,\eta}-\lambda v_{\eta}^{*}t-\frac{1}{\eta}\sum_{i=1}^{N_{\eta}\left(s\right)}v_{i,\eta}+\lambda v_{\eta}^{*}s\Big|\F_{s}\right]\\
 & =\E\left[\frac{1}{\eta}\sum_{i=N_{\eta}\left(s\right)+1}^{N_{\eta}\left(t\right)}v_{i,\eta}+\lambda v_{\eta}^{*}\left(s-t\right)\Big|\F_{s}\right]\\
 & =\E\left[\frac{1}{\eta}\sum_{i=N_{\eta}\left(s\right)+1}^{N_{\eta}\left(t\right)}v_{i,\eta}\Big|\F_{s}\right]+\lambda v_{\eta}^{*}\left(s-t\right)\\
 & =\E\left[\E\left[\frac{1}{\eta}\sum_{i=N_{\eta}\left(s\right)+1}^{N_{\eta}\left(t\right)}v_{i,\eta}\Bigg|N_{\eta}\left(t\right),\F_{s}\right]\Big|\F_{s}\right]+\lambda v_{\eta}^{*}\left(s-t\right)\\
 & =\E\left[\left(N_{\eta}\left(t\right)-N_{\eta}\left(s\right)\right)\frac{1}{\eta}v_{\eta}^{*}\Big|\F_{s}\right]+\lambda v_{\eta}^{*}\left(s-t\right),
\end{align*}

where the last line follows from the assumption that $v_{i}'s$ are
i.i.d conditional on $N_{\eta}\left(t\right)$. Since the interarrival
times are i.i.d, then $N_{\eta}\left(t\right)-N_{\eta}\left(s\right)$
conditional on $\F_{s}$ has the same distribution as $N_{\eta}\left(t-s\right)$.
So, 
\begin{align*}
\E\left[\M_{\bowtie}^{\eta}(t)-\M_{\bowtie}^{\eta}(s)\Big|\F_{s}\right] & =\E\left[\left(N_{\eta}\left(t\right)-N_{\eta}\left(s\right)\right)\frac{1}{\eta}v_{\eta}^{*}\Big|\F_{s}\right]+\lambda v_{\eta}^{*}\left(s-t\right)\\
 & =\frac{1}{\eta}v_{\eta}^{*}\E\left[N_{\eta}\left(t-s\right)\right]+\lambda v_{\eta}^{*}\left(s-t\right)\\
 & =\frac{1}{\eta}v_{\eta}^{*}\lambda\eta\left(t-s\right)+\lambda v_{\eta}^{*}\left(s-t\right)\\
 & =0.
\end{align*}
This completes the proof. 
\end{proof}

\subsubsection{\label{app:Gcpconverge} Proof of the Lemma \ref{lem:Gcpconverge}.}

\GcpconvergeLemma*

\begin{proof}
 By Lemma \ref{lem:MgMG}
\[
\M_{\bowtie}^{\eta}(t):=\frac{1}{\eta}\sum_{i=1}^{N_{\eta}\left(t\right)}v_{i,\eta}-\lambda v_{\eta}^{*}t
\]
is a martingale.

We will use the Borel-Cantelli lemma to show almost sure convergence.
We start by defining the following events for some $\epsilon>0$ and
integers $\eta=k,k+1,k+2,..$. and some starting finite index $k$:
\[
\A_{\eta}:=\left\{ \sup_{0\leq t\leq T}\left|\M_{\bowtie}^{\eta}(t)\right|>\epsilon\right\} .
\]
We start indexing $\eta$ at $k$ because it will help make moment
generating functions valid later as well as simplify later inequalities.
It does not affect using Borel-Cantelli because it only excludes a
finite number of beginning indices. 
\begin{align*}
\P\left\{ \A_{\eta}\right\}  & =\P\left\{ \sup_{0\leq t\leq T}\left|\M_{\bowtie}^{\eta}(t)\right|>\epsilon\right\} \\
 & \leq\P\left\{ \sup_{0\leq t\leq T}\M_{\bowtie}^{\eta}(t)\geq\epsilon\right\} +\P\left\{ \inf_{0\leq t\leq T}\M_{\bowtie}^{\eta}(t)\leq-\epsilon\right\} \\
 & =\P\left\{ \sup_{0\leq t\leq T}\M_{\bowtie}^{\eta}(t)\geq\epsilon\right\} +\P\left\{ \sup_{0\leq u\leq T}-\M_{\bowtie}^{\eta}(u)\leq\epsilon\right\} \\
 & =\P\left\{ \sup_{0\leq t\leq T}e^{\kappa\M_{\bowtie}^{\eta}(t)}\geq e^{\kappa\epsilon}\right\} +\P\left\{ \sup_{0\leq u\leq T}e^{-\kappa\M_{\bowtie}^{\eta}(t)}\leq e^{\kappa\epsilon}\right\} ,
\end{align*}
where we choose $\kappa=\frac{\sqrt{\eta}}{\epsilon}$ because it
will lead us to convergence later. Since $\M_{\bowtie}^{\eta}(t)$
is a martingale, $-\M_{\bowtie}^{\eta}(t)$ is also a martingale.
Also, $e^{\kappa x}$is a convex function. Then both $e^{\kappa\M_{\bowtie}^{\eta}(t)}$
and $e^{-\kappa\M_{\bowtie}^{\eta}(t)}$are submartingales. Then,
by the Doob's maximal inequality for continuous-time submartingales:
\begin{align*}
\P\left\{ \A_{\eta}\right\}  & \leq\frac{\E\left[e^{\kappa\M_{\bowtie}^{\eta}(T)}\right]}{e^{\kappa\epsilon}}+\frac{\E\left[e^{-\kappa\M_{\bowtie}^{\eta}(T)}\right]}{e^{\kappa\epsilon}}\\
 & =\frac{1}{e^{\kappa\epsilon}}\left(\E\left[e^{\kappa\frac{1}{\eta}\sum_{i=1}^{N_{\eta}\left(T\right)}v_{i,\eta}-\lambda v_{\eta}^{*}T\kappa}\right]+\E\left[e^{-\kappa\frac{1}{\eta}\sum_{i=1}^{N_{\eta}\left(T\right)}v_{i,\eta}+\lambda v_{\eta}^{*}T\kappa}\right]\right)
\end{align*}
To calculate these two terms, observe that 
\[
\E\left[e^{\kappa\frac{1}{\eta}\sum_{i=1}^{N_{\eta}\left(T\right)}v_{i,\eta}}\right]=\E\left[\E\left[e^{\kappa\frac{1}{\eta}\sum_{i=1}^{N_{\eta}\left(T\right)}v_{i,\eta}}\Big|N_{\eta}\left(T\right)\right]\right].
\]
For the inner expectation, we have 
\begin{align*}
\E\left[e^{\kappa\frac{1}{\eta}\sum_{i=1}^{N_{\eta}\left(T\right)}v_{i,\eta}}\Big|N_{\eta}\left(T\right)=n\right] & =\E\left[\prod_{i=1}^{n}e^{\kappa\frac{1}{\eta}v_{i,\eta}}\right]\\
 & =\prod_{i=1}^{n}\E\left[e^{\frac{\kappa}{\eta}v_{i,\eta}}\right]\\
 & =\prod_{i=1}^{n}\zeta\left(\frac{\kappa}{\eta}\right)\\
 & =\left(\zeta\left(\frac{\kappa}{\eta}\right)\right)^{n},
\end{align*}
where the second line follows from the assumption that $v_{i,\eta}$
are i.i.d conditioned on $N_{\eta}^{b}\left(t\right)$, and in the
second-to-last equality, and we define the moment generating
function (MGF) of $v_{i}$ by $\zeta\left(\frac{\kappa}{\eta}\right):=\E\left[e^{\frac{\kappa}{\eta}v_{i,\eta}}\right]$,
assuming it exists.

Then, 
\begin{align*}
\E\left[\E\left[e^{\kappa\frac{1}{\eta}\sum_{i=1}^{N_{\eta}\left(T\right)}v_{i}}\Big|N_{\eta}\left(T\right)\right]\right] & =\E\left[\left(\zeta\left(\frac{\kappa}{\eta}\right)\right)^{N_{\eta}\left(T\right)}\right]\\
 & =\text{\ensuremath{\sum_{n=0}^{\infty}}}\left(\zeta\left(\frac{\kappa}{\eta}\right)\right)^{n}\P\left(N_{\eta}\left(T\right)=n\right)\\
 & =\text{\ensuremath{\sum_{n=0}^{\infty}}}\left(\zeta\left(\frac{\kappa}{\eta}\right)\right)^{n}\frac{\left(\lambda\eta T\right)^{n}e^{-\lambda\eta T}}{n!}\\
 & =e^{-\lambda\eta T}\text{\ensuremath{\sum_{n=0}^{\infty}}}\frac{\left(\zeta\left(\frac{\kappa}{\eta}\right)\lambda\eta T\right)^{n}}{n!}\\
 & =e^{-\lambda\eta T}e^{\zeta\left(\frac{\kappa}{\eta}\right)\lambda\eta T}\\
 & \leq e^{-\lambda\eta T}e^{\left(\zeta\left(0\right)+\zeta^{'}\left(0\right)\frac{\kappa}{\eta}+\frac{Q}{2}\left(\frac{\kappa}{\eta}\right)^{2}\right)\lambda\eta T}\\
 & =e^{-\lambda\eta T}e^{\left(1+v_{\eta}^{*}\frac{\kappa}{\eta}+\frac{Q}{2}\left(\frac{\kappa}{\eta}\right)^{2}\right)\lambda\eta T}\\
 & =e^{\lambda v_{\eta}^{*}\kappa T+\frac{Q\kappa^{2}\lambda T}{2\eta}}
\end{align*}
where the second step follows from $N_{\eta}\left(T\right)\sim Poiss(\lambda\eta T)$,
and the the third-to-last line follows from the Taylor's theorem and
$c$ is the constant in the upper error bound term of Taylor's Theorem.
Simplifying, we have 
\begin{align*}
\E\left[e^{\kappa\frac{1}{\eta}\sum_{i=1}^{N_{\eta}\left(T\right)}v_{i,\eta}}\right] & \leq e^{\lambda v_{\eta}^{*}\kappa T+\frac{Q\kappa^{2}\lambda T}{2\eta}}.
\end{align*}
Similarly, by an analogous derivation 
\begin{align*}
\E\left[e^{-\kappa\frac{1}{\eta}\sum_{i=1}^{N_{\eta}\left(T\right)}v_{i,\eta}}\right] & \leq e^{-\lambda v_{\eta}^{*}\kappa T-\frac{Q\kappa^{2}\lambda T}{2\eta}}.
\end{align*}
So, as $\eta\rightarrow\infty,$ 
\begin{align*}
\P\left\{ \A_{\eta}\right\}  & \leq\frac{1}{e^{\kappa\epsilon}}\left(\E\left[e^{\kappa\frac{1}{\eta}\sum_{i=1}^{N_{\eta}\left(T\right)}v_{i,\eta}-\lambda v_{\eta}^{*}T\kappa}\right]+\E\left[e^{-\kappa\frac{1}{\eta}\sum_{i=1}^{N_{\eta}\left(T\right)}v_{i,\eta}+\lambda v_{\eta}^{*}T\kappa}\right]\right)\\
 & =\frac{1}{e^{\kappa\epsilon}}\left(e^{\lambda v_{\eta}^{*}\kappa T+\frac{Q\kappa^{2}\lambda T}{2\eta}-\lambda v_{\eta}^{*}T\kappa}+e^{-\lambda v^{*}\kappa T-\frac{Q\kappa^{2}\lambda T}{2\eta}+\lambda v_{\eta}^{*}T\kappa}\right)\\
 & =\frac{1}{e^{\kappa\epsilon}}\left(e^{0}+e^{0}\right)\\
 & =\frac{2}{e^{\kappa\epsilon}}
\end{align*}
We plug in $\kappa=\frac{\sqrt{\eta}}{\epsilon},$ 
\begin{align*}
\P\left\{ \A_{\eta}\right\}  & \leq\frac{2}{e^{\sqrt{\eta}}}
\end{align*}
To use Borel-Cantelli, we take the sum 
\begin{align*}
\sum_{\eta=1}^{\infty}\P\left\{ \A_{\eta}\right\}  & \leq2\sum_{\eta=1}^{\infty}\frac{1}{e^{\sqrt{\eta}}}\\
 & =2\sum_{\eta=1}^{\infty}e^{-\sqrt{\eta}}
\end{align*}
 Therefore, 
\[
\sum_{\eta=k}^{\infty}\P\left\{ \A_{\eta}\right\} <\infty.
\]
By the Borel-Cantelli theorem, 
\begin{gather*}
\P\left\{ \sup_{0\leq t\leq T}\left|\frac{1}{\eta}\sum_{i=1}^{N_{\eta}\left(T\right)}v_{i,\eta}-\lambda v_{\eta}^{*}t\right|\geq\epsilon\ \text{i.o.}\right\} =\P\left\{ \sup_{0\leq t\leq T}\left|\M_{\bowtie}^{\eta}\left(t\right)\right|\geq\epsilon\ \text{i.o.}\right\} =0.
\end{gather*}
Since this holds for any $\epsilon>0$, this implies 
\[
\sup_{0\leq t\leq T}\left|\frac{1}{\eta}\sum_{i=1}^{N_{\eta}\left(T\right)}v_{i,\eta}-\lambda v_{\eta}^{*}t\right|\underset{a.s.}{\to}0.
\]
\end{proof}
\subsubsection{\label{app:ReturnAS} Proof of the Lemma \ref{lem:ReturnAS}}

\ReturnASLemma*

\begin{proof}
This result follows from applying the Borel-Cantelli lemma to desired
processes, which was already done in Lemma \ref{lem:Gcpconverge}. We
need only show that the summation term $\frac{1}{\eta}\sum_{j=1}^{N_{\eta}^{b}(t)}\left(1-p_{j}\right)b_{j}\left\{ u>a_{j,\eta}+s_{j}\right\} ,$
in particular, fits the assumptions of Lemma \ref{lem:Gcpconverge} for
the Borel-Cantelli lemma to be applied. To this end, observe when
conditioned on $N_{\eta}^{b}\left(t\right)=n$, the $\left\{ a_{j,\eta}\right\} ^{n}$
are the order statistics of $n$ i.i.d uniform random variables on
$\left[0,t\right]$. Let $\left\{ u_{j}\right\} _{j=1}^{n}$ be $n$
independent $\left[0,1\right]$ uniform random variables, and $u_{\left(j\right)}$
is the $j$-th order statistics of $\left\{ u_{j}\right\} _{j=1}^{n}$
. Then, $a_{j,\eta}$ has the same distribution as $tu_{\text{\ensuremath{\left(j\right)}}}$.
We define $v_{j,\eta}$ as the random variable $\left(1-p_{j}\right)b_{j}\left\{ t>a_{j,\eta}+s_{j}\right\} $.
So, the summation conditional on $N_{\eta}^{b}\left(t\right)=n$ can
be rearranged as follows:

\begin{align*}
\frac{1}{\eta}\sum_{j=1}^{N_{\eta}^{b}(t)=n}v_{j,\eta} & =\frac{1}{\eta}\sum_{j=1}^{N_{\eta}^{b}(t)=n}\left(1-p_{j}\right)b_{j}\left\{ t>a_{j,\eta}+s_{j}\right\} \\
 & =\frac{1}{\eta}\sum_{j=1}^{n}\left(1-p_{j}\right)b_{j}\left\{ t>tu_{\left(j\right)}+s_{j}\right\} \\
 & =\frac{1}{\eta}\sum_{\sigma\left(j\right)=1}^{n}\left(1-p_{\sigma\left(j\right)}\right)b_{\sigma\left(j\right)}\left\{ t>tu_{j}+s_{\sigma\left(j\right)}\right\} 
\end{align*}
where $\sigma\left(j\right)$ in the second-to-last line is the permutation
that reverses the ordering on $u_{\left(j\right)}$ into $u_{j}$
.

Since $p_{\sigma\left(j\right)},b_{\sigma\left(j\right),}u_{j},s_{\sigma\left(j\right)}$
are all i.i.d random variables, it can be shown that the summands
$v_{j,\eta}$ are i.i.d conditional on $N_{\eta}^{b}\left(t\right)$.

We then compute the conditional expectation of the $v_{j,\eta}$ and
denote it by $v_{j,\eta}^{*}$:

\begin{align*}
v_{j,\eta}^{*} & =\E\left[v_{j,\eta}\Bigg|N_{\eta}^{b}(t)=n\right]\\
 & =\E\left[\left(1-p_{\sigma\left(j\right)}\right)b_{\sigma\left(j\right)}\left\{ t>tu_{j}+s_{\sigma\left(j\right)}\right\} \right]\\
 & =\left(1-p^{*}\right)b^{*}\E\left[\E\left[\left\{ t-tu_{j}>s_{\sigma\left(j\right)}\right\} \Bigg|u_{j}\right]\right]\\
 & =\left(1-p^{*}\right)b^{*}\E\left[F_{s}\left(t-tu_{j}\right)\right]\\
 & =\left(1-p^{*}\right)b^{*}\int_{0}^{1}F_{s}\left(t-tx\right)dx.
\end{align*}
Let$w=tx$ and $dw=tdx$, and then, 
\begin{align*}
v_{j,\eta}^{*} & =\left(1-p^{*}\right)b^{*}\int_{0}^{t}F_{s}\left(t-w\right)\frac{1}{t}dw.
\end{align*}
Change the variable again by substituting $t-w$ with $v$, so $dv=-dw$
and we have
\begin{align*}
v_{j,\eta}^{*} & =\left(1-p^{*}\right)b^{*}\int_{t}^{0}F_{s}\left(v\right)\frac{1}{t}\left(-1\right)dv\\
 & =\frac{1}{t}\left(1-p^{*}\right)b^{*}\int_{0}^{t}F_{s}\left(v\right)dv\text{,}
\end{align*}
which shows the value of $v_{j,\eta}^{*}$ only depends on $\eta.$ 

Thus, by Lemma \ref{lem:Gcpconverge} , for any fixed time T, as $\eta\rightarrow\infty$,
\[
\sup_{0\leq t\leq T}\left|\frac{1}{\eta}\sum_{j=1}^{N_{\eta}^{b}\left(t\right)}v_{j,\eta}-\text{\ensuremath{\lambda_{b}v_{\eta}^{*}t}}\right|\underset{a.s.}{\to}0.
\]
More precisely,

\[
\sup_{0\leq t\leq T}\left|\frac{1}{\eta}\sum_{j=1}^{N_{\eta}^{b}\left(t\right)}\left(1-p_{j}\right)b_{j}\left\{ t>a_{j,\eta}+s_{j}\right\} -\lambda_{b}\left(1-p^{*}\right)b^{*}\int_{0}^{t}F_{s}\left(v\right)dv\right|\underset{a.s.}{\to}0.
\]
This completes the proof. 
\end{proof}

\subsection{\label{app:partial}Proof of the results in Section \ref{sec:SkrkandParitial}.}
\subsubsection{Proof of the Corollary \ref{col:partial}}

\partialCorollary*

\begin{proof}
This result follows directly from Theorem \ref{thm:partial}. Consider
our Infinite Acceptance Model (Equation \eqref{eq:MIR}) as an instance
within the class of simple bail processes (Definition \ref{def:simple}).
Using Theorem \ref{thm:partial}, we explicitly decompose each component
of the Skorokhod output process for Infinite Acceptance Process below.

\begin{align*}
M^{P^{*},R^{\infty}}\left(t\right) & :=\phi\left[M^{\infty,R^{\infty}}\right]\left(t\right)\\
 & =G^{\infty,R^{\infty}}\left(t\right)-\sum_{i=1}^{N^{b}\left(t\right)}b_{i}\wedge M^{P^{*},R^{\infty}}\left(a_{i}-\right)\\
 & =\underbrace{M_{0}+\sum_{i=1}^{N^{d}\left(t\right)}d_{i}+\sum_{i=1}^{N^{b}\left(t\right)}p_{i}b_{i}\left\{ a_{i}+s_{i}\leq t\right\} }_{G^{\infty,R^{\infty}}\left(t\right)}-\sum_{i=1}^{N^{b}\left(t\right)}b_{i}\wedge M^{P^{*},R^{\infty}}\left(a_{i}-\right)\\
 & =M_{0}+\sum_{i=1}^{N^{d}\left(t\right)}d_{i}-\underbrace{\sum_{i=1}^{N^{b}\left(t\right)}b_{i}\wedge M^{P^{*},R^{\infty}}\left(a_{i}-\right)}_{\text{Bail Returns Partially Fulfilled}}+\underbrace{\sum_{i=1}^{N^{b}\left(t\right)}p_{i}b_{i}\left\{ a_{i}+s_{i}\leq t\right\} }_{\text{Bail Returns from Fully Requested Amount}}
\end{align*}
This completes the proof.
\end{proof}

\subsection{Proof of the Results in Section \ref{subsubsec:BRConverge}.}

\subsubsection{Supporting Lemmas for the Lemma \ref{lem:BCP}.}

\begin{restatable}[label={lem:BMG1}]{lemma}{BMG1Lemma}
\label{lem:BMG1}The process 
\[
\M_{\eta,1}\left(t\right):=\frac{1}{\eta}\sum_{j=1}^{N_{\eta}^{b}(t)}b_{j}\left\{ M_{\eta}^{B,R}\left(a_{j}-\right)\geq b_{j}\right\} -\text{\ensuremath{\frac{1}{\eta}}\ensuremath{\sum_{j=1}^{N_{\eta}^{b}( t)}H\left(M_{\eta}^{B,R}\left(a_{j}-\right)\right)}}
\]
is a martingale with respect to the filtration $\left\{ \F_{s}\right\} $
which contains all bail arrivals and bail request sizes until time
$s$:
\[
\F_{s}:=\sigma\left(N_{\eta}^{b}\left( s\right),\left\{ b_{j}\right\} _{j=1}^{N_{\eta}^{b}\left(s\right)},\left\{ a_{j}\right\} _{j=1}^{N_{\eta}^{b}\left( s\right)}\right).
\]
\end{restatable}
\begin{proof}
It suffices to check the martingale property for any $s<t$: in other
words, $\E\left[\M_{\eta,1}\left(t\right)-\M_{\eta,{1}}\left(s\right)\Big|\F_{s}\right]=0$.

{\footnotesize
\begin{align*}
\E\left[\M_{\eta,b_{1}}\left(t\right)-\M_{\eta,b_{1}}\left(s\right)\Big|\F_{s}\right] & =\E\left[\frac{1}{\eta}\sum_{j=1}^{N_{\eta}^{b}(t)}b_{j}\left\{ M_{\eta}^{B,R}\left(a_{j}-\right)\geq b_{j}\right\} -\text{\ensuremath{\frac{1}{\eta}\sum_{j=1}^{N_{\eta}^{b}(t)}H\left(M_{\eta}^{B,R}\left(a_{j}-\right)\right)}}\Bigg|\F_{s}\right]\\
 & \quad-\E\left[\left(\frac{1}{\eta}\sum_{j=1}^{N_{\eta}^{b}(s)}b_{j}\left\{ M_{\eta}^{B,R}\left(a_{j}-\right)\geq b_{j}\right\} -\text{\ensuremath{\frac{1}{\eta}\sum_{j=1}^{N_{\eta}^{b}(s)}H\left(M_{\eta}^{B,R}\left(a_{j}-\right)\right)}}\right)\Bigg|\F_{s}\right]\\
 & =\E\left[\frac{1}{\eta}\sum_{j=N_{\eta}^{b}(s)+1}^{N_{\eta}^{b}(t)}b_{j}\left\{ M_{\eta}^{B,R}\left(a_{j}-\right)\geq b_{j}\right\} -\text{\ensuremath{\frac{1}{\eta}\sum_{j=N_{\eta}^{b}(s)+1}^{N_{\eta}^{b}(t)}H\left(M_{\eta}^{B,R}\left(a_{j}-\right)\right)}}\Bigg|\F_{s}\right]\\
 & =\E\left[\frac{1}{\eta}\sum_{j=N_{\eta}^{b}(s)+1}^{N_{\eta}^{b}(t)}b_{j}\left\{ M_{\eta}^{B,R}\left(a_{j}-\right)\geq b_{j}\right\} -\text{\ensuremath{H\left(M_{\eta}^{B,R}\left(a_{j}-\right)\right)}}\Bigg|\F_{s}\right]\\
 & =\E\left[\frac{1}{\eta}\E\left[\sum_{j=N_{\eta}^{b}(s)+1}^{N_{\eta}^{b}(t)}b_{j}\left\{ M_{\eta}^{B,R}\left(a_{j}-\right)\geq b_{j}\right\} -\text{\ensuremath{H\left(M_{\eta}^{B,R}\left(a_{j}-\right)\right)}}\Bigg|N_{\eta}^{b}(t),\F_{s}\right]\Bigg|\F_{s}\right]\\
 & =\E\left[\frac{1}{\eta}\sum_{j=N_{\eta}^{b}(s)+1}^{N_{\eta}^{b}(t)}\E\left[b_{j}\left\{ M_{\eta}^{B,R}\left(a_{j}-\right)\geq b_{j}\right\} -\text{\ensuremath{H\left(M^{B,R,\eta}\left(a_{j}-\right)\right)}}\Bigg|N^{b}(\eta t),\F_{s}\right]\Bigg|\F_{s}\right]\\
 & =\scalebox{0.8}{$\E\left[\frac{1}{\eta}\sum_{j=N_{\eta}^{b}(s)+1}^{N_{\eta}^{b}(t)}\E\left[\E\left[b_{j}\left\{ M_{\eta}^{B,R}\left(a_{j}-\right)\geq b_{j}\right\} -\text{\ensuremath{H\left(M_{\eta}^{B,R}\left(a_{j}-\right)\right)\Bigg|M_{\eta}^{B,R}\left(a_{j}-\right),N_{\eta}^{b}(t),\F_{s}}}\right]\Bigg|N_{\eta}^{b}(t),\F_{s}\right]\Bigg|\F_{s}\right]$}\\
 & =\E\left[\frac{1}{\eta}\sum_{j=N_{\eta}^{b}(s)+1}^{N_{\eta}^{b}(t)}\E\left[H\left(M_{\eta}^{B,R}\left(a_{j}-\right)\right)-H\left(M_{\eta}^{B,R}\left(a_{j}-\right)\right)\Bigg|N_{\eta}^{b}(t),\F_{s}\right]\Bigg|\F_{s}\right]\\
 & =0.
\end{align*}
}This completes the proof.
\end{proof}

\begin{restatable}[label={lem:BMG3}]{lemma}{BMG3Lemma}
\label{lem:BMG3}The process 
\[
\M_{\eta,2}\left(t\right):=\text{\ensuremath{\frac{1}{\eta}\int_{0}^{t}H\left(M_{\eta}^{B,R}\left(u\right)\right)}}dN_{\eta}\left(u\right)-\lambda_{b}\int_{0}^{t}H\left(M_{\eta}^{B,R}\left(u\right)\right)du
\]
is a martingale with respect to the filtration $\left\{ \F_{s}\right\} $
which contains all bail arrivals and bail request sizes until time
$s$:
\[
\F_{s}:=\sigma\left(N_{\eta}^{b}\left(s\right),\left\{ b_{j}\right\} _{j=1}^{N_{\eta}^{b}\left(s\right)},\left\{ a_{j}\right\} _{j=1}^{N_{\eta}^{b}\left(s\right)}\right).
\]
\end{restatable}
\begin{proof}

It suffices to check the Martingale property for any $s<t$: in other
words, $\E\left[\M_{\eta,2}(t)-\M_{\eta,2}(s)\Big|\F_{s}\right]=0$.

\begin{gather*}
\E\left[\M_{\eta,2}\left(t\right)-\M_{\eta,2}\left(s\right)\Big|\F_{s}\right]\\ \displaybreak[1]
=\E\left[\ensuremath{\frac{1}{\eta}\int_{s}^{t}H\left(M_{\eta}^{B,R}\left(u\right)\right)}dN_{\eta}\left(u\right)-\lambda_{b}\int_{s}^{t}H\left(M_{\eta}^{B,R}\left(u\right)\right)du\Big|\F_{s}\right]\\
=\E\left[\underbrace{\ensuremath{\frac{1}{\eta}\int_{s}^{t}H\left(M_{\eta}^{B,R}\left(u\right)\right)}dN_{\eta}\left(u\right)}_{\text{A}}\right]-\E\left[\lambda_{b}\int_{s}^{t}H\left(M_{\eta}^{B,R}\left(u\right)\right)du\right]
\end{gather*}
 Denote change $\Delta_{j}N_{\eta}=N_{\eta}\left(j\delta\right)-N_{\eta}\left(\left(j-1\right)\delta\right)$,
and $dt$ are approximated with the length $\delta$.
By definition of the stochastic integral , we have:
\begin{align*}
A & =\ensuremath{\frac{1}{\eta}\int_{s}^{t}H\left(M_{\eta}^{B,R}\left(u\right)\right)}dN_{\eta}\left(u\right)\\
 & =\frac{1}{\eta}\lim_{\delta\rightarrow0}\sum_{j=\left\lfloor \frac{s}{\delta}\right\rfloor +1}^{\left\lfloor \frac{t}{\delta}\right\rfloor }H\left(M_{\eta}^{B,R}\left(\left(j-1\right)\delta\right)\right)\Delta_{j}N_{\eta}
\end{align*}
 So that,
\[
\E\left[A\right]=\E\left[\frac{1}{\eta}\lim_{\delta\rightarrow0}\sum_{j=\frac{s}{\delta}+1}^{\left\lfloor \frac{t}{\delta}\right\rfloor }H\left(M_{\eta}^{B,R}\left(\left(j-1\right)\delta\right)\right)\Delta_{j}N_{\eta}\right],
\]
 Observe that for every $\delta>0,$ 
\[
\sum_{j=\left\lfloor \frac{s}{\delta}\right\rfloor +1}^{\left\lfloor \frac{t}{\delta}\right\rfloor }H\left(M_{\eta}^{B,R}\left(\left(j-1\right)\delta\right)\right)\Delta_{j}N_{\eta}\leq\sum_{j=1}^{N_{\eta}\left(t\right)}b^{*}
\]
because H$\left(\cdot\right)\leq b^{*}$ by definition. In addition,
$\E\left[\sum_{j=1}^{N_{\eta}\left(t\right)}b^{*}\right]=\lambda\eta tb^{*}<\infty$.
By the dominated convergence theorem, we interchange
the expectation and limit and proceed with the calculation as follows:
\begin{align*}
\E\left[A\right] & =\E\left[\frac{1}{\eta}\lim_{\delta\rightarrow0}\sum_{j=\frac{s}{\delta}+1}^{\left\lfloor \frac{t}{\delta}\right\rfloor }H\left(M_{\eta}^{B,R}\left(\left(j-1\right)\delta\right)\right)\Delta_{j}N_{\eta}\right]\text{\text{(Dominated Convergence Theorem)}}\\
 & =\frac{1}{\eta}\lim_{\delta\rightarrow0}\E\left[\sum_{j=\frac{s}{\delta}+1}^{\left\lfloor \frac{t}{\delta}\right\rfloor }H\left(M_{\eta}^{B,R}\left(\left(j-1\right)\delta\right)\right)\Delta_{j}N_{\eta}\right]\\  \displaybreak[1]
 & =\frac{1}{\eta}\lim_{\delta\rightarrow0}\sum_{j=\frac{s}{\delta}+1}^{\left\lfloor \frac{t}{\delta}\right\rfloor }\E\left[H\left(M_{\eta}^{B,R}\left(\left(j-1\right)\delta\right)\right)\Delta_{j}N_{\eta}\right]\\ 
 & =\frac{1}{\eta}\lim_{\delta\rightarrow0}\sum_{j=\frac{s}{\delta}+1}^{\left\lfloor \frac{t}{\delta}\right\rfloor }\E\left[\E\left[H\left(M_{\eta}^{B,R}\left(\left(j-1\right)\delta\right)\right)\Delta_{j}N_{\eta}\Bigg|\F_{\left(j-1\right)\delta}\right]\right]\\
 & =\frac{1}{\eta}\lim_{\delta\rightarrow0}\sum_{j=\frac{s}{\delta}+1}^{\left\lfloor \frac{t}{\delta}\right\rfloor }\E\left[H\left(M_{\eta}^{B,R}\left(\left(j-1\right)\delta\right)\right)\E\left[\Delta_{j}N_{\eta}\Bigg|\F_{\left(j-1\right)\delta}\right]\right]\\
 & =\frac{1}{\eta}\lim_{\delta\rightarrow0}\sum_{j=\frac{s}{\delta}+1}^{\left\lfloor \frac{t}{\delta}\right\rfloor }\E\left[H\left(M_{\eta}^{B,R}\left(\left(j-1\right)\delta\right)\right)\lambda_{b}\eta\delta\right]\\
 & =\E\left[\lim_{\delta\rightarrow0}\sum_{j=\left\lfloor \frac{s}{\delta}\right\rfloor +1}^{\left\lfloor \frac{t}{\delta}\right\rfloor }H\left(M_{\eta}^{B,R}\left(\left(j-1\right)\delta\right)\right)\lambda_{b}\delta\right]\\
 & =\E\left[\int_{s}^{t}H\left(M_{\eta}^{B,R}\left(u\right)\right)\lambda_{b}du\right],
\end{align*}
where the third-to-last line used the Poisson property, so $\Delta_{j}N_{\eta}=N_{\eta}\left(j\delta\right)-N_{\eta}\left(\left(j-1\right)\delta\right)$
is a Poisson random variable with rate $\eta\lambda_{b}\delta$.

Therefore,
\begin{gather*}
\E\left[\M_{\eta,2}\left(t\right)-\M_{\eta,2}^ {}\left(s\right)\Big|\F_{s}\right]\\
=\E\left[A\right]-\E\left[\lambda_{b}\int_{s}^{t}H\left(M_{\eta}^{B,R}\left(u\right)\right)du\right]\\
=\E\left[\int_{s}^{t}H\left(M_{\eta}^{B,R}\left(a_{j}-\right)\right)\lambda_{b}dt\right]-\E\left[\lambda_{b}\int_{s}^{t}H\left(M_{\eta}^{B,R}\left(u\right)\right)du\right]\\
=0
\end{gather*}
This completes the proof. 
\end{proof}

The lemma below will be used to prove both centered bail
component convergence in probability in Lemma \ref{lem:BCP}, and
centered return component convergence in probability in Lemma \ref{lem:RCP}.

\begin{restatable}[label={lem:GB1}]{lemma}{GB1Lemma}
\label{lem:GB1}Consider an arbitrary $\eta$-dependent process $M_{1}\left(t\right)$
defined as
\begin{align*}
M_{1}\left(t\right) & :=\sum_{j=1}^{N_{\eta}\left(t\right)}J_{1}\left(a_{j}\right)-J_{2}\left(a_{j}\right),
\end{align*}
where both $\left\{ J_{1}\left(t\right)\right\} _{t\geq0}$ and $\left\{ J_{2}\left(t\right)\right\} _{t\geq0}$
are $\F_{t-}$-adapted processes, $\left\{ a_{j}\right\} $ are the
arrival times of an arbitrary Poisson process $\left\{ N_{\eta}\left(t\right)\right\} _{t\geq0}$
with rate $\eta\lambda$, $t$ is an arbitrary end-time, and $\eta$
is a scaling constant. For some set of i.i.d jump sizes $\left\{ b_{j}\right\} $,
assume
\begin{itemize}
\item $0\leq J_{1}\left(a_{j}\right)\leq b_{j}$, $0\leq J_{2}\left(a_{j}\right)\leq\E\left[b_{j}\right],J_{1}(t)=J_{2}(t)=0\text{\ensuremath{\text{ for}} }t\neq\left\{ a_{j}\right\} ,$ 
\item $\E\left[J_{1}\left(t\right)-J_{2}\left(t\right)\Bigg|\F_{t^{'}}\right]=0\text{ for any }t^{'}<t,$
\item for any fixed $t$, the value of the jump size processes $J_{1}\left(t\right)$
and $J_{2}\left(t\right)$ are independent of whether an arrival actually
happens at the time $\left\{ t=a_{j}\text{ for some }j\right\} $. 
\end{itemize}
Then, 
\begin{gather*}
\E\left[M_{1}\left(t\right)^{2}\right]\leq2\lambda\eta t\E\left[b_{j}^{2}\right].
\end{gather*}

\end{restatable}
\begin{proof}
Observe,

\begin{align*}
\left(M_{1}\left(t\right)\right)^{2} & =\left(\sum_{j=1}^{N_{\eta}\left(t\right)}J_{1}\left(a_{j}-\right)-J_{2}\left(a_{j}-\right)\right)^{2}\\
 & =\left(\int_{0}^{t}(J_{1}\left(s\right)-J_{2}\left(s\right))dN_{\eta}\left(s\right)\right)^{2}\\
 & =\left(\int_{0}^{t}J_{1}\left(s\right)-J_{2}\left(s\right)dN_{\eta}\left(s\right)\right)\left(\int_{0}^{t}J_{1}\left(s\right)-J_{2}\left(s\right)dN_{\eta}\left(s\right)\right)\\
 & =\int_{0}^{t}\left(\int_{0}^{t}J_{1}\left(u\right)-J_{2}\left(u\right)dN_{\eta}\left(u\right)\right)\left(J_{1}\left(s\right)-J_{2}\left(s\right)\right)dN_{\eta}\left(s\right)\\
 & =\int_{0}^{t}\left[\int_{0}^{t}\left(J_{1}\left(u\right)-J_{2}\left(u\right)\right)\left(J_{1}\left(s\right)-J_{2}\left(s\right)\right)dN_{\eta}\left(u\right)\right]dN_{\eta}\left(s\right)\text{ (Double integral over full square \ensuremath{\left[0,t\right]^{2})}}
\end{align*}

Notice we can partition the square $\left[0,t\right]^{2}$ into the
upper triangle $u>s$, the lower triangle $u<s,$ and the diagonal
$u=s$. Note that the diagonal may not vanish due to $N_{\eta}\left(t\right)$
being a jump process.

\begin{align*}
\left(M_{1}\left(t\right)\right)^{2} & =\int_{0}^{t}\left[\int_{0}^{s-}\left(J_{1}\left(u\right)-J_{2}\left(u\right)\right)\left(J_{1}\left(s\right)-J_{2}\left(s\right)\right)dN_{\eta}\left(u\right)\right]dN_{\eta}\left(s\right)\text{ (Restrict to lower triangle: \ensuremath{u<s})}\\
 & \quad+\int_{0}^{t}\left[\int_{s}^{t}\left(J_{1}\left(u\right)-J_{2}\left(u\right)\right)\left(J_{1}\left(s\right)-J_{2}\left(s\right)\right)dN_{\eta}\left(u\right)\right]dN_{\eta}\left(s\right)\text{ (Add upper triangle: \ensuremath{u\geq s})}\\
 & =\int_{0}^{t}\left[\int_{0}^{s-}\left(J_{1}\left(u\right)-J_{2}\left(u\right)\right)\left(J_{1}\left(s\right)-J_{2}\left(s\right)\right)dN_{\eta}\left(u\right)\right]dN_{\eta}\left(s\right)\text{}\\
 & \quad+\int_{0}^{t}\left[\int_{0}^{u}\left(J_{1}\left(u\right)-J_{2}\left(u\right)\right)\left(J_{1}\left(s\right)-J_{2}\left(s\right)\right)dN_{\eta}\left(s\right)\right]dN_{\eta}\left(u\right)\text{ (Change orders of integration)}\\
 & =\int_{0}^{t}\left[\int_{0}^{s-}\left(J_{1}\left(u\right)-J_{2}\left(u\right)\right)\left(J_{1}\left(s\right)-J_{2}\left(s\right)\right)dN_{\eta}\left(u\right)\right]dN_{\eta}\left(s\right)\text{ }\\
 & \quad+\int_{0}^{t}\left[\int_{0}^{s-}\left(J_{1}\left(u\right)-J_{2}\left(u\right)\right)\left(J_{1}\left(s\right)-J_{2}\left(s\right)\right)dN_{\eta}\left(u\right)\right]dN_{\eta}\left(s\right)\text{ (Switch \ensuremath{u,s} dummies)}\\
 & \quad+\int_{0}^{t}\left[\int_{s}^{s}\left(J_{1}\left(u\right)-J_{2}\left(u\right)\right)\left(J_{1}\left(s\right)-J_{2}\left(s\right)\right)dN_{\eta}\left(u\right)\right]dN_{\eta}\left(s\right)\text{(Non-vanishing diagonal )}\\
 & =2\int_{0}^{t}\left[\int_{0}^{s-}\left(J_{1}\left(u\right)-J_{2}\left(u\right)\right)\left(J_{1}\left(s\right)-J_{2}\left(s\right)\right)dN_{\eta}\left(u\right)\right]dN_{\eta}\left(s\right)\text{ (Combine first two terms)}\\
 & \quad+\int_{0}^{t}\left[\int_{s}^{s}\left(J_{1}\left(u\right)-J_{2}\left(u\right)\right)\left(J_{1}\left(s\right)-J_{2}\left(s\right)\right)dN_{\eta}\left(u\right)\right]dN_{\eta}\left(s\right)\text{ }
\end{align*}

To manipulate this stochastic integral, denote $\Delta_{i}N_{\eta}=N_{\eta}\left(i\delta\right)-N_{\eta}\left(\left(i-1\right)\delta\right)$
for a very small $\delta$. Then, by definition, the first term of
the double triangle's is

{\scriptsize
\begin{gather*}
2\E\left[\int_{0}^{t}\left[\int_{0}^{s-}\left(J_{1}\left(u\right)-J_{2}\left(u\right)\right)\left(J_{1}\left(s\right)-J_{2}\left(s\right)\right)dN_{\eta}\left(u\right)\right]dN_{\eta}\left(s\right)\right]\\
=2\E\left[\lim_{\delta\rightarrow0}\sum_{i=1}^{\left\lfloor \frac{t}{\delta}\right\rfloor }\sum_{j=1}^{i-1}\left(J_{1}\left(j\delta\right)-J_{2}\left(j\delta\right)\right)\left(J_{1}\left(i\delta\right)-J_{2}\left(i\delta\right)\right)\Delta_{j}N_{\eta}\Delta_{i}N_{\eta}\right]\\
=2\lim_{\delta\rightarrow0}\sum_{i=1}^{\left\lfloor \frac{t}{\delta}\right\rfloor }\sum_{j=1}^{i-1}\E\left[\left(J_{1}\left(j\delta\right)-J_{2}\left(j\delta\right)\right)\left(J_{1}\left(i\delta\right)-J_{2}\left(i\delta\right)\right)\Delta_{j}N_{\eta}\Delta_{i}N_{\eta}\right]\\
=2\lim_{\delta\rightarrow0}\sum_{i=1}^{\left\lfloor \frac{t}{\delta}\right\rfloor }\sum_{j=1}^{i-1}\E\left[\E\left[\left(J_{1}\left(j\delta\right)-J_{2}\left(j\delta\right)\right)\left(J_{1}\left(i\delta\right)-J_{2}\left(i\delta\right)\right)\Delta_{j}N_{\eta}\Delta_{i}N_{\eta}\Bigg|\F_{j\delta}\right]\right]\\
=2\lim_{\delta\rightarrow0}\sum_{i=1}^{\left\lfloor \frac{t}{\delta}\right\rfloor }\sum_{j=1}^{i-1}\E\left[\E\left[\left(J_{1}\left(i\delta\right)-J_{2}\left(i\delta\right)\right)\Delta_{i}N_{\eta}\Bigg|\F_{j\delta}\right]\left(J_{1}\left(j\delta\right)-J_{2}\left(j\delta\right)\right)\Delta_{j}N_{\eta}\right]\\
=2\lim_{\delta\rightarrow0}\sum_{i=1}^{\left\lfloor \frac{t}{\delta}\right\rfloor }\sum_{j=1}^{i-1}\E\left[\underbrace{\E\left[\left(J_{1}\left(i\delta\right)-J_{2}\left(i\delta\right)\right)\Bigg|\F_{j\delta}\right]}_{\text{}=0\text{ by assumption.}}\E\left[\Delta_{i}N_{\eta}\Bigg|\F_{j\eta\delta}\right]\left(J_{1}\left(j\eta\delta\right)-J_{2}\left(j\eta\delta\right)\right)\Delta_{j}N_{\eta}\right],
\end{gather*}
}where the last step follows from the fact that
$\Delta_{i}N_{\eta}$ counts the arrivals within $N_{\eta}\left(i\delta\right)-N_{\eta}\left(\left(i-1\right)\delta\right)$
and, by assumption, these only depends on the arrival times, which
are independent of jump sizes $J_{1}\left(i\delta\right)$ and $J_{2}\left(i\delta\right)$
for small enough $\delta$.

For the second term, 
\begin{gather*}
\E\left[\int_{0}^{t}\left[\int_{s}^{s}\left(J_{1}\left(u\right)-J_{2}\left(u\right)\right)\left(J_{1}\left(s\right)-J_{2}\left(s\right)\right)dN_{\eta}\left(u\right)\right]dN_{\eta}\left(s\right)\right]\\
=\lim_{\delta\rightarrow0}\sum_{i=1}^{\left\lfloor \frac{t}{\delta}\right\rfloor }\sum_{j=i}^{i}\E\left[\left(J_{1}\left(j\delta\right)-J_{2}\left(j\delta\right)\right)\left(J_{1}\left(i\delta\right)-J_{2}\left(i\delta\right)\right)\Delta_{j}N_{\eta}\Delta_{i}N_{\eta}\right]\\
=\lim_{\delta\rightarrow0}\sum_{i=1}^{\left\lfloor \frac{t}{\delta}\right\rfloor }\E\left[\left(J_{1}\left(i\delta\right)-J_{2}\left(i\delta\right)\right)^{2}\left(\Delta_{i}N_{\eta}\right)^{2}\right]\\
=\lim_{\delta\rightarrow0}\sum_{i=1}^{\left\lfloor \frac{t}{\delta}\right\rfloor }\E\left[\E\left[\left(J_{1}\left(i\delta\right)-J_{2}\left(i\delta\right)\right)^{2}\left(\Delta_{i}N_{\eta}\right)^{2}\Bigg|\F_{\left(i-1\right)\delta}\right]\right]
\end{gather*}

By assumption, $\left(\Delta_{i}N_{\eta}\right)^{2}$ is independent
of $\left(J_{1}\left(i\delta\right)-J_{2}\left(i\delta\right)\right)^{2}$
. Therefore,
\begin{gather*}
\E\left[\int_{0}^{t}\left[\int_{s}^{s}\left(J_{1}\left(u\right)-J_{2}\left(u\right)\right)\left(J_{1}\left(s\right)-J_{2}\left(s\right)\right)dN_{\eta}\left(u\right)\right]dN_{\eta}\left(s\right)\right]\\
=\lim_{\delta\rightarrow0}\sum_{i=1}^{\left\lfloor \frac{t}{\delta}\right\rfloor }\E\left[\E\left[\left(J_{1}\left(i\delta\right)-J_{2}\left(i\delta\right)\right)^{2}\Bigg|\F_{\left(i-1\right)\delta}\right]\E\left[\left(\Delta_{i}N_{\eta}\right)^{2}\Bigg|\F_{\left(i-1\right)\delta}\right]\right]\\
=\lim_{\delta\rightarrow0}\sum_{i=1}^{\left\lfloor \frac{t}{\delta}\right\rfloor }\E\left[\E\left[\left(J_{1}\left(i\delta\right)-J_{2}\left(i\delta\right)\right)^{2}\Bigg|\F_{\left(i-1\right)\delta}\right]\right]\lambda\eta\delta\\
=\lim_{\delta\rightarrow0}\sum_{i=1}^{\left\lfloor \frac{t}{\delta}\right\rfloor }\E\left[\left(J_{1}\left(i\delta\right)-J_{2}\left(i\delta\right)\right)^{2}\right]\\
\leq\lim_{\delta\rightarrow0}\sum_{i=1}^{\left\lfloor \frac{t}{\delta}\right\rfloor }\E\left[\left(J_{1}\left(i\delta\right)\right)^{2}+\left(J_{2}\left(i\delta\right)\right)^{2}\right]\lambda\eta\delta\\
\text{(since \ensuremath{J_{1}} and \ensuremath{J_{2}} are non-negative)}\\
\leq\lambda\eta\lim_{\delta\rightarrow0}\delta\sum_{i=1}^{\left\lfloor \frac{t}{\delta}\right\rfloor }\E\left[\left(\underbrace{J_{1}\left(i\delta\right)}_{\leq b_{N_{\eta}\left(i\delta\right)}}\right)^{2}\right]+\E\left[\left(\underbrace{J_{2}\left(i\delta\right)}_{\leq\E\left[b_{N_{\eta}\left(i\delta\right)}\right]}\right)^{2}\right]\\
\leq\lambda\eta\lim_{\delta\rightarrow0}\delta\sum_{i=1}^{\left\lfloor \frac{t}{\delta}\right\rfloor }\E\left[b_{N_{\eta}\left(i\delta\right)}^{2}\right]+\left(\E\left[b_{N_{\eta}\left(i\delta\right)}\right]\right)^{2}\text{ (by assumtions of \ensuremath{J_{1}} and \ensuremath{J_{2}})}\\
=\lambda\eta\int_{0}^{t}\left(\E\left[b_{}^{2}\right]+\left(\E\left[b_{}\right]\right)^{2}\right)ds\\
=\lambda\eta t\left(\E\left[b^{2}\right]+\left(\E\left[b\right]\right)^{2}\right)\\
\leq\lambda\eta t\left(\E\left[b^{2}\right]+\E\left[b_{}^{2}\right]\right)\\
=2\lambda\eta t\E\left[b_{}^{2}\right]
\end{gather*}
Where the second-to-last line used Jensen's inequality, and the generic
random variable $b$ has the same distribution as the jump sizes $\left\{ b_{j}\right\} $.
This completes the proof.
\end{proof}

The following lemma will be used to prove both centered bail component convergence in probability in Lemma \ref{lem:BCP}, and centered return component convergence in probability in Lemma \ref{lem:RCP}.

\begin{restatable}[label={lem:GB2}]{lemma}{GB2Lemma}
\label{lem:GB2}
Consider an arbitrary $\eta$-dependent process $M_{2}\left(t\right)$
defined as

\begin{align*}
M_{2}\left(t\right) & :=\text{\ensuremath{\frac{1}{\eta}\int_{0}^{t}J_{2}\left(u\right)}}dN_{\eta}\left(u\right)-\lambda_{b}\int_{0}^{t}J_{2}\left(u\right)du
\end{align*}
 where $\left\{ J_{2}\left(t\right)\right\} _{t\geq0}$ is a $\F_{t-}$
adapted process, the arbitrary Poisson process $\left\{ N_{\eta}\left(t\right)\right\} _{t\geq0}$
has arrival rate $\eta\lambda$, $t$ is an arbitrary end-time, and
$\eta$ is a scaling constant. For some set of i.i.d jump sizes $\left\{ b_{j}\right\} $,
assume
\begin{itemize}
\item $0\leq J_{2}\left(a_{j}\right)\leq\E\left[b_{j}\right],J_{2}(t)=0\text{\ensuremath{\text{ for}} }t\neq\left\{ a_{j}\right\} ,$ 

\item for any fixed $t$, the value of the jump size processes $J_{2}\left(t\right)$ are independent of whether an arrival actually
happens at the time $\left\{ t=a_{j}\text{ for some }j\right\} $. 
\end{itemize}
Then ,

\[
\E\left[\left(M_{2}\left(t\right)\right)^{2}\right]\leq\left(\E\left[b_{j}\right]\right)^{2}\frac{\lambda t}{\eta}
\]

\end{restatable}

\begin{proof}
Notice 
\begin{align*}
\left(M_{2}\left(t\right)\right)^{2} & =\left(\text{\ensuremath{\frac{1}{\eta}\int_{0}^{t}J_{2}\left(u\right)}}dN_{\eta}\left(u\right)-\lambda\int_{0}^{t}J_{2}\left(u\right)du\right)^{2}\\
 & =\left(\ensuremath{\int_{0}^{t}\frac{1}{\eta}J_{2}\left(u\right)}dN_{\eta}\left(u\right)-\lambda J_{2}\left(u\right)du\right)^{2}\text{(Linearity of the stochastic integral)}\\
 & =\left(\ensuremath{\int_{0}^{t}J_{2}\left(u\right)}\left(\frac{1}{\eta}dN_{\eta}\left(u\right)-\lambda du\right)\right)^{2}\\
 & =\ensuremath{\int_{0}^{t}\ensuremath{\int_{0}^{t}J_{2}\left(s\right)}\left(\frac{1}{\eta}dN_{\eta}\left(s\right)-\lambda ds\right)J_{2}\left(u\right)}\left(\frac{1}{\eta}dN_{\eta}\left(u\right)-\lambda du\right)\text{ (Double integral over full square \ensuremath{\left[0,t\right]^{2})}}
\end{align*}
Observe we can partition the square $\left[0,t\right]^{2}$ into the
upper triangle $u>s$, the lower triangle $u<s,$ and the diagonal
$u=s$. Note that the diagonal may not vanish due to $N_{\eta}\left(t\right)$
being a jump process.

\begin{align*}
\left(M_{2}\left(t\right)\right)^{2} & =\ensuremath{\int_{0}^{t}\ensuremath{\int_{0}^{u-}J_{2}\left(s\right)}\left(\frac{1}{\eta}dN_{\eta}\left(s\right)-\lambda ds\right)J_{2}\left(u\right)}\left(\frac{1}{\eta}dN_{\eta}\left(u\right)-\lambda du\right)\text{ (Lower triangle)}\\
 & \quad+\ensuremath{\int_{0}^{t}\ensuremath{\int_{u}^{t}J_{2}\left(s\right)}\left(\frac{1}{\eta}dN_{\eta}\left(s\right)-\lambda ds\right)J_{2}\left(u\right)}\left(\frac{1}{\eta}dN_{\eta}\left(u\right)-\lambda du\right)\text{ (Upper triangle)}\\
 & =\ensuremath{\int_{0}^{t}\ensuremath{\int_{0}^{u-}J_{2}\left(s\right)}\left(\frac{1}{\eta}dN_{\eta}\left(s\right)-\lambda ds\right)J_{2}\left(u\right)}\left(\frac{1}{\eta}dN_{\eta}\left(u\right)-\lambda du\right)\\
 & \quad+\int_{0}^{t}\ensuremath{\int_{0}^{s}J_{2}\left(u\right)}\left(\frac{1}{\eta}dN_{\eta}\left(u\right)-\lambda du\right)J_{2}\left(s\right)\left(\frac{1}{\eta}dN_{\eta}\left(s\right)-\lambda ds\right)\text{ (Change orders of integration)}\\
 & =\ensuremath{\int_{0}^{t}\ensuremath{\int_{0}^{u-}J_{2}\left(s\right)}\left(\frac{1}{\eta}dN_{\eta}\left(s\right)-\lambda ds\right)J_{2}\left(u\right)}\left(\frac{1}{\eta}dN_{\eta}\left(u\right)-\lambda du\right)\\
 & \quad+\ensuremath{\int_{0}^{t}\ensuremath{\int_{0}^{u-}J_{2}\left(s\right)}\left(\frac{1}{\eta}dN_{\eta}\left(s\right)-\lambda ds\right)J_{2}\left(u\right)}\left(\frac{1}{\eta}dN_{\eta}\left(u\right)-\lambda du\right)\text{\text{ (Switch \ensuremath{u,s} dummies)}}\\
 & \quad+\ensuremath{\int_{0}^{t}\ensuremath{\int_{u}^{u}J_{2}\left(s\right)}\left(\frac{1}{\eta}dN_{\eta}\left(s\right)-\lambda ds\right)J_{2}\left(u\right)}\left(\frac{1}{\eta}dN_{\eta}\left(u\right)-\lambda du\right)\text{\text{(Non-vanishing diagonal )}}\\
 & =\ensuremath{2\int_{0}^{t}\ensuremath{\int_{0}^{u-}J_{2}\left(s\right)}\left(\frac{1}{\eta}dN_{\eta}\left(s\right)-\lambda_{b}ds\right)J_{2}\left(u\right)}\left(\frac{1}{\eta}dN_{\eta}\left(u\right)-\lambda du\right)\text{\text{ (Combine first two terms)}}\\
 & \quad+\ensuremath{\int_{0}^{t}\ensuremath{\int_{u}^{u}J_{2}\left(s\right)}\left(\frac{1}{\eta}dN_{\eta}\left(s\right)-\lambda ds\right)J_{2}\left(u\right)}\left(\frac{1}{\eta}dN_{\eta}\left(u\right)-\lambda du\right)
\end{align*}

Denote the change $\Delta_{i}N_{\eta}:=N_{\eta}\left(i\delta\right)-N_{\eta}\left(\left(i-1\right)\delta\right)$.
Both $du$ and $dt$ are approximated with the length $\delta$.

Then for the first term,
\begin{gather*}
\E\left[\ensuremath{2\int_{0}^{t}\ensuremath{\int_{0}^{u-}J_{2}\left(s\right)}\left(\frac{1}{\eta}dN_{\eta}\left(s\right)-\lambda ds\right)J_{2}\left(u\right)}\left(\frac{1}{\eta}dN_{\eta}\left(u\right)-\lambda du\right)\right]\\
=2\E\left[\lim_{\delta\rightarrow0}\sum_{i=1}^{\left\lfloor \frac{t}{\delta}\right\rfloor }\sum_{j=1}^{i-1}J_{2}\left(j\delta\right)\left(\frac{1}{\eta}\Delta_{j}N_{\eta}-\lambda\delta\right)J_{2}\left(i\delta\right)\left(\frac{1}{\eta}\Delta_{i}N_{\eta}-\lambda\delta\right)\right]\\
=2\lim_{\delta\rightarrow0}\E\left[\sum_{i=1}^{\left\lfloor \frac{t}{\delta}\right\rfloor }\sum_{j=1}^{i-1}J_{2}\left(j\delta\right)\left(\frac{1}{\eta}\Delta_{j}N_{\eta}-\lambda\delta\right)J_{2}\left(i\delta\right)\left(\frac{1}{\eta}\Delta_{i}N_{\eta}-\lambda\delta\right)\right]\\
=2\lim_{\delta\rightarrow0}\sum_{i=1}^{\left\lfloor \frac{t}{\delta}\right\rfloor }\sum_{j=1}^{i-1}\E\left[J_{2}\left(j\delta\right)\left(\frac{1}{\eta}\Delta_{j}N_{\eta}-\lambda\delta\right)J_{2}\left(i\delta\right)\left(\frac{1}{\eta}\Delta_{i}N_{\eta}-\lambda\delta\right)\right]\\ \displaybreak[1]
=2\lim_{\delta\rightarrow0}\sum_{i=1}^{\left\lfloor \frac{t}{\delta}\right\rfloor }\sum_{j=1}^{i-1}\E\left[\E\left[J_{2}\left(j\delta\right)\left(\frac{1}{\eta}\Delta_{j}N_{\eta}-\lambda\delta\right)J_{2}\left(i\delta\right)\left(\frac{1}{\eta}\Delta_{i}N_{\eta}-\lambda\delta\right)\Bigg|\F_{j\delta}\right]\right]\\
=2\lim_{\delta\rightarrow0}\sum_{i=1}^{\left\lfloor \frac{t}{\delta}\right\rfloor }\sum_{j=1}^{i-1}\E\left[J_{2}\left(j\delta\right)\left(\frac{1}{\eta}\Delta_{j}N_{\eta}-\lambda\delta\right)\E\left[J_{2}\left(i\delta\right)\left(\frac{1}{\eta}\Delta_{j}N_{\eta}-\lambda\delta\right)\Bigg|\F_{j\delta}\right]\right]\\
=2\lim_{\delta\rightarrow0}\sum_{i=1}^{\left\lfloor \frac{t}{\delta}\right\rfloor }\sum_{j=1}^{i-1}\E\left[J_{2}\left(j\delta\right)\left(\frac{1}{\eta}\Delta_{j}N_{\eta}-\lambda\delta\right)\E\left[J_{2}\left(i\delta\right)\Bigg|\F_{i\delta}\right]\E\left[\underbrace{\left(\frac{1}{\eta}\Delta_{i}N_{\eta}-\lambda_{b}\delta\right)}_{=0\text{ by definition}}\Bigg|\F_{j\delta}\right]\right]\\
=0
\end{gather*}
where the second-to-last line follows from the fact that $\Delta_{i}N_{\eta}$
counts the arrival within $N_{\eta}\left(i\delta\right)-N_{\eta}\left(\left(i-1\right)\delta\right)$.
By assumption, these only depends on the arrival times, which are
independent of jump sizes $J_{2}\left(i\delta\right)$ for small enough
$\delta$. 

For the second term,
\begin{gather*}
\E\left[\ensuremath{\int_{0}^{t}\ensuremath{\int_{u}^{u}J_{2}\left(s\right)}\left(\frac{1}{\eta}dN_{\eta}\left(s\right)-\lambda ds\right)J_{2}\left(u\right)}\left(\frac{1}{\eta}dN_{\eta}\left(u\right)-\lambda du\right)\right]\\
=\E\left[\lim_{\delta\rightarrow0}\sum_{i=1}^{\left\lfloor \frac{t}{\delta}\right\rfloor }\sum_{j=i}^{i}J_{2}\left(j\delta\right)\left(\frac{1}{\eta}\Delta_{j}N_{\eta}-\lambda\delta\right)J_{2}\left(i\delta\right)\left(\frac{1}{\eta}\Delta_{i}N_{\eta}-\lambda\delta\right)\right]\\
=\E\left[\lim_{\delta\rightarrow0}\sum_{i=1}^{\left\lfloor \frac{t}{\delta}\right\rfloor }\left(J_{2}\left(i\delta\right)\right)^{2}\left(\frac{1}{\eta}\Delta_{i}N_{\eta}-\lambda\delta\right)^{2}\right]\\ \displaybreak[1]
=\lim_{\delta\rightarrow0}\sum_{i=1}^{\left\lfloor \frac{t}{\delta}\right\rfloor }\E\left[\left(J_{2}\left(i\delta\right)\right)^{2}\left(\frac{1}{\eta}\Delta_{i}N_{\eta}-\lambda\delta\right)^{2}\right]\\
=\lim_{\delta\rightarrow0}\sum_{i=1}^{\frac{t}{\delta}}\E\left[\E\left[\left(J_{2}\left(i\delta\right)\right)^{2}\left(\frac{1}{\eta}\Delta_{i}N_{\eta}-\lambda\delta\right)^{2}\Bigg|\F_{\left(i-1\right)\delta}\right]\right]\\
=\lim_{\delta\rightarrow0}\sum_{i=1}^{\frac{t}{\delta}}\E\left[\E\left[\left(J_{2}\left(i\delta\right)\right)^{2}\Bigg|\F_{\left(i-1\right)\delta}\right]\text{\ensuremath{\E\left[\left(\frac{1}{\eta}\Delta_{i}N_{\eta}-\lambda\delta\right)^{2}\Bigg|\F_{\left(i-1\right)\delta}\right]}}\right]\\
=\lim_{\delta\rightarrow0}\sum_{i=1}^{\frac{t}{\delta}}\E\left[\E\left[\left(J_{2}\left(i\delta\right)\right)^{2}\Bigg|\F_{\left(i-1\right)\delta}\right]\text{\ensuremath{\E\left[\left(\frac{1}{\eta}\Delta_{i}N_{\eta}-\lambda\delta\right)^{2}\right]}}\right]\\
=\lim_{\delta\rightarrow0}\sum_{i=1}^{\frac{t}{\delta}}\E\left[\left(\underbrace{J_{2}\left(i\delta\right)}_{\leq\E\left[b_{N_{\eta}(i\delta)}\right]}\right)^{2}\right]\E\left[\left(\frac{1}{\eta}\Delta_{i}N_{\eta}-\lambda\delta\right)^{2}\right]\\
\leq\lim_{\delta\rightarrow0}\sum_{i=1}^{\frac{t}{\delta}}\left(\E\left[b\right]\right)^{2}\E\left[\left(\frac{1}{\eta}\Delta_{i}N_{\eta}-\lambda\delta\right)^{2}\right]
\end{gather*}
Because of the property of Poisson process, $\Delta_{j}N_{\eta}=N_{\eta}\left(i\delta\right)-N\left(\left(i-1\right)\delta\right)$
is a Poisson random variable with rate $\lambda_{b}\eta\delta$. And
the generic random variable $b$ has the same distribution as the
jump sizes $\left\{ b_{j}\right\} $. Then,

\begin{align*}
\E\left[\left(\frac{1}{\eta}\Delta_{i}N_{\eta}-\lambda\delta\right)^{2}\right] & =\E\left[\frac{1}{\eta^{2}}\left(\Delta_{i}N_{\eta}\right)^{2}-\frac{2\lambda}{\eta^{2}}\Delta_{i}N_{\eta}\delta+\lambda^{2}\delta^{2}\right]\\
 & =\frac{1}{\eta^{2}}\left(\var\left(\Delta_{i}N_{\eta}\right)-\E\left[\Delta_{i}N_{\eta}\right]^{2}\right)-\frac{2\lambda}{\eta^{2}}\E\left[\Delta_{i}N_{\eta}\right]\delta+\lambda^{2}\delta^{2}\\
 & =\frac{1}{\eta^{2}}\left(\lambda\eta\delta-\lambda^{2}\eta^{2}\delta^{2}\right)-\frac{2\lambda_{b}}{\eta^{2}}\lambda\eta\delta^{2}+\lambda^{2}\delta^{2}
\end{align*}
So, we have, 
\begin{align*}
\E\left[\left(\frac{1}{\eta}\Delta_{i}N_{\eta}-\lambda\delta\right)^{2}\right] & =\frac{\lambda\delta}{\eta}
\end{align*}
Therefore,

\begin{align*}
\E\left[\left(M_{2}\left(t\right)\right)^{2}\right] & \leq\lim_{\delta\rightarrow0}\sum_{i=1}^{\left\lfloor \frac{t}{\delta}\right\rfloor }\left(\E\left[b\right]\right)^{2}\E\left[\left(\frac{1}{\eta}\Delta_{i}N-\lambda\delta\right)^{2}\right]\\
 & =\lim_{\delta\rightarrow0}\sum_{i=1}^{\left\lfloor \frac{t}{\delta}\right\rfloor }\left(\E\left[b\right]\right)^{2}\frac{\lambda}{\eta}\delta\\
 & =\int_{0}^{t}\left(\E\left[b\right]\right)^{2}\frac{\lambda}{\eta}ds\\
 & =\left(\E\left[b\right]\right)^{2}\frac{\lambda t}{\eta}
\end{align*}
 This completes the proof.
\end{proof}

\subsubsection{\label{app:ProbConvergence}Proof of the Lemma \ref{lem:ProbConvergence}}
\ProbConvergenceLemma*
\begin{proof}
We define the following events for some $\epsilon>0$,

\[
\A_{\eta}:=\left\{ \sup_{0\leq u\leq T}\left|M_{\eta}\left(u\right)\right|>\epsilon\right\} 
\]

Notice that, since $M\left(u\right)$ is a martingale and $x^{2}$
is a convex function, $\left|M_{\eta}\left(T\right)\right|^{2}$is
a submartingale. Then, by the Doob's supremum inequality for continuous-submartingales.

\begin{align*}
\P\left(\A_{\eta}\right) & =\P\left\{ \sup_{0\leq u\leq T}\left|M_{\eta}\left(u\right)\right|>\epsilon\right\} \\
 & \leq\frac{\E\left[\left|M_{\eta}\left(T\right)\right|^{2}\right]}{\epsilon^{2}}\\
 & =\frac{\E\left[\left(M_{\eta}\left(T\right)\right)^{2}\right]}{\epsilon^{2}}
\end{align*}

Then we have

\begin{align*}
\P\left(\A_{\eta}\right) & \leq\frac{\E\left[\left(M_{\eta}\left(T\right)\right)^{2}\right]}{\epsilon^{2}}\\
 & \leq\frac{C}{\epsilon^{2}\eta}
\end{align*}
where last line follows from the assumption that $\E\left[\left(M_{\eta}\left(T\right)\right)^{2}\right]\leq\frac{C}{\eta}$.
This suffices to show $\sup_{0\leq u\leq T}\left|M_{\eta}\left(u\right)\right|\underset{P}{\rightarrow}0.$ 
\end{proof}

\subsubsection{\label{app:BCP}Proof of the Lemma \ref{lem:BCP}.}

\BCPLemma*

\begin{proof}
Notice $\M_{b,\eta}^{B,R}\left(t\right)$ can be decomposed into martingales
in Lemma \ref{lem:BMG1} and Lemma \ref{lem:BMG3}:

\begin{align*}
\M_{\eta,b}^{B,R}\left(t\right) & =\left[\underbrace{\frac{1}{\eta}\sum_{j=1}^{N_{\eta}^{b}(t)}b_{j}\left\{ M_{\eta}^{B,R}\left(a_{j}-\right)\geq b_{j}\right\} -\text{\ensuremath{\frac{1}{\eta}}\ensuremath{\sum_{j=1}^{N_{\eta}^{b}(t))}H\left(M_{\eta}^{B,R}\left(a_{j}-\right)\right)}}}_{c_{1}^{\eta}\left(t\right)}\right]\\
 & \quad+\left[\underbrace{\text{\ensuremath{\frac{1}{\eta}}\ensuremath{\sum_{j=1}^{N_{\eta}^{b}(t)}H\left(M_{\eta}^{B,R}\left(a_{j}-\right)\right)}}-\frac{1}{\eta}\int_{0}^{t}H\left(M_{\eta}^{B,R}\left(u\right)\right)dN_{\eta}\left(u\right)}_{\text{\ensuremath{=0\text{ by definition of stochastic intergrals}}}}\right]\\
 & \quad+\left[\underbrace{\text{\ensuremath{\frac{1}{\eta}\int_{0}^{t}H\left(M_{\eta}^{B,R}\left(u\right)\right)}}dN_{\eta}\left(u\right)-\lambda_{b}\int_{0}^{t}H\left(M_{\eta}^{B,R}\left(u\right)\right)du}_{c_{2}^{\eta}\left(t\right)}\right],
\end{align*}
therefore $\M_{b,\eta}^{B,R}\left(t\right)$ is a martingale because
it is a sum of martingales.

Our goal now is to apply Lemma \ref{lem:GB1} on the term $c_{1}^{\eta}\left(t\right)$
and apply Lemma \ref{lem:GB2} on the term $c_{2}^{\eta}\left(t\right)$
for the following jump sizes 

\begin{align*}
\text{\ensuremath{J_{1}\left(t\right)}	 } & :=b_{N_{\eta}\left(t\right)}\left\{ b_{N_{\eta}\left(t\right)}\leq M_{\eta}^{B,R^{B}}\left(t-\right)\right\} \\
J_{2}\left(t\right) & :=\int_{0}^{M_{\eta}^{B,R^{B}}\left(t-\right)}xf_{b}\left(x\right)dx.
\end{align*}

Notice that $b_{N\left(t\right)}$ in $J_{1}\left(t\right)$
is known at time $t-$ , as its acceptance at time $t$ depends on
the comparison with $M_{\eta}^{B,R^{B}}\left(t-\right)$ at time $t-$.
Hence, both $\left\{ J_{1}\left(t\right)\right\} _{t\geq0}$ and $\left\{ J_{2}\left(t\right)\right\} _{t\geq0}$
are $\F_{t-}$-adapted processes. It is routine to check the assumptions
of Lemma \ref{lem:GB1} for these choices of $J_{1}\left(t\right)$
and $J_{2}\left(t\right)$ are satisfied. We now rewrite $c_{1}^{\eta}\left(t\right)$
as follows:

\begin{align*}
c_{1}^{\eta}\left(t\right) & =\frac{1}{\eta}\sum_{j=1}^{N_{\eta}^{b}(t)}b_{j}\left\{ M_{\eta}^{B,R}\left(a_{j}-\right)\geq b_{j}\right\} -\text{\ensuremath{\frac{1}{\eta}}\ensuremath{\sum_{j=1}^{N_{\eta}^{b}(t))}H\left(M_{\eta}^{B,R}\left(a_{j}-\right)\right)}}\\
 & =\frac{1}{\eta}\sum_{j=1}^{N_{\eta}^{b}(t)}b_{N_{\eta}\left(a_{j}\right)}\left\{ b_{N_{\eta}\left(a_{j}\right)}\leq M_{\eta}^{B,R^{B}}\left(a_{j}-\right)\right\} -\text{\ensuremath{\frac{1}{\eta}}\ensuremath{\sum_{j=1}^{N_{\eta}^{b}(t))}H\left(M_{\eta}^{B,R}\left(a_{j}-\right)\right)}}\\
 & =\frac{1}{\eta}\sum_{j=1}^{N_{\eta}\left(t\right)}J_{1}\left(a_{j}-\right)-J_{2}\left(a_{j}-\right).
\end{align*}
We can now apply Lemma \ref{lem:GB1} on the term $c_{1}^{\eta}\left(t\right)$
to generate a upper bound for the second moment :

\[
\E\left[\left(c_{1}^{\eta}\left(t\right)\right)^{2}\right]\leq\frac{1}{\eta^{2}}\cdot\left(2\lambda_{b}\eta t\E\left[b_{i}^{2}\right]\right)=\frac{2\lambda_{b}t\E\left[b_{i}^{2}\right]}{\eta}.
\]
Moreover, we can rewrite $c_{1}^{\eta}\left(t\right)$
in terms of $J_{2}\left(t\right)$ as follows:

\begin{align*}
c_{2}^{\eta}\left(t\right) & \text{\ensuremath{=\frac{1}{\eta}\int_{0}^{t}H\left(M_{\eta}^{B,R}\left(u\right)\right)}}dN_{\eta}\left(u\right)-\lambda_{b}\int_{0}^{t}H\left(M_{\eta}^{B,R}\left(u\right)\right)du\\
 & \text{\ensuremath{=\frac{1}{\eta}\int_{0}^{t}J_{2}\left(u\right)}}dN_{\eta}\left(u\right)-\lambda_{b}\int_{0}^{t}J_{2}\left(u\right)du.
\end{align*}
It is routine to check the assumptions of Lemma
\ref{lem:GB2} for the choice of $J_{2}\left(t\right)$ are satisfied.
We can now apply Lemma \ref{lem:GB2} on the term $c_{2}^{\eta}\left(t\right)$
to generate a upper bound for the second moment :

\begin{align*}
\E\left[\left(c_{2}^{\eta}\left(t\right)\right)^{2}\right] & \leq\E\left[b_{j}^{2}\right]\frac{\lambda_{b}t}{\eta}.
\end{align*}
By Lemma \ref{lem:ProbConvergence},

\[
\sup_{t\in[0,T]}\left|c_{1}^{\eta}\left(t\right)\right|\underset{P}{\rightarrow}0,
\]

\[
\sup_{t\in[0,T]}\left|c_{2}^{\eta}\left(t\right)\right|\underset{P}{\rightarrow}0.
\]

Then by the triangle inequality and the fact that the sum of terms
converging in probability also converges in probability, it follows
that:
\[
\sup_{t\in[0,T]}\left|\M_{b,\eta}^{B,R}\left(t\right)\right|\underset{P}{\rightarrow}0.
\]

This completes the proof.
\end{proof}

\subsubsection{Supporting Lemmas for the Lemma \ref{lem:RCP}}

\begin{lemma}
\label{lem:RMG1} The process 
\[
\M_{\eta,3}^{B,R}\left(t\right):=\frac{1}{\eta}\sum_{j=1}^{N_{\eta}^{b}\left(t\right)}\left(1-p_{j}\right)\cdot b_{j}\cdot\left\{ M_{\eta}^{B,R}\left(a_{j}-\right)\geq b_{j}\right\} \cdot\left\{ t>a_{j}^{\eta}+s_{j}\right\} -\frac{1}{\eta}\sum_{j=1}^{N_{\eta}^{b}\left(t\right)}\left(1-p^{*}\right)H\left(M_{\eta}^{B,R}\left(a_{j}-\right)\right)\cdot F_{s}\left(t-a_{j}^{\eta}\right)
\]
is a martingale with respect to the filtration $\left\{ \F_{s}\right\} $
which contains all bail arrivals and bail request sizes until time
$s$:

\[
\F_{s}:=\sigma\left(N_{\eta}^{b}\left(s\right),\left\{ b_{j}\right\} _{j=1}^{N_{\eta}^{b}\left(s\right)},\left\{ a_{j}\right\} _{j=1}^{N_{\eta}^{b}\left(s\right)},\left\{ s_{j}\right\} _{j=1}^{N_{\eta}^{b}\left(s\right)}\right).
\]
\end{lemma}

\begin{proof}
It suffices to check the martingale property for any $s<t$: in other
words, $\E\left[\M_{\eta,3}^{B,R}\left(t\right)-\M_{\eta,3}^{B,R}\left(s\right)\Big|\F_{s}\right]=0$.

\begin{gather*}
\E\left[\M_{\eta,3}^{B,R}\left(t\right)-\M_{\eta,3}^{B,R}\left(s\right)\Big|\F_{s}\right]\\
=\underbrace{\E\left[\frac{1}{\eta}\sum_{j=N_{\eta}^{b}\left(s\right)+1}^{N_{\eta}^{b}\left(t\right)}\left(1-p_{j}\right)\cdot b_{j}\cdot\left\{ M_{\eta}^{B,R}\left(a_{j}^{\eta}-\right)\geq b_{j}\right\} \cdot\left\{ t>a_{j}^{\eta}+s_{j}\right\} \Bigg|\F_{s}\right]}_{A}\\
-\underbrace{\E\left[\frac{1}{\eta}\sum_{j=N_{\eta}^{b}\left(s\right)+1}^{N_{\eta}^{b}\left(t\right)}\left(1-p^{*}\right)H\left(M_{\eta}^{B,R}\left(a_{j}^{\eta}-\right)\right)\cdot F_{s}\left(t-a_{j}^{\eta}\right)\Bigg|\F_{s}\right]}_{B}
\end{gather*}
We first look into term $A,$
\begin{align*}
A & =\E\left[\frac{1}{\eta}\sum_{j=N_{\eta}^{b}\left(s\right)+1}^{N_{\eta}^{b}\left(t\right)}\left(1-p_{j}\right)\cdot b_{j}\cdot\left\{ M_{\eta}^{B,R}\left(a_{j}^{\eta}-\right)\geq b_{j}\right\} \cdot\left\{ t>a_{j}^{\eta}+s_{j}\right\} \Bigg|\F_{s}\right]\\
 & =\E\left[\E\left[\frac{1}{\eta}\sum_{j=N_{\eta}^{b}\left(s\right)+1}^{N_{\eta}^{b}\left(t\right)}\left(1-p_{j}\right)\cdot b_{j}\cdot\left\{ M_{\eta}^{B,R}\left(a_{j}^{\eta}-\right)\geq b_{j}\right\} \cdot\left\{ t>a_{j}^{\eta}+s_{j}\right\} \Bigg|\F_{s},N_{\eta}^{b}\left(t\right)\right]\Bigg|\F_{s}\right]\\
 & =\E\left[\frac{1}{\eta}\sum_{j=N_{\eta}^{b}\left(s\right)+1}^{N_{\eta}^{b}\left(t\right)}\E\left[\left(1-p_{j}\right)\cdot b_{j}\cdot\left\{ M_{\eta}^{B,R}\left(a_{j}^{\eta}-\right)\geq b_{j}\right\} \cdot\left\{ t>a_{j}^{\eta}+s_{j}\right\} \Bigg|\F_{s},N_{\eta}^{b}\left(t\right)\right]\Bigg|\F_{s}\right]\\
 & =\E\left[\frac{1}{\eta}\sum_{j=N_{\eta}^{b}\left(s\right)+1}^{N_{\eta}^{b}\left(t\right)}\E\left[\left(1-p_{j}\right)\cdot b_{j}\cdot\left\{ M_{\eta}^{B,R}\left(a_{j}^{\eta}-\right)\geq b_{j}\right\} \cdot\left\{ t>a_{j}^{\eta}+s_{j}\right\} \Bigg|\F_{s},N_{\eta}^{b}\left(t\right)\right]\Bigg|\F_{s}\right]\\
 & =\E\left[\frac{1}{\eta}\sum_{j=N_{\eta}^{b}\left(s\right)+1}^{N_{\eta}^{b}\left(t\right)}\E\left[\left(1-p_{j}\right)\cdot b_{j}\cdot\left\{ M_{\eta}^{B,R}\left(a_{j}^{\eta}-\right)\geq b_{j}\right\} \cdot\left\{ t>a_{j}^{\eta}+s_{j}\right\} \Bigg|\F_{s},N_{\eta}^{b}\left(t\right)\right]\Bigg|\F_{s}\right]\\
 & =\E\left[\frac{1}{\eta}\sum_{j=N_{\eta}^{b}\left(s\right)+1}^{N_{\eta}^{b}\left(t\right)}\E\left[\left(1-p_{j}\right)\cdot b_{j}\cdot\left\{ M_{\eta}^{B,R}\left(a_{j}^{\eta}-\right)\geq b_{j}\right\} \cdot\left\{ t>a_{j}^{\eta}+s_{j}\right\} \Bigg|\F_{s},N_{\eta}^{b}\left(t\right)\right]\Bigg|\F_{s}\right]\\  & =\scalebox{0.8}{$\E\left[\frac{1}{\eta}\sum_{j=N_{\eta}^{b}\left(s\right)+1}^{N_{\eta}^{b}\left(t\right)}\E\left[\E\left[\left(1-p_{j}\right)\cdot b_{j}\cdot\left\{ M_{\eta}^{B,R}\left(a_{j}^{\eta}-\right)\geq b_{j}\right\} \cdot\left\{ t>a_{j}^{\eta}+s_{j}\right\} \Bigg|\F_{s},N_{\eta}^{b}\left(t\right),M_{\eta}^{B,R}\left(a_{j}^{\eta}-\right),a_{j}^{\eta}\right]\Bigg|\F_{s},N_{\eta}^{b}\left(t\right)\right]\Bigg|\F_{s}\right]$}\\
 & =\scalebox{0.8}{$
\E\left[\frac{1}{\eta}\sum_{j=N_{\eta}^{b}\left(s\right)+1}^{N_{\eta}^{b}\left(t\right)}\E\left[\E\left[\left(1-p_{j}\right)\right]\E\left[b_{j}\right]\E\left[\left\{ M_{\eta}^{B,R}\left(a_{j}^{\eta}-\right)\geq b_{j}\right\} \Bigg|\F_{s},N_{\eta}^{b}\left(t\right),M_{\eta}^{B,R}\left(a_{j}^{\eta}-\right),a_{j}^{\eta}\right]\E\left[\left\{ t>a_{j}^{\eta}+s_{j}\right\} \Bigg|a_{j}^{\eta}\right]\Bigg|\F_{s},N_{\eta}^{b}\left(t\right)\right]\Bigg|\F_{s}\right]$}\\
 & =\E\left[\frac{1}{\eta}\sum_{j=N_{\eta}^{b}\left(s\right)+1}^{N_{\eta}^{b}\left(t\right)}\left(1-p^{*}\right)H\left(M_{\eta}^{B,R}\left(a_{j}^{\eta}-\right)\right)\cdot F_{s}\left(t-a_{j}^{\eta}\right)\Bigg|\F_{s}\right]\\
 & =B.
\end{align*}

This completes the proof.
\end{proof}

\begin{lemma}
\label{lem:RMG2}The process 

\[
\M_{\eta,4}^{B,R}\left(t\right):=\frac{1}{\eta}\int_{0}^{t}\left(1-p^{*}\right)H\left(M_{\eta}^{B,R}\left(u\right)\right)\cdot F_{s}\left(t-u\right)dN_{\eta}^{b}\left(u\right)-\lambda_{b}\int_{0}^{t}\left(1-p^{*}\right)H\left(M_{\eta}^{B,R}\left(u\right)\right)\cdot F_{s}\left(t-u\right)du
\]
is a martingale with respect to the filtration $\left\{ \F_{s}\right\} $
which contains all bail arrivals and bail request sizes until time
$s$:

\[
\F_{s}:=\sigma\left(N_{\eta}^{b}\left(s\right),\left\{ b_{j}\right\} _{j=1}^{N_{\eta}^{b}\left(s\right)},\left\{ a_{j}\right\} _{j=1}^{N_{\eta}^{b}\left(s\right)},\left\{ s_{j}\right\} _{j=1}^{N_{\eta}^{b}\left(s\right)}\right).
\]
\end{lemma}
\begin{proof}
It suffices to check the Martingale property for any $s<t$: in other
words, $\E\left[\M_{\eta,4}^{B,R}\left(t\right)-\M_{\eta,4}^{B,R}\left(s\right)\Big|\F_{s}\right]=0$.

\begin{gather*}
\E\left[\M_{\eta,4}^{B,R}\left(t\right)-\M_{\eta,4}^{B,R}\left(s\right)\Big|\F_{s}\right]\\
=\underbrace{\E\left[\frac{1}{\eta}\int_{\text{\ensuremath{s}}}^{t}\left(1-p^{*}\right)H\left(M_{\eta}^{B,R}\left(u\right)\right)\cdot F_{s}\left(t-u\right)dN_{\eta}^{b}\left(u\right)\Big|\F_{s}\right]}_{A}\\
-\underbrace{\E\left[\lambda_{b}\int_{s}^{t}\left(1-p^{*}\right)H\left(M_{\eta}^{B,R}\left(u\right)\right)\cdot F_{s}\left(t-u\right)du\Big|\F_{s}\right]}_{B}
\end{gather*}
 We look into the term $A$ first.
 
\begin{align*}
A & =\E\left[\frac{1}{\eta}\int_{\text{\ensuremath{s}}}^{t}\left(1-p^{*}\right)H\left(M_{\eta}^{B,R}\left(u\right)\right)\cdot F_{s}\left(t-u\right)dN_{\eta}^{b}\left(u\right)\Big|\F_{s}\right]\\
 & =\frac{1}{\eta}\E\left[\lim_{\delta\rightarrow0}\sum_{j=\left\lfloor \frac{s}{\delta}\right\rfloor }^{\left\lfloor \frac{t}{\delta}\right\rfloor }\left(1-p^{*}\right)H\left(M_{\eta}^{B,R}\left(\left(j-1\right)\delta\right)\right)\cdot F_{s}\left(t-j\delta\right)\Delta_{j}N_{\eta}^{b}\Big|\F_{s}\right]\\
 & =\frac{1}{\eta}\lim_{\delta\rightarrow0}\sum_{j=\left\lfloor \frac{s}{\delta}\right\rfloor }^{\left\lfloor \frac{t}{\delta}\right\rfloor }\E\left[\left(1-p^{*}\right)H\left(M_{\eta}^{B,R}\left(\left(j-1\right)\delta\right)\right)\cdot F_{s}\left(t-j\delta\right)\Delta_{j}N_{\eta}^{b}\Big|\F_{s}\right]\text{}\\
 & =\frac{1}{\eta}\lim_{\delta\rightarrow0}\sum_{j=\left\lfloor \frac{s}{\delta}\right\rfloor }^{\left\lfloor \frac{t}{\delta}\right\rfloor }\E\left[\E\left[\left(1-p^{*}\right)H\left(M_{\eta}^{B,R}\left(\left(j-1\right)\delta\right)\right)\cdot F_{s}\left(t-j\delta\right)\Delta_{j}N_{\eta}^{b}\Big|\F_{s},M_{\eta}^{B,R}\left(\left(j-1\right)\delta\right)\right]\Big|\F_{s}\right]\\\displaybreak[1]
 & =\frac{1}{\eta}\lim_{\delta\rightarrow0}\sum_{j=\left\lfloor \frac{s}{\delta}\right\rfloor }^{\left\lfloor \frac{t}{\delta}\right\rfloor }\E\left[\left(1-p^{*}\right)H\left(M_{\eta}^{B,R}\left(\left(j-1\right)\delta\right)\right)\cdot F_{s}\left(t-j\delta\right)\E\left[\Delta_{j}N_{\eta}^{b}\Big|\F_{s},M_{\eta}^{B,R}\left(\left(j-1\right)\delta\right)\right]\Big|\F_{s}\right]\\
 & =\frac{1}{\eta}\lim_{\delta\rightarrow0}\sum_{j=\left\lfloor \frac{s}{\delta}\right\rfloor }^{\left\lfloor \frac{t}{\delta}\right\rfloor }\E\left[\left(1-p^{*}\right)H\left(M_{\eta}^{B,R}\left(\left(j-1\right)\delta\right)\right)\cdot F_{s}\left(t-j\delta\right)\lambda_{b}\eta\delta\Big|\F_{s}\right]\\
 & =\frac{1}{\eta}\E\left[\lim_{\delta\rightarrow0}\sum_{j=\left\lfloor \frac{s}{\delta}\right\rfloor }^{\left\lfloor \frac{t}{\delta}\right\rfloor }\left(1-p^{*}\right)H\left(M_{\eta}^{B,R}\left(\left(j-1\right)\delta\right)\right)\cdot F_{s}\left(t-j\delta\right)\lambda_{b}\eta\delta\Big|\F_{s}\right]\\
 & =\E\left[\lambda_{b}\int_{s}^{t}\left(1-p^{*}\right)H\left(M_{\eta}^{B,R}\left(u\right)\right)\cdot F_{s}\left(t-j\delta\right)du\Big|\F_{s}\right]\\
 & =B.
\end{align*}
This completes the proof.
\end{proof}

\subsubsection{\label{app:RCP}Proof of the Lemma \ref{lem:RCP}.}

\RCPLemma*

\begin{proof}
Notice $\M_{\eta,r}^{B,R}\left(t\right)$ can be decomposed into martingales
in Lemma \ref{lem:RMG1} and Lemma \ref{lem:RMG2}:

\begin{align*}
\M_{\eta,r}^{B,R}\left(t\right) & =\scalebox{0.8}{$\left[\underbrace{\frac{1}{\eta}\sum_{j=1}^{N_{\eta}^{b}\left(t\right)}\left(1-p_{j}\right)\cdot b_{j}\cdot\left\{ M_{\eta}^{B,R}\left(a_{j}^{\eta}-\right)\geq b_{j}\right\} \cdot\left\{ t>a_{j}^{\eta}+s_{j}\right\} -\frac{1}{\eta}\sum_{j=1}^{N_{\eta}^{b}\left(t\right)}\left(1-p^{*}\right)H\left(M^{B,R}\left(a_{j}-\right)\right)\cdot F_{s}\left(t-a_{j}^{\eta}\right)}_{c_{1}^{\eta}\left(t\right)}\right]$}\\ \displaybreak[1]
 & \scalebox{0.8}{$\quad+\left[\underbrace{\frac{1}{\eta}\sum_{j=1}^{N_{\eta}^{b}\left(t\right)}\left(1-p^{*}\right)H\left(M_{\eta}^{B,R}\left(a_{j}^{\eta}-\right)\right)\cdot F_{s}\left(t-a_{j}^{\eta}\right)-\frac{1}{\eta}\int_{0}^{t}\left(1-p^{*}\right)H\left(M_{\eta}^{B,R}\left(u\right)\right)\cdot F_{s}\left(t-u\right)dN_{\eta}^{b}\left(u\right)}_{0\text{ by definition of stochastic integrals }}\right]$}\\
 & \scalebox{0.8}{$\quad+\left[\underbrace{\frac{1}{\eta}\int_{0}^{t}\left(1-p^{*}\right)H\left(M_{\eta}^{B,R}\left(u\right)\right)\cdot F_{s}\left(t-u\right)dN_{\eta}^{b}\left(u\right)-\lambda_{b}\int_{0}^{t}\left(1-p^{*}\right)H\left(M_{\eta}^{B,R}\left(u\right)\right)\cdot F_{s}\left(t-u\right)du}_{c_{2}^{\eta}\left(t\right)}\right]$},
\end{align*}
therefore $\M_{\eta,r}^{B,R}\left(t\right)$ is a martingale because
it is a sum of martingales. 

Similar to \ref{lem:BCP}, we now cam apply Lemma \ref{lem:GB1}
on the term $c_{1}^{\eta}\left(t\right)$ and Lemma \ref{lem:GB2}
on the term $c_{2}^{\eta}\left(t\right)$ to generate a upper bound
for the second moment of each:

\begin{align*}
\E\left[c_{1}^{\eta}\left(t\right)\right] & \leq\frac{2\lambda_{b}t\E\left[b_{i}^{2}\right]}{\eta},\\
\E\left[c_{2}^{\eta}\left(t\right)\right] & \leq\E\left[b_{j}\right]^{2}\frac{\lambda_{b}t}{\eta}.
\end{align*}

By Lemma \ref{lem:ProbConvergence},

\[
\sup_{t\in[0,T]}\left|c_{1}^{\eta}\left(t\right)\right|\underset{P}{\rightarrow}0,
\]

\[
\sup_{t\in[0,T]}\left|c_{2}^{\eta}\left(t\right)\right|\underset{P}{\rightarrow}0.
\]

Then by the triangle inequality and the fact that the sum of terms
converging in probability also converges in probability, it follows
that:
\[
\sup_{t\in[0,T]}\left|\M_{\eta,r}^{B,R}\left(t\right)\right|\underset{P}{\rightarrow}0.
\]

This completes the proof.
\end{proof}

\subsection{Proof of the Results in Section \ref{subsubsec:FLConverge}}

\subsubsection{\label{app:BRFL}Proof of the Theorem \ref{thm:BRFL}. }

\BRFLTheorem*

\begin{proof}
By definition, we want to show that for all $\epsilon>0$, 
\begin{align*}
P\left(\sup_{t\in\left[0,T\right]}\left|M_{\eta}^{B,R}\left(t\right)-m^{B,R}\left(t\right)\right|\leq\epsilon\right) & \to1.
\end{align*}
 Denote the compensator
\begin{align*}
\Pi\left(t\right) & :=M_{0}+\lambda d^{*}t-\lambda\int_{0}^{t}H\left(M_{\eta}^{B,R}\left(u\right)\right)du\\
 & \quad+\left(1-p^{*}\right)\cdot\lambda_{b}\cdot\int_{0}^{t}H\left(M_{\eta}^{B,R}(u)\right)F_{s}\left(t-u\right)du.
\end{align*}

Observe, by the triangle inequality,
\begin{align*}
\sup_{t\in\left[0,T\right]}\left|M_{\eta}^{B,R}\left(t\right)-m^{B,R}\left(t\right)\right| & =\sup_{t\in\left[0,T\right]}\left|M_{\eta}^{B,R}\left(t\right)-\Pi\left(t\right)+\Pi\left(t\right)-m^{B,R}\left(t\right)\right|\\
 & \leq\sup_{t\in\left[0,T\right]}\left|M_{\eta}^{B,R}\left(t\right)-\Pi\left(t\right)\right|+\sup_{t\in\left[0,T\right]}\left|\Pi\left(t\right)-m^{B,R}\left(t\right)\right|.
\end{align*}

Define the event $\mathcal{E}_{\epsilon^{\prime}}^{\eta}:=\left\{ \sup_{t\in\left[0,T\right]}\left|M_{\eta}^{B,R}\left(t\right)-\Pi\left(t\right)\right|\leq\epsilon^{\prime}\right\} $.
If the event $\mathcal{E}_{\epsilon^{\prime}}^{\eta}$ occurs, then

\begin{gather*}
\sup_{t\in\left[0,T\right]}\left|M_{\eta}^{B,R}\left(t\right)-m^{B,R}\left(t\right)\right|\\
\leq\epsilon^{\prime}+\sup_{t\in\left[0,T\right]}\left|\Pi\left(t\right)-m^{B,R}\left(t\right)\right|\\
\text{\ensuremath{\leq}}\epsilon^{\prime}+\sup_{t\in\left[0,T\right]}\left|-\lambda\int_{0}^{t}H\left(M_{\eta}^{B,R}\left(u\right)\right)du+\lambda\int_{0}^{t}H\left(m^{B,R}\left(u\right)\right)du\right|\\
\quad+\left(1-p^{*}\right)\cdot\lambda_{b}\cdot\int_{0}^{t}\left[H\left(M_{\eta}^{B,R}(u)\right)-H\left(m^{B,R}(u)\right)\right]F_{s}\left(t-u\right)du.
\end{gather*}

We use the following notations for simplicity:

\begin{align*}
H_{b}^{M}\left(t\right) & :=H\left(M_{\eta}^{B,R}\left(t\right)\right)\\
H_{b}^{m}\left(t\right) & :=H\left(m^{B,R}(t)\right)
\end{align*}

Then,
{\small
\begin{gather*}
\sup_{t\in\left[0,T\right]}\left|M_{\eta}^{B,R}\left(t\right)-m^{B,R}\left(t\right)\right|\\
\leq\epsilon^{\prime}+\sup_{t\in\left[0,T\right]}\left|\Pi\left(t\right)-m^{B,R}\left(t\right)\right|\\
\text{\ensuremath{\leq}}\epsilon^{\prime}+\sup_{t\in\left[0,T\right]}\left|-\lambda_{b}\int_{0}^{t}H_{b}^{M}\left(u\right)du+\lambda_{b}\int_{0}^{t}H_{b}^{m}\left(u\right)du+\left(1-p^{*}\right)\cdot\lambda_{b}\cdot\int_{0}^{t}\left[H_{b}^{M}\left(u\right)-H_{b}^{m}\left(u\right)\right]F_{s}\left(t-u\right)du\right|\\
=\epsilon^{\prime}+\sup_{t\in\left[0,T\right]}\left|\lambda_{b}\int_{0}^{t}H_{b}^{m}\left(u\right)-H_{b}^{M}\left(u\right)du+\left(1-p^{*}\right)\cdot\lambda_{b}\cdot\int_{0}^{t}\left[H_{b}^{M}\left(u\right)-H_{b}^{m}\left(u\right)\right]F_{s}\left(t-u\right)du\right|\\
=\epsilon^{\prime}+\sup_{t\in\left[0,T\right]}\left|\lambda_{b}\int_{0}^{t}\left[H_{b}^{m}\left(u\right)-H_{b}^{M}\left(u\right)\right]\left[\left(1-p^{*}\right)\cdot F_{s}\left(t-u\right)+1\right]du\right| \\
\leq\epsilon^{\prime}+\sup_{t\in\left[0,T\right]}\left\{ \lambda_{b}\int_{0}^{t}\left|\left[H_{b}^{m}\left(u\right)-H_{b}^{M}\left(u\right)\right]\left[\left(1-p^{*}\right)\cdot F_{s}\left(t-u\right)+1\right]\right|du\right\} \\
\leq\epsilon^{\prime}+\sup_{t\in\left[0,T\right]}\left\{ 2\lambda_{b}\int_{0}^{t}\left|H_{b}^{m}\left(u\right)-H_{b}^{M}\left(u\right)\right|du\right\} 
\end{gather*}
}
Since $H$ is $L$-Lipschitz, we have

\begin{align*}
\sup_{t\in\left[0,T\right]}\left|M_{\eta}^{B,R}\left(t\right)-m^{B,R}\left(t\right)\right| & \leq\epsilon^{\prime}+\lambda_{b}L\sup_{t\in\left[0,T\right]}\int_{0}^{t}\left|M_{\eta}^{B,R}\left(u\right)-m^{B,R}\left(u\right)\right|du\\
 & =\epsilon^{\prime}+\lambda_{b}L\int_{0}^{T}\left|M_{\eta}^{B,R}\left(u\right)-m^{B,R}\left(u\right)\right|du\\
 & \leq\epsilon^{\prime}+\lambda_{b}L\int_{0}^{T}\sup_{t\in\left[0,u\right]}\left|M_{\eta}^{B,R}\left(t\right)-m^{B,R}\left(t\right)\right|du
\end{align*}
Then, by Gronwall's Lemma,
\begin{gather*}
\sup_{t\in\left[0,T\right]}\left|M_{\eta}^{B,R}\left(t\right)-m^{B,R}\left(t\right)\right|\leq\epsilon^{\prime}e^{\lambda LT}.
\end{gather*}
If $\mathcal{E}_{\epsilon^{\prime}}^{\eta}$ happens, then $\sup_{t\in\left[0,T\right]}\left|M^{\eta}\left(t\right)-m\left(t\right)\right|\leq\epsilon^{\prime}e^{\lambda LT}$
must be true. So,
\begin{gather*}
P\left(\sup_{t\in\left[0,T\right]}\left|M_{\eta}^{B,R}\left(t\right)-m^{B,R}\left(t\right)\right|\leq\epsilon^{\prime}e^{\lambda LT}\right)\geq P\left(\mathcal{E}_{\epsilon^{\prime}}^{\eta}\right).
\end{gather*}
For any $\epsilon$, we can choose $\epsilon^{\prime}=\frac{\epsilon}{e^{\lambda LT}}.$
Then,
\begin{align*}
P\left(\sup_{t\in\left[0,T\right]}\left|M_{\eta}^{B,R}\left(t\right)-m^{B,R}\left(t\right)\right|\leq\epsilon\right) & =P\left(\sup_{t\in\left[0,T\right]}\left|M^{\eta}\left(t\right)-m\left(t\right)\right|\leq\epsilon^{\prime}e^{\lambda Lt}\right)\\
 & \geq P\left(\mathcal{E}_{\epsilon^{\prime}}^{\eta}\right)
\end{align*}
Recall that $\mathcal{E}_{\epsilon^{\prime}}^{\eta}:=\left\{ \sup_{t\in\left[0,T\right]}\left|M_{\eta}^{B,R}\left(t\right)-\Pi\left(t\right)\right|\leq\epsilon^{\prime}\right\} .$
We now bound $\sup_{t\in\left[0,T\right]}\left|M_{\eta}^{B,R}\left(t\right)-\Pi\left(t\right)\right|$
in the definition of this event:

{\small
\begin{gather*}
\sup_{t\in\left[0,T\right]}\left|M_{\eta}^{B,R}\left(t\right)-\Pi\left(t\right)\right|\\
\leq\sup_{t\in\left[0,T\right]}\left|\underbrace{\sum_{i=1}^{N_{\eta}^{d}\left(t\right)}d_{i}-\lambda d^{*}\eta t}_{A}\right|+\sup_{t\in\left[0,T\right]}\left|\underbrace{\sum_{j=1}^{N_{\eta}^{d}\left(t\right)}b_{j}\cdot\left\{ M_{\eta}^{B,R}\left(a_{j}^{\eta}-\right)\geq b_{j}\right\} -\lambda\int_{0}^{t}H\left(M_{\eta}^{B,R}\left(u\right)\right)du}_{B}\right| \\
\quad+\sup_{t\in\left[0,T\right]}\left|\underbrace{\sum_{j=1}^{N_{\eta}^{d}\left(t\right)}\left(1-p_{j}\right)\cdot b_{j}\cdot\left\{ M_{\eta}^{B,R}\left(a_{j}^{\eta}-\right)\geq b_{j}\right\} \cdot\left\{ t>a_{j}^{\eta}+s_{j}\right\} -\left(1-p^{*}\right)\cdot\lambda_{b}\cdot\int_{0}^{t}H\left(m^{B,R}(u)\right)F_{s}\left(t-u\right)du}_{C}\right|
\end{gather*}
}{\small\par}

By Lemma \ref{lem:Gcpconverge} , $\sup_{t\in\left[0,T\right]}\left\{ \left|A\right|\right\} \underset{a.s.}{\rightarrow}0$
as $\eta\rightarrow\infty.$ 

By Lemma \ref{lem:BCP} and Lemma \ref{lem:RCP}, as $\eta\rightarrow\infty$
we have

\[
\sup_{t\in\left[0,T\right]}\left|B\right|\underset{P}{\rightarrow}0,
\]

\[
\sup_{t\in\left[0,T\right]}\left|C\right|\underset{P}{\rightarrow}0.
\]
Therefore, by the triangle inequality and the fact that the sum of
terms converging in probability also converges in probability, we
know that as $\eta\to\infty$ and any $\epsilon^{\prime}>0$,
\begin{gather*}
P\left(\mathcal{E}_{\epsilon^{\prime}}^{\eta}\right)\to1.
\end{gather*}
So, it follows that
\begin{gather*}
P\left(\sup_{t\in\left[0,T\right]}\left|M_{\eta}^{B,R}\left(t\right)-m^{B,R}\left(t\right)\right|\leq\epsilon\right)\geq P\left(\mathcal{E}_{\epsilon^{\prime}}^{\eta}\right)\to1.
\end{gather*}
This completes the proof.
\end{proof}
\subsection{Detailed Description of Example\ref{ex:returnsmess}\label{app:exreturnsmess}}
\ReturnsExample*
In this simplified example, we start with one donation of 5 units at the initial time $t_0=0$, denoted by $d_1$, and consider two bail requests,  $b_j$ for $j=1,2$, in this example. Both defendants have the same trial time of $s_j\equiv2$ and same poundage rate of $p_j\equiv0$, meaning the entire bail amount is returned if fully fulfilled ($r_j=(1-p_j)b_j=b_j$). 

\par \textbf{At time $t_0 = 0$}, both the Skorokhod approximation $M^{P^*,R^P}(t)$ and the blocking process $M^{B,R^B}(t)$ receive donations of 5 units. \par \textbf{At time $t_1=1$}, the first defendant requests a bail of $b_1 = 6$. At this moment:\begin{itemize}
    \item The infinite acceptance model accepts it and reaches $M^{\infty,R^{\infty}}(t_1) = -1$ with unlimited credit to borrow.
    \item The Skorokhod approximation transforms this via the Skorokhod map:\[M^{P^*,R^P}(t_1)=\phi[M^{\infty,R^{\infty}}](t_1)=M^{\infty,R^{\infty}}(t_1)-\inf_{s\leq t_1}\{0,M^{\infty,R^{\infty}}(t_1)\}=-1-(-1)=0.\] So it partially fulfills the first defendant by 5 units.
    \item In contrast, the blocking process $M^{B,R^B}(t_1)$ entirely rejects $b_1$ due to insufficient funds, maintaining its value at 5.
\end{itemize} 
\par \textbf{At time $t_3 = 3$}, the second defendant arrives, requesting bail $b_2 = 4$:\begin{itemize}
    \item Again, the infinite acceptance model accepts it and becomes negative: $M^{\infty,R^{\infty}}(t_3) = -5$.
    \item Correspondingly, the Skorokhod approximation yields: \[M^{P^*,R^P}(t_3)=\phi[M^{\infty,R^{\infty}}](t_3)=M^{\infty,R^{\infty}}(t_3)-\inf_{s\leq t_3}\{0,M^{\infty,R^{\infty}}(t_3)\}=-5-(-5)=0.\] So, it gives the second defendant no units.
    \item The blocking process, however, accepts $b_2$ since $M^{B,R^B}(t_2) = 5 > 4 = b_2$. Consequently, $M^{B,R^B}(t_3) = 1$.
\end{itemize}
\textbf{At time $t_4 = 4$}, the first defendant returns the bail amount received, denoted by $r_1$, \begin{itemize}
    \item The Skorokhod approximation receives the full amount of $b_1$, just as it would under the infinite acceptance model. Thus, $r_1=b_1=6$. This result is proven in our main Theorem \ref{thm:partial}. This can be checked by the formula:
     \[M^{P^*,R^P}(t_4)=\phi[M^{\infty,R^{\infty}}](t_4)=M^{\infty,R^{\infty}}(t_4)-\inf_{s\leq t_4}\{0,M^{\infty,R^{\infty}}(s)\}=1-(-5)=6.\]
    \item In contrast,the blocking process remains unchanged at $M^{B,R^B}(t_4) = 1$ as it never accepted $b_1$.
\end{itemize}
\textbf{At time $t_6 = 6$}, the second defendant returns the bail amount received, denoted by $r_2$, \begin{itemize}
    \item By the same logic as at time $t_4=4$, the Skorokhod approximation receives the full returned amount, so $r_2=b_2=4$.Thus, $M^{P^*,R^P}(t_6)=10.$
    \item The blocking process also receives full returned amount $r_2=b_2=4$, so $M^{B,R^B}(t_6)=5.$
\end{itemize}
\textbf{Throughout this timeline:}\begin{itemize}
    \item At $t=t_0$, $M^{P^*,R^P}(t) = M^{B,R^B}(t)$.
    \item Between $t=t_1$ and $t=t_3$,$M^{P^*,R^P}(t) < M^{B,R^B}(t)$.
    \item From $t=t_4$ to $t=t_6$, $M^{P^*,R^P}(t) > M^{B,R^B}(t)$.   
\end{itemize}
\subsection{Proof of the Results in Section \ref{subsecion:OrdNR}}
In this subsection, we will prove the relationships in Equation \eqref{eq:AllOrderNR}. This will be done individually in Propositions \ref{prop:OrdBSkrNR}, \ref{prop:OrdBSkrInfNR} and \ref{prop:OrdSkrPNR}. The full ordering will be shown again in Theorem \ref{thm:orderNR}.

\subsubsection{Supporting Propositions for Theorem \ref{thm:orderNR}.}
\begin{proposition}\label{prop:OrdBSkrNR}
Consider the processes in Equations \ref{eq:M-B} and \ref{eq:M-infty-star}. Then, the
following is true: 
\begin{gather*}
M^{B}\left(t\right)\geq M^{\infty*}\left(t\right).
\end{gather*}
\end{proposition}

\begin{proof}
Recall the definition 
\begin{gather*}
M^{\infty*}\left(t\right)=M^{\infty}\left(t\right)-\inf_{s\leq t}\left(0,M^{\infty}\left(s\right)\right).
\end{gather*}
If $\min_{s\leq t}\left\{ 0,M^{\infty}\left(s\right)\right\} \geq0$,
then, in this case,
\begin{gather*}
M^{B}\left(t\right)=M^{\infty*}\left(t\right)=M^{\infty}\left(t\right).
\end{gather*}
by the same argument in the base case of the proof
of the Theorem \ref{thm:partial}. In short, by Corollary \ref{cor:partial},
if the Skorokhod shift $\inf_{s\leq t}\left\{ 0,M^{\infty}\left(s\right)\right\} =0,$
then both processes $M^{B}\left(\cdot\right)$ and $M^{\infty}\left(\cdot\right)$
at times before $t$ have been larger than all bail requests at the
time of their arrival. Thus, $M^{B}\left(\cdot\right)$ and $M^{\infty}\left(\cdot\right)$
are equal at all times before time $t.$

In the other case where $\min_{s\leq t}\left\{ 0,M^{\infty}\left(s\right)\right\} <0$,
then 
\begin{gather*}
M^{\infty*}\left(t\right)=M^{\infty}\left(t\right)-\inf_{s\leq t}\left\{ 0,M^{\infty}\left(s\right)\right\} .
\end{gather*}
Because $\left[0,t\right]$ is a closed interval, by the
extreme value theorem there exists an $s_{t}^{*}\leq t$ such that
\begin{gather*}
M^{\infty}\left(s_{t}^{*}\right)=\inf_{s\leq t}\left\{ 0,M^{\infty}\left(s\right)\right\} .
\end{gather*}
So, 
\begin{gather*}
M^{\infty*}\left(t\right)=M^{\infty}\left(t\right)-M^{\infty}\left(s_{t}^{*}\right).
\end{gather*}
Observe, 
\begin{gather*}
M^{B}\left(t\right)=M^{B}\left(t\right)-M^{B}\left(s_{t}^{*}\right)+M^{B}\left(s_{t}^{*}\right).
\end{gather*}
Note that $M^{B}\left(s_{t}^{*}\right)\geq0$ because it is the \emph{blocking}
process, so 
\begin{gather*}
M^{B}\left(t\right)\geq M^{B}\left(t\right)-M^{B}\left(s_{t}^{*}\right).
\end{gather*}
Then,
\begin{align*}
M^{B}\left(t\right) & \geq M^{B}\left(t\right)-M^{B}\left(s_{t}^{*}\right)\\
 & =\sum_{i=N^{d}\left(s_{t}^{*}\right)+1}^{N^{d}\left(t\right)}d_{i}-\sum_{j=N^{b}\left(s_{t}^{*}\right)+1}^{N^{b}\left(t\right)}b_{j}\left\{ M^{B}\left(a_{j}-\right)\geq b_{j}\right\} \\\displaybreak[1]
 & \geq\sum_{i=N^{d}\left(s_{t}^{*}\right)+1}^{N^{d}\left(t\right)}d_{i}-\sum_{j=N^{b}\left(s_{t}^{*}\right)+1}^{N^{b}\left(t\right)}b_{j}\\
 & =M^{\infty}\left(t\right)-M^{\infty}\left(s_{t}^{*}\right)\\
 & =M^{\infty*}\left(t\right).
\end{align*}
This completes the proof.
\end{proof}

\begin{proposition}\label{prop:OrdSkrPNR}
Consider the processes in s \ref{eq:M-infty-star} and \ref{eq:M-P}. Then,
the following is true: 

\[
M^{\infty*}\left(t\right)=M^{P}\left(t\right)
\]
\end{proposition}
\begin{proof}
This result directly follows by apply Theorem \ref{thm:partial}.
\end{proof}
\begin{proposition}\label{prop:OrdBSkrInfNR}

Consider the processes in Equations \ref{eq:M-infty-star} and \ref{eq:M-infty} . Then,
the following is true: 
\begin{gather*}
M^{\infty*}\left(t\right)\geq M^{\infty}\left(t\right).
\end{gather*}
\end{proposition}
\begin{proof}
It suffices to show that 
\begin{gather*}
M^{\infty*}\left(t\right)-M^{\infty}\left(t\right)\geq0.
\end{gather*}
The result follows from the following observation: 
\begin{gather*}
M^{\infty*}\left(t\right)-M^{\infty}\left(t\right)=-\min_{s\leq t}\left(0,M^{\infty}\left(s\right)\right)\geq0.
\end{gather*}
\end{proof}

\subsubsection{Proof of  Theorem \ref{thm:orderNR}.\label{app:OrdNR}}
\orderNRTheorem*
\begin{proof}
This follows from combining Propositions  \ref{prop:OrdBSkrNR}, \ref{prop:OrdSkrPNR} and \ref{prop:OrdBSkrInfNR}.  
\end{proof}

\subsection{Proof of the Results in Section \ref{subsection:OrdR}.}
In this subsection, we will prove the relationships in Equation \eqref{eq:AllOrdR}. This will be done individually in Propositions \ref{prop:OrdPInf} and \ref{prop:OrdSkrP}, and the full ordering will be shown again in Theorem \ref{thm:OrdR}.
\subsubsection{Supporting Propositions for Theorem \ref{thm:OrdR}.}

\begin{proposition}\label{prop:OrdSkrP}
Consider the two processes $M^{P^{*},R^{\infty}}\left(t\right)$ and
$M^{P,R^{P}}\left(t\right)$ defined in Equations 
\eqref{eq:SkrR} and \eqref{eq:Partial}. The followings holds: 
\begin{gather*}
M^{P^{*},R^{\infty}}\left(t\right)\geq M^{P,R^{P}}\left(t\right).
\end{gather*}
\end{proposition}
\begin{proof}
By using induction, we wish to show $M^{P^{*},R^{\infty}}\left(t\right)\geq M^{P,R^{P}}\left(t\right)$
for all $t\in\left[a_{k}+,a_{k+1}-\right]$ for all $k\geq0$ where
we use the convention $a_{0}=0$.

For the base case of $k=0$, for $t\in\left[a_{0}+,a_{1}-\right]=\left[0,a_{1}-\right]$,
there are no bail requests yet. So, 
\begin{gather*}
M^{P^{*},R^{\infty}}\left(t\right)=M_{0}+\sum_{i=1}^{N^{d}\left(t\right)}d_{i}=M^{P,R^{P}}\left(t\right).
\end{gather*}

Now we do the induction. Assume that $M^{P^{*},R^{\infty}}\left(t\right)\geq M^{P,R^{P}}\left(t\right)$
for all $t\in\left[a_{l}+,a_{l+1}-\right]$. We will show it is also
true for $t\in\left[a_{l+1}+,a_{l+2}-\right]$. So, by assumption,
\begin{gather*}
M^{P^{*},R^{\infty}}\left(a_{l+1}-\right)\geq M^{P,R^{P}}\left(a_{l+1}-\right).
\end{gather*}
Now, observe that for any $t\in\left[a_{l+1}+,a_{l+2}-\right]$, 
\begin{align*}
M^{P^{*},R^{\infty}}\left(t\right) & =M_{0}+\sum_{i=1}^{N^{d}\left(t\right)}d_{i}+\sum_{i=1}^{\infty}p_{i}b_{i}\left\{ a_{i}+s_{i}\leq t\right\} -\sum_{i=1}^{N^{b}\left(t\right)}b_{i}\wedge M^{*,R^{\infty}}\left(a_{i}-\right)\\
 & =M^{P^{*},R^{\infty}}\left(a_{l+1}+\right)+\underbrace{\sum_{i=N^{d}\left(a_{l+1}\right)}^{N^{d}\left(t\right)}d_{i}+\sum_{i=1}^{\infty}p_{i}b_{i}\left\{ a_{l+1}+<a_{i}+s_{i}\leq t\right\} }_{:=(\star)}\\
M^{P,R^{P}}\left(t\right) & =M_{0}+\sum_{i=1}^{N^{d}\left(t\right)}d_{i}+\sum_{i=1}^{\infty}p_{i}\left(b_{i}\wedge M^{P,R^{P}}\left(a_{i}-\right)\right)\left\{ a_{i}+s_{i}\leq t\right\} -\sum_{i=1}^{N^{b}\left(t\right)}b_{i}\wedge M^{P,R^{P}}\left(a_{i}-\right)\\
 & =M^{P,R^{P}}\left(a_{l+1}+\right)+\sum_{i=N^{d}\left(a_{l+1}+\right)}^{N^{d}\left(t\right)}d_{i}+\sum_{i=1}^{\infty}p_{i}\left(b_{i}\wedge M^{P,R^{P}}\left(a_{i}-\right)\right)\left\{ a_{l+1}+<a_{i}+s_{i}\leq t\right\} \\
 & \leq M^{P,R^{P}}\left(a_{l+1}+\right)+\underbrace{\sum_{i=N^{d}\left(a_{l+1}+\right)}^{N^{d}\left(t\right)}d_{i}+\sum_{i=1}^{\infty}p_{i}b_{i}\left\{ a_{l+1}+<a_{i}+s_{i}\leq t\right\} }_{\left(\star\right)}.
\end{align*}
So, it suffices to prove that $M^{P,R^{P}}\left(a_{l+1}+\right)\leq M^{P^{*},R^{\infty}}\left(a_{l+1}+\right)$.
Observe, 
\begin{align*}
M^{P^{*},R^{\infty}}\left(a_{l+1}+\right) & =M^{P^{*},R^{\infty}}\left(a_{l+1}-\right)-b_{l+1}\wedge M^{P^*,R^{\infty}}\left(a_{l+1}-\right)\\
M^{P,R^{P}}\left(a_{l+1}+\right) & =M^{P,R^{P}}\left(a_{l+1}-\right)-b_{l+1}\wedge M^{P,R^{P}}\left(a_{l+1}-\right).
\end{align*}
Recall by the induction hypothesis $M^{P,R^{P}}\left(a_{l+1}-\right)\leq M^{P^*,R^{\infty}}\left(a_{l+1}-\right)$.
Now, we consider different cases. 
\begin{itemize}
\item If $b_{l+1}\leq M^{P,R^{P}}\left(a_{l+1}-\right)\leq M^{P^*,R^{\infty}}\left(a_{l+1}-\right)$,
then 
\begin{align*}
M^{P^{*},R^{\infty}}\left(a_{l+1}+\right) & =M^{P^{*},R^{\infty}}\left(a_{l+1}-\right)-b_{l+1}\\
M^{P,R^{P}}\left(a_{l+1}+\right) & =M^{P,R^{P}}\left(a_{l+1}-\right)-b_{l+1}.
\end{align*}
So $M^{P,R^{P}}\left(a_{l+1}+\right)\leq M^{P^{*},R^{\infty}}\left(a_{l+1}+\right)$
by the induction hypothesis. 
\item If $M^{P,R^{P}}\left(a_{l+1}-\right)\leq M^{P^*,R^{\infty}}\left(a_{l+1}-\right)\leq b_{l+1}$,
then 
\begin{align*}
M^{P^{*},R^{\infty}}\left(a_{l+1}+\right) & =M^{P^{*},R^{\infty}}\left(a_{l+1}-\right)-M^{P^{*},R^{\infty}}\left(a_{l+1}-\right)=0\\
M^{P,R^{P}}\left(a_{l+1}+\right) & =M^{P,R^{P}}\left(a_{l+1}-\right)-M^{P,R^{P}}\left(a_{l+1}-\right)=0.
\end{align*}
So, $M^{P,R^{P}}\left(a_{l+1}+\right)=M^{P^{*},R^{\infty}}\left(a_{l+1}+\right)$
in this case. 
\item If $M^{P,R^{P}}\left(a_{l+1}-\right)\leq b_{l+1}\leq M^{*,R^{\infty}}\left(a_{l+1}-\right)$,
then 
\begin{align*}
M^{P^{*},R^{\infty}}\left(a_{l+1}+\right) & =M^{P^{*},R^{\infty}}\left(a_{l+1}-\right)-b_{l+1}\geq0\\
M^{P,R^{P}}\left(a_{l+1}+\right) & =M^{P,R^{P}}\left(a_{l+1}-\right)-M^{P,R^{P}}\left(a_{l+1}-\right)=0.
\end{align*}
So $M^{P,R^{P}}\left(a_{l+1}+\right)\leq M^{P^{*},R^{\infty}}\left(a_{l+1}+\right)$
in this case. 
\end{itemize}
So, we have proven that $M^{P,R^{P}}\left(a_{l+1}+\right)\leq M^{P^{*},R^{\infty}}\left(a_{l+1}+\right)$
in all cases. This completes the induction. 
\end{proof}
\begin{proposition}\label{prop:OrdPInf}
Consider the two processes $M^{P^{*},R^{\infty}}\left(t\right)$ and
$M^{P,R^{P}}\left(t\right)$ defined in Equations. The followings
holds: \\
\[
M^{P,R^{P}}\left(t\right)\geq M^{\infty,R^{\infty}}\left(t\right).
\]
\end{proposition}
\begin{proof}
It suffices to show that $0\leq M^{P,R^{P}}\left(t\right)-M^{\infty,R^{\infty}}\left(t\right)$.
We rewrite the inequality as follows:

\begin{gather*}
0\leq\sum_{j=1}^{\infty}p_{j}\left(b_{j}\wedge M^{P,R^{P}}\left(a_{j}-\right)\right)\left\{ a_{j}+s_{j}\leq t\right\} -\sum_{i=1}^{\infty}b_{j}\wedge M^{P,R^{P}}\left(a_{j}-\right)\left\{ a_{j}\leq t\right\} \\
-\sum_{j=1}^{\infty}p_{j}b_{j}\left\{ a_{j}+s_{j}\leq t\right\} +\sum_{j=1}^{\infty}b_{j}\left\{ a_{j}\leq t\right\} ,
\end{gather*}

Next, we rearrange the terms,

\begin{align*}
0 & \leq\underbrace{\sum_{i=1}^{\infty}p_{i}\left(b_{i}\wedge M^{P,R^{P}}\left(a_{i}-\right)\right)\left\{ a_{i}+s_{i}\leq t\right\} -\sum_{i=1}^{\infty}p_{i}b_{i}\left\{ a_{i}+s_{i}\leq t\right\} }_{\text{A}}\\
 & \quad\underbrace{+\sum_{i=1}^{\infty}b_{i}\left\{ a_{i}\leq t\right\} -\sum_{i=1}^{\infty}b_{i}\wedge M^{P,R^{P}}\left(a_{i}-\right)\left\{ a_{i}\leq t\right\} }_{\text{B}}.
\end{align*}
We now analyze each term separately. First consider term A:

\begin{align*}
A & =\sum_{i=1}^{\infty}\left(b_{i}\wedge M^{P,R^{P}}\left(a_{i}-\right)-b_{j}\right)p_{i}\left\{ a_{i}+s_{i}\leq t\right\} \\
 & =\sum_{i=1}^{\infty}-\left(b_{i}-M^{P,R^{P}}\left(a_{i}-\right)\right)^{+}p_{i}\left\{ a_{i}+s_{i}\leq t\right\} .
\end{align*}
Next, we examine term B:

\begin{align*}
B & =\sum_{i=1}^{\infty}\left(b_{i}-b_{i}\wedge M^{P,R^{P}}\left(a_{i}-\right)\right)\left\{ a_{i}\leq t\right\} \\
 & =\sum_{i=1}^{\infty}\left(b_{i}-M^{P,R^{P}}\left(a_{i}-\right)\right)^{+}\left\{ a_{i}\leq t\right\} .
\end{align*}
 Observe from above, $0\leq p_{i}\leq1$ , and $\left\{ a_{i}+s_{i}\leq t\right\} $
, $\left\{ a_{i}\leq t\right\} $ are indicator functions, each taking
values either 0 or 1. Thus, 

\[
\left(b_{i}-M^{P,R^{P}}\left(a_{i}-\right)\right)^{+}\left\{ a_{i}\leq t\right\} -\left(b_{i}-M^{P,R^{P}}\left(a_{i}-\right)\right)^{+}p_{i}\left\{ a_{i}+s_{i}\leq t\right\} \geq0.
\]
Therefore, combining terms A and B gives: 
\begin{align*}
A+B &  =M^{P,R^{P}}\left(t\right)-M^{\infty,R^{\infty}}\left(t\right)\geq0
\end{align*}
\end{proof}
\subsubsection{Proof of Theorem \ref{thm:OrdR}.\label{app:OrdR}}
\OrdRTheorem*
\begin{proof}
This follows from combining Propositions \ref{prop:OrdSkrP} and \ref{prop:OrdPInf}.   
\end{proof}

\end{document}